\documentclass[a4paper,11pt,twoside]{article}
\usepackage{pst-fill,pst-grad}
\usepackage{textcomp}
\usepackage[english]{babel}
\usepackage[utf8x]{inputenc}
\usepackage{graphicx}
\usepackage{multicol}
\usepackage{float}
\usepackage{fancyhdr}
\usepackage[matrix,arrow,curve]{xy}
\usepackage{pstricks} 
\usepackage{amsmath,amsfonts,verbatim,afterpage,theorem,euscript,mathrsfs,amssymb}
\usepackage{array}
\usepackage{dsfont}
\usepackage{hyperref}
\newcommand{\mysection}{\setcounter{equation}{0} \section}

\numberwithin{equation}{section}
\newtheorem{Theorem}{Theorem}[section]
\newtheorem{Lemma}{Lemma}[section]
\newtheorem{Proposition}{Proposition}[section]
\newtheorem{Corollary}{Corollary}[section]
\newtheorem{Remark}{Remark}[section]
\newtheorem{Definition}{Definition}[section]
\newcommand{\vu}{\vec{u}}
\newcommand{\vw}{\vec{\omega}}
\newcommand{\vv}{\vec{v}}
\newcommand{\vf}{\vec{f}}
\newcommand{\vg}{\vec{g}}
\newcommand{\vV}{\vec{V}}
\newcommand{\vA}{\vec{A}}
\newcommand{\vB}{\vec{B}}

\newcommand{\vD}{\vec{D}}
\newcommand{\vE}{\vec{E}}
\newcommand{\vF}{\vec{F}}

\newcommand{\vH}{\vec{H}}
\newcommand{\vphi}{\vec{\phi}}
\newcommand{\vUU}{\vec{\mathcal{U}}}
\def \p{\mathfrak{p}}
\def \q{\mathfrak{q}}
\newcommand{\vn}{\vec{\nabla}}

\def \R{\mathbb{R}^3}
\def \M{\mathcal{M}}
\title{\bf A new approach for the regularity of weak solutions of the 3D Boussinesq system}
\usepackage{authblk}
\author{Diego Chamorro\footnote{\emph{diego.chamorro@univ-evry.fr}} }
\affil{\footnotesize LaMME, Univ. Evry, CNRS\\ Universit\'e Paris-Saclay, 91025, Evry, France.}
\author{Claudiu M\^indril\u{a}\footnote{\emph{mindrila@karlin.mff.cuni.cz}} }
\affil{\footnotesize  Department of Mathematical Analysis\\ Faculty of Mathematics and Physics, Charles University\\ Sokolovsk\'a 83, 18675 Praha 8 - Karl\'in}

\begin{document}
\maketitle
\begin{scriptsize}
\abstract{We address here the problem of regularity for weak solutions of the 3D Boussinesq equation. By introducing the new notion of partial suitable solutions, which imposes some conditions over the velocity field only, we show a local gain of regularity for the two variables $\vu$ and $\theta$.}\\[2mm]
\end{scriptsize}
\mysection{Introduction}
We consider in this article the 3D Boussinesq equation and for $\vu:[0,+\infty[\times \R\longrightarrow \R$ a divergence free velocity field, $p:[0,+\infty[\times \R\longrightarrow \mathbb{R}$ a pressure and $\theta:[0,+\infty[\times \R\longrightarrow \mathbb{R}$ the temperature, we have the system
\begin{equation}\label{Equation_Boussinesq_intro}
\begin{cases}
\partial_t\vu= \Delta \vu -(\vu\cdot \vn)\vu -\vn p+\theta e_3,\quad div(\vu)=0,\\[2mm]
\partial_t\theta=\Delta\theta -\vu\cdot \vn\theta,\\[2mm]
\vu(0,x)=\vu_0(x),\; div(\vu_0)=0, \quad \theta(0,x)=\theta_0(x),
\end{cases}
\end{equation}
where $e_3=(0,0,1)^t$. Note that if $\theta\equiv 0$ we then recover the usual 3D Navier-Stokes system which contains many challenging open problems (see the book \cite{PGLR1} for a detailed and up to date treatment of the Navier-Stokes equations).\\

The Boussinesq system (in 2D or 3D) was extensively studied from many points of view, see \cite{CDiB80}, \cite{Chae06} and \cite{Danchin08} for existence results, \cite{Brandolese12} for the large time behavior,  \cite{Gancedo017}, \cite{Guo96} and \cite{Hu015} for some regularity properties. See also \cite{Hmidi}, \cite{Lazar22} and \cite{Ye17} and the references there in for other problems related to the Boussinesq equations.\\

In this work we will focus on the 3D case which, by its proximity to the Navier-Stokes problem, is slightly more delicate to handle (in terms of regularity issues) than the 2D case. The Boussinesq system is of course different than the Navier-Stokes equations: indeed, the presence of the temperature $\theta$ in the first equation of (\ref{Equation_Boussinesq_intro}) and the coupling via the drift term $\vu \cdot \vn \theta$ in the second equation induce some interesting modifications. For example, applying the divergence operator to the first equation in (\ref{Equation_Boussinesq_intro}) we obtain the following relationship
\begin{equation}\label{Equation_Pressure_intro}
-\Delta p = div \big(div (\vu\otimes \vu)\big)-\partial_{x_3}\theta,
\end{equation}
which gives an equation for the pressure that is distinct than in the case of the Navier-Stokes problem since in order to recover the pressure $p$ we need \emph{both} variable $\vu$ and $\theta$. Another difference is related to the energy inequalities: from some initial data $\vu_0, \theta_0\in L^2(\R)$ and for some fixed time $0<T^*<+\infty$ we can construct weak solutions $(\vu, \theta)\in L^\infty([0,T^*], L^2(\R))\cap L^2([0,T^*], \dot{H}^1(\R))$ that satisfy the following energy inequalities (valid for all $0<t<T^*$)
\begin{eqnarray}
\|\vu(t,\cdot)\|_{L^2}^2+2\int_{0}^t\|\vn\otimes\vu(s, \cdot)\|_{L^2}^2ds&\leq &C(\|\vu_0\|_{L^2}^2+t^2\|\theta_0\|_{L^2}^2)\label{Energy_Inequalities_U}\\
\|\theta(t,\cdot)\|_{L^2}^2+2\int_{0}^t\|\vn\theta(s, \cdot)\|_{L^2}^2ds&\leq &\|\theta_0\|_{L^2}^2,\label{Energy_InequalitiesTheta}
\end{eqnarray}
but here the estimate (\ref{Energy_Inequalities_U}) is \emph{not} uniform in time and this fact can cause some problems when studying global solutions. See for example \cite{Brandolese12}, \cite{CDiB80} and \cite{Danchin08} for more details concerning existence issues.\\

Concerning regularity theories, the celebrated Caffarelli-Kohn-Nirenberg criterion for the Navier-Stokes system (see \cite{CKN}) was extended in \cite{Guo96} to the Boussinesq equation. This theory is based in the notion of \emph{suitable solutions} which are, roughly speaking, weak solutions that satisfy a local energy inequality. Let us note that in the mentioned work \cite{Guo96}, the suitability condition is imposed to the variables $\vu$ and $\theta$.\\

We will show here that this condition over the two variables is redundant and that it is enough to consider some behavior for the velocity field $\vu$ only: this fact will lead us to the notion of \emph{partial suitable solutions} and with the help of this concept we will see how to deduce a gain of regularity for \emph{both} variables. \\

As we are interested in the behavior of the solutions on small neighborhood of points, we will consider  parabolic balls which are defined in the following manner:
\begin{equation}\label{Def_BolaParabolica}
Q_R(t, x)= ]t-R^2, t+R^2[\times B_R(x),
\end{equation}
for some radius $0<R<1$ such that $t-R^2>0$ with $0<t<T^*$, where $T^*$ is a fixed bounded time for which we have the estimates (\ref{Energy_Inequalities_U}) and (\ref{Energy_InequalitiesTheta}):  we thus have $(\vu, \theta)\in L^\infty([0,T^*], L^2(\R))\cap L^2([0,T^*], \dot{H}^1(\R))$ and also $(\vu, \theta)\in L^\infty_t L^2_x(Q_R)\cap L^2_t \dot{H}_x^1(Q_R)$. In this general framework, we have the following definition:
\begin{Definition}[Partial suitable solutions]\label{Partial_suitable_solutions}
Consider  $\vu,\theta \in L_t^{\infty}L_x^2(Q_R)\cap L_t^2\dot H_x^1(Q_R)$ two functions that satisfy the equation (\ref{Equation_Boussinesq_intro}) in the weak sense over the set $Q_R$. Assume moreover that we have the following local information over the pressure: $p\in L^{\frac{3}{2}}_{t,x}(Q_R)$. We will say that $(\vu,p,\theta)$ is a \emph{partial suitable solution} for the 3D Boussinesq equation (\ref{Equation_Boussinesq_intro}) if the distribution $\mu$ given by the expression
\begin{eqnarray}\label{Def_Mu}
\mu=-\partial_t|\vu|^2+\Delta|\vu|^2-2|\vn\otimes\vu|^2-div\left((|\vu|^2+2p)\vu\right)+\theta e_3\cdot\vu,
\end{eqnarray}
is a non-negative locally finite measure on $Q_R$.
\end{Definition}
The condition $p\in L^{\frac{3}{2}}_{t,x}$ is rather classical (see \cite[Chapter 13]{PGLR1}) and the contribution of the variable $\theta$ via the relationship (\ref{Equation_Pressure_intro}) does not cause any interference as we have  $\theta \in L_t^{\infty}L_x^2(Q_R)\cap L_t^2\dot H_x^1(Q_R)$.\\

This notion of \emph{partial} suitable solution is useful to deduce some local energy inequalities for the velocity field $\vu$, but we do not need to impose a similar condition to the variable $\theta$, which can be seen at this stage as a non-divergence free external force. Indeed, by a separate study of each of these variables we will obtain our main result:
\begin{Theorem}\label{Teo_HolderRegularity}
Consider $(\vu, p, \theta)$ a partial suitable solution for the 3D Boussinesq equation (\ref{Equation_Boussinesq_intro}) in the sense of the Definition \ref{Partial_suitable_solutions}.
There exists a small constant  $0<\epsilon^*\ll1$  such that if for some point $(t_0,x_0)\in Q_R$ we have
\begin{equation}\label{Condition_PetitEpsilon}
\limsup_{r\to 0}\frac{1}{r}\int _{]t_0-r^2, t_0+r^2[\times B(x_0,r)}|\vn \otimes \vu|^2dx ds<\epsilon^*,
\end{equation}
then, the solution $(\vu,\theta)$ is H\"older continuous in time and space for some exponent $0<\alpha \ll1$ in a small neighborhood of $(t_0,x_0)$.
\end{Theorem}

Some remarks are in order here. First note that besides the partial suitable condition we only need a mild behavior for the gradient of the velocity field $\vu$ (stated in the hypothesis (\ref{Condition_PetitEpsilon}) above), and thus no particular constraint is asked for the variable $\theta$. Next remark that in this context we can obtain a (local) gain of regularity for both variables $\vu$ and $\theta$: when dealing with regularity issues we can thus observe that the variable $\vu$ \emph{dominates} the variable $\theta$, in the sense that we can deduce some regularity information on the variable $\theta$ from the behavior of the variable $\vu$. This fact seems (to the best of our knowledge) to be new in the study of the regularity properties for the 3D Boussinesq system.

Note also that the gain of regularity is stated in terms of H\"older spaces (in time and space) over a small neighborhood of points $(t_0, x_0)$ where we have (\ref{Condition_PetitEpsilon}), thus it will be quite natural to use parabolic H\"older spaces which are defined in the expression (\ref{Def_Parabolic_Holder}) below. Finally, let us mention that the points $(t_0, x_0)\in [0,T]\times \R$ for which we have (\ref{Condition_PetitEpsilon}) are called \emph{regular points} and following \cite{Guo96} (or adapting the ideas of \cite{Scheffer} or \cite[Section 13.10]{PGLR1}) it can be proven that the parabolic 1-Hausdorff measure of the set of \emph{singular points} is null. \\

The outline of the article is the following. In Section \ref{Secc_HolderMorrey} we recall the definition of the parabolic H\"older spaces as well as the notion of parabolic Morrey spaces. These spaces, although completely absent in the statement of Theorem \ref{Teo_HolderRegularity},  are a powerful tool when studying problems related to regularity in PDEs (see the key Lemma \ref{Lemma_parabolicHolder} below) and in this article we will use them in a systematic manner. In Section \ref{Sec_Propriedades_Vu1} we study the variable $\vu$  considering the variable $\theta$ as an external force and we will obtain a gain of information for $\vu$ in terms of Morrey spaces. Section \ref{Secc_GainTheta} is devoted to the study of the variable $\theta$ and we will see how to obtain a gain of integrability (also stated in terms of Morrey spaces) for this variable. Finally, in Section \ref{Secc_Final}, gathering all the information available on the variables $\vu$ and $\theta$ we will prove Theorem \ref{Teo_HolderRegularity}.
\section{Parabolic H\"older and Morrey spaces}\label{Secc_HolderMorrey}

We will consider the homogeneous space $(\mathbb{R}\times \R, d, \mu)$ where $d$ is the parabolic distance given by 
$d\big((t,x), (s,y)\big)=|t-s|^{\frac{1}{2}}+|x-y|$ and where $\mu$ is the usual Lebesgue measure $d\mu=dxdt$. We then define the homogeneous (parabolic) H\"older spaces $\dot{\mathcal{C}}^\alpha(\mathbb{R}\times \R, \R)$ with $0<\alpha<1$ by the usual condition:
\begin{equation}\label{Def_Parabolic_Holder}
\|\vphi\|_{\dot{\mathcal{C}}^\alpha}=\underset{(t,x)\neq (s,y)}{\sup}\frac{|\vphi(t,x)-\vphi(s,y)|}{\left(|t-s|^{\frac{1}{2}}+|x-y|\right)^\alpha}<+\infty,
\end{equation}
and it is with respect to this functional space that we will obtain the regularity gain announced.\\

Now, for $1< p\leq q<+\infty$, the (parabolic) Morrey spaces $\mathcal{M}_{t,x}^{p,q}(\mathbb{R}\times \R)$ are defined as the set of measurable functions $\vphi:\mathbb{R}\times\R\longrightarrow \R$ that belong to the space $(L^p_{t,x})_{loc}$ such that $\|\vphi\|_{M_{t,x}^{p,q}}<+\infty$ where
\begin{equation}\label{DefMorreyparabolico}
\|\vphi\|_{\mathcal{M}_{t,x}^{p,q}}=\underset{x_{0}\in \R, t_{0}\in \mathbb{R}, r>0}{\sup}\left(\frac{1}{r^{5(1-\frac{p}{q})}}\int_{|t-t_{0}|<r^{2}}\int_{B(x_{0},r)}|\vphi(t,x)|^{p}dxdt\right)^{\frac{1}{p}}.
\end{equation}
Morrey spaces appear to be very convenient functional spaces when dealing with regularity issues as it was pointed out in \cite{Kukavica}, \cite{PGLR1} and \cite{Robinson}.\\

We present now some well-known facts:
\begin{Lemma}[H\"older inequalities]\label{lemma_Product}
\begin{itemize}
\item[]
\item[1)]If $\vf, \vg:\mathbb{R} \times \R\longrightarrow \R$ are two functions such that $\vf\in \M_{t,x}^{p,q} (\mathbb{R} \times \R)$ and $\vg\in L^{\infty}_{t,x} (\mathbb{R} \times \R)$, then for all $1\leq p\leq q<+\infty$ we have
$\|\vf\cdot\vg\|_{\M_{t,x}^{p,q}}\leq  C\|\vf\|_{\M_{t,x}^{p, q}} \|\vg\|_{L^{\infty}_{t,x}}$.
\item[2)] Let $1\leq p_0 \leq q_0 <+\infty$, $1\leq p_1\leq q_1<+\infty$ and $1\leq p_2\leq q_2<+\infty$. If $\tfrac{1}{p_1}+\tfrac{1}{p_2}= \frac{1}{p_0}$ and $\tfrac{1}{q_1}+\tfrac{1}{q_2}=\tfrac{1}{q_0}$, then for two measurable functions $\vf, \vg:\mathbb{R} \times \R\longrightarrow \R$ such that $\vf\in \mathcal{M}^{p_1, q_1}_{t,x}(\mathbb{R} \times \R)$ and $\vg\in \mathcal{M}^{p_2, q_2}_{t,x}(\mathbb{R} \times \R)$, we have the following H\"older inequality in Morrey spaces 
$$\|\vf\cdot \vg\|_{\mathcal{M}^{p_0, q_0}_{t,x}}\leq \|\vf\|_{\mathcal{M}^{p_1, q_1}_{t,x}}\|\vg\|_{\mathcal{M}^{p_2, q_2}_{t,x}}.$$
\end{itemize}
\end{Lemma} 
\begin{Lemma}[Localization]\label{lemma_locindi}
Let $\Omega$ be a bounded set of $\mathbb{R} \times \R$. If we have $1\leq p_0 \leq q_0$, $1\leq p_1\leq q_1$ with the condition $q_0 \leq q_1<+\infty$ and if the function $\vf:\mathbb{R} \times \R\longrightarrow \R$  belongs to the space $\M_{t,x}^{p_1,q_1}(\mathbb{R} \times \R)$ then we have the following localization property 
$$\|\mathds{1}_{\Omega}\vf\|_{\M_{t,x}^{p_0, q_0}} \leq C\|\mathds{1}_{\Omega}\vf\|_{\M_{t,x}^{p_1,q_1}}\leq C\|\vf\|_{\M_{t,x}^{p_1,q_1}}.$$
\end{Lemma} 

In our work, the notion of parabolic Riesz potential (and its properties) will be crucial and for some index $0<\mathfrak{a}<5$ we define the parabolic Riesz potential $\mathcal{L}_{\mathfrak{a}}$ of a locally integrable function $\vec f: \mathbb{R}\times\mathbb{R}^3\to \mathbb{R}^3$ by
\begin{equation}\label{Def_ParabolicRieszPotential}
\mathcal{L}_{\mathfrak{a}}(\vf)(t,x)=\int_{\mathbb{R}}\int_{\mathbb{R}^3}
\frac{1}{(|t-s|^{\frac{1}{2}}+|x-y|)^{5-\mathfrak{a}}}\vec{f}(s,y)dyds.
\end{equation}
Then, we have the following property in Morrey spaces
\begin{Lemma}[Adams-Hedberg inequality]\label{Lemma_Hed}
If $0<\mathfrak{a}<\frac{5}{q}$, $1<p\le q<+\infty$ and $\vf\in \M_{t,x}^{p,q}(\mathbb{R} \times \R)$,
then for $\lambda=1-\frac{\mathfrak{a}q}{5}$ we have the following boundedness property in Morrey spaces:
\begin{equation*}
\|\mathcal{L}_{\mathfrak{a}}(\vf)\|_{\M_{t,x}^{\frac{p}{\lambda},\frac{q}{\lambda}}}\le C\|\vf\|_{\M_{t,x}^{p,q}}. 
\end{equation*}
\end{Lemma}
The three lemmas above constitute our main tools in Morrey spaces. For a more detailed study of these functional spaces we refer to \cite{Adams} and \cite{PGLR1}.

\section{A partial gain of information for the variable $\vu$}\label{Sec_Propriedades_Vu1}
In this section we will only focus our study in the variable $\vu$ and its equation: 
\begin{equation*}
\partial_t \vu  = \Delta \vu  -(\vu \cdot \vn)\vu-\vn p +\theta e_3,\qquad div(\vu)=0. 
\end{equation*}
Here, the variable $\theta$ can be seen as an external force for which we have the information $\theta\in L^\infty([0,T], L^2(\R))\cap L^2([0,T],\dot{H}^1(\R))$, note that by interpolation we also have
\begin{equation}\label{Information_Interpolation_Theta}
\theta \in L^{\frac{10}{3}}_{t,x}([0,T]\times \R).
\end{equation}
In order to obtain a gain of information in the variable $\vu$, we will first consider some estimates for the pressure and then we will deduce inequalities for the velocity field $\vu$.
\subsection{Study of the Pressure}
Since the pressure $p$ satisfies the equation (\ref{Equation_Pressure_intro}) and it depends on the velocity field $\vu$ and on the temperature $\theta$, we first need to establish some controls for it. Although the equation for the pressure is in the case of the 3D Boussinesq equation relatively similar to the equation of the pressure for the 3D Navier-Stokes one, there are some differences that must be treated carefully.\\ 

We introduce the following quantities:  for a point $(t,x)\in \mathbb{R}\times \R$ and for a real parameter $r>0$ we write
\begin{equation}\label{Def_Invariants}
\begin{split}
\mathbb{A}_r(t,x)&=\sup_{t-r^2<s<t+r^2}\frac{1}{r}\int_{B(x,r)}|\vu(s,y)|^2dy, 
\qquad\qquad \alpha_r(t,x)=\frac{1}{r} \int _{Q_r(t,x)}|\vn \otimes \vu(s,y)|^2dyds,\\
\mathbb{B}_r(t,x)&=\frac{1}{r^2} \int _{Q_r(t,x)}|\vu(s,y)|^3dyds, 
\qquad\qquad \qquad\mathcal{P}_r(t,x)=\frac{1}{r^2} \int _{Q_r(t,x)}|p (s,y)|^{\frac{3}{2}}dyds,
\end{split}
\end{equation}
and when the context is clear we will simply write $\mathbb{A}_r=\mathbb{A}_r(t,x)$. Note that the previous quantities correspond to the information $L^\infty_tL^2_x$, $L^2_t\dot{H}^1_x$, $L^3_{t,x}$ and $L^{\frac32}_{t,x}$, in particular we have the following relationships:
\begin{eqnarray}\label{Formules_Definitions}
r\mathbb{A}_r=\|\vu\|_{L^\infty_tL^2_x(Q_r)}^2, \qquad r\alpha_r=\|\vn\otimes\vu\|_{L^2_{t,x}(Q_r)}^2,\qquad r^2 \mathcal{P}_r=\|p\|_{L^{\frac{3}{2}}_{t,x}(Q_r)}^{\frac{3}{2}}.
\end{eqnarray}
With these quantities above we will deduce now a useful estimate for the pressure.
\begin{Proposition}[Pressure estimate]\label{Proposition_EstimatePressure}
Under the hypotheses of Theorem \ref{Teo_HolderRegularity} and with the notations given in (\ref{Def_Invariants}) for any $0<r<\frac{\rho}{2}<R$ we have the inequality
\begin{equation}\label{Estimate_Pressure}
\mathcal{P}_r^{\frac{2}{3}}\le C\biggl(\left(\frac{\rho}{r}\right)\left(\mathbb{A}_\rho
\alpha_\rho\right)^{\frac{1}{2}}+\left(\frac{r}{\rho}\right)^{\frac{3}{2}}\rho^{\frac{5}{6}}\|\theta\|_{L_t^{\infty}L_x^{2} (Q_\rho)}+\left(\frac{r}{\rho}\right)^{\frac{2}{3}} 
\mathcal{P}_\rho^{\frac{2}{3}} \biggr).
\end{equation}
\end{Proposition}
\noindent {\bf Proof.} For proving the inequality above, we will start by the following estimate (where $\sigma$ is a real number such that $0<\sigma<\frac{1}{2}$):
\begin{equation}\label{Estimation_Pression_Intermediaire}
\|p\|_{L_{t,x}^{\frac{3}{2}}(Q_\sigma)} \leq C \left( \sigma^{\frac{1}{3}}\|\vu\|_{L_t^{\infty}L_x^{2} (Q_1)} \|\vn \otimes \vu\|_{L_{t,x}^{2} (Q_1)} +\sigma^{\frac{17}{6}}\|\theta\|_{L_t^{\infty}L_x^{2} (Q_1)}+
\sigma^{2}\|p\|_{L_{t,x}^{\frac{3}{2}} (Q_1)}\right),
\end{equation}
here $Q_\sigma$ and $Q_1$ are parabolic balls of radius $\sigma$ and $1$ respectively (the definition of such balls is given in (\ref{Def_BolaParabolica})).\\ 

To obtain this inequality (\ref{Estimation_Pression_Intermediaire}) we introduce $\eta : \R \longrightarrow [0, 1]$ a smooth function supported in the ball $B_1$ such that $\eta \equiv 1$ on the ball $B_{\frac35}$ and  $\eta \equiv 0$ outside the ball $B_{\frac45}$, note also that over the set $Q_\sigma$ we have the identity $p=\eta p$. Thus, by a straightforward calculation we have
$$ - \Delta (\eta p) = -\eta \Delta p + (\Delta \eta)p - 2  \sum^3_{i= 1}\partial_i ((\partial_i \eta) p),$$
from which we deduce the inequality
\begin{equation*}
\|p\|_{L_{t,x}^{\frac{3}{2}}(Q_\sigma)} =\|\eta p\|_{L_{t,x}^{\frac{3}{2}}(Q_\sigma)} \leq \left\|\frac{\big(-\eta \Delta p\big)}{(-\Delta)}\right\|_{L_{t,x}^{\frac{3}{2}}(Q_\sigma)} + \left\|\frac{(\Delta \eta) p}{(-\Delta)}\right\|_{L_{t,x}^{\frac{3}{2}}(Q_\sigma)} +2\sum^3_{i= 1}\left\|\frac{\partial_i ( (\partial_i \eta) p)}{(-\Delta)}\right\|_{L_{t,x}^{\frac{3}{2}}(Q_\sigma)}.
\end{equation*}
For the first term above, since we have the equation (\ref{Equation_Pressure_intro}) for the pressure we can write
\begin{eqnarray}
\|p\|_{L_{t,x}^{\frac{3}{2}}(Q_\sigma)}&\leq  &C\underbrace{\left\|\frac{1}{(-\Delta)}\Big(\eta \,  \sum^3_{i,j= 1}\partial_i \partial_j (u_i u_j)  \Big)\right\|_{L_{t,x}^{\frac 32}(Q_\sigma)}}_{(p_1)} +\underbrace{\left\|\frac{1}{(-\Delta)}\Big(\eta\partial_{x_3}\theta\Big)\right\|_{L_{t,x}^{\frac 32}(Q_\sigma)}}_{(p_2)}\label{FormuleEtaPression}\\
&&+\underbrace{\left\|\frac{(\Delta \eta) p}{(-\Delta)}\right\|_{L_{t,x}^{\frac{3}{2}}(Q_\sigma)} }_{(p_3)} +2\sum^3_{i= 1}\underbrace{\left\|\frac{\partial_i ( (\partial_i \eta) p)}{(-\Delta)}\right\|_{L_{t,x}^{\frac{3}{2}}(Q_\sigma)}}_{(p_4)}.\notag
\end{eqnarray}
We will study each term above separately.
\begin{itemize}
\item For the term $(p_1)$ in (\ref{FormuleEtaPression}), if we denote by $C_{i,j} = u_i (u_j - (u_j)_1)$ where $ (u_j)_1$ is the average of $u_j$ over the ball of radius $1$, since $\vu$ is divergence free, we have the formula $\displaystyle{\sum^3_{i,j= 1}}\partial_i \partial_j (u_i u_j)=\displaystyle{\sum^3_{i,j= 1}}\partial_i \partial_j C_{i,j}$ and thus we can write
\begin{eqnarray}
(p_1)&=&\left\|\frac{1}{(-\Delta)}\Big(\eta \,  \sum^3_{i,j= 1}\partial_i \partial_j  C_{i,j}  \Big)\right\|_{L_{t,x}^{\frac 32}(Q_\sigma)}\notag\\
&\leq& \sum^3_{i,j= 1}\left\|\frac{1}{(-\Delta)} \left(\partial_i \partial_j(\eta C_{i,j} ) - \partial_i \big((\partial_j \eta) C_{i,j} \big) - \partial_j \big((\partial_i \eta) C_{i,j} \big) + 2(\partial_i \partial_j \eta) C_{i,j}\right)\right\|_{L_{t,x}^{\frac 32}(Q_\sigma)}.\label{FormulePression_I}
\end{eqnarray}
Denoting by $\mathcal{R}_i=\frac{\partial_i}{\sqrt{-\Delta}}$ the usual Riesz transforms on $\mathbb{R}^3$, by the boundedness of these operators in Lebesgue spaces and using the support properties of the auxiliary function $\eta$, we have for the first term above (in the space variable):
\begin{eqnarray*}
\left\|\frac{\partial_i \partial_j}{(-\Delta)} \eta C_{i,j}(t,\cdot) \right\|_{L^{\frac32} (B_\sigma)} &\leq &\|\mathcal{R}_i \mathcal{R}_j (\eta C_{i,j} )(t,\cdot) \|_{L^{\frac32} (\R)}  \leq  C\|\eta C_{i,j}(t,\cdot)\|_{L^{\frac32} (B_1)}\\ 
&\leq &C\|u_i(t,\cdot) \|_{L^{2}(B_1)} \|u_j(t,\cdot) - (u_j)_1\|_{L^{6} (B_1)}\\ &\leq &C\|\vu(t,\cdot)\|_{L^{2} (B_1)} \|\vn \otimes \vu(t,\cdot)\|_{L^{2} (B_1)},
\end{eqnarray*}
where we used H\"older and Poincaré inequalities in the last line. Now taking the $L^{\frac32}$-norm in the time variable of the previous inequality we obtain:
\begin{equation}\label{p11}
\left\|\frac{\partial_i \partial_j}{(-\Delta)} \eta C_{i,j} \right\|_{L_{t,x}^{\frac32} (Q_\sigma)}\leq C\sigma^{\frac{1}{3}}\|\vu\|_{L^{\infty}_tL^{2}_x (Q_1)} \|\vn \otimes \vu\|_{L_{t,x}^{2} (Q_1)}. 
\end{equation}
For the other terms of (\ref{FormulePression_I}) we note that $\partial_i \eta$ vanishes on $B_{\frac35} \cup B^c_{\frac45}$ and since $B_\sigma \subset B_{\frac12}\subset B_{\frac{3}{5}}$, using the integral representation for the operator $\frac{\partial_i}{(- \Delta )}$ we have for the second term of (\ref{FormulePression_I}) the estimate
\begin{eqnarray}
\left\|\frac{\partial_i}{(- \Delta )}\big((\partial_j\eta)C_{i,j}\big)(t,\cdot)\right\|_{L^{\frac32}(B_\sigma)} &\leq &C\sigma^2\left\|\frac{\partial_i}{(- \Delta )}\big((\partial_j\eta)C_{i,j}\big)(t,\cdot)\right\|_{L^{\infty}(B_\sigma)}\notag\\
&\leq &C \, \sigma^2 \left\|\int_{\{\frac35<|y|<\frac45\}} \frac{x_i - y_i }{|x-y|^3} \big((\partial_j \eta) C_{i,j} \big)(t,y)  \, dy\right\|_{L^{\infty}(B_\sigma)},\label{KernelEstimate1}
\end{eqnarray}
since this is a convolution with a bounded kernel (due to the support properties in the variables $x$ and $y$), by the Young inequalities we have
\begin{eqnarray*}
& \leq &C \, \sigma^2 \|C_{i,j}(t,\cdot)\|_{L^{1} (B_1)} \leq C \, \sigma^2 \|u_i(t,\cdot) \|_{L^{2} (B_1)} \|u_j(t,\cdot) - (u_j)_1\|_{L^{2}(B_1)} \\
&\leq& C \,  \|\vu(t,\cdot)\|_{L^{2} (B_1)}\|\vn \otimes \vu(t,\cdot)\|_{L^{2} (B_1)},
\end{eqnarray*}
where we used the same ideas as in (\ref{p11}) before and the fact that $0<\sigma<1$. Thus taking the $L^{\frac32}$-norm in the time variable, we obtain
\begin{equation}\label{p12}
\left\|\frac{\partial_i}{(- \Delta )}\big((\partial_j\eta)C_{i,j}\big)\right\|_{L^{\frac32}_{t,x}(Q_\sigma)} \leq C\sigma^{\frac{1}{3}}\|\vu\|_{L^{\infty}_tL^{2}_x (Q_1)} \|\vn \otimes \vu\|_{L_{t,x}^{2} (Q_1)}. 
\end{equation}
A symmetric argument gives 
\begin{equation}\label{p13}
\left\|\frac{\partial_{j}}{(- \Delta )}\big((\partial_i\eta)C_{i,j}\big)\right\|_{L^{\frac32}_{t,x}(Q_\sigma)} \leq C\sigma^{\frac{1}{3}}\|\vu\|_{L^{\infty}_tL^{2}_x (Q_1)} \|\vn \otimes \vu\|_{L_{t,x}^{2} (Q_1)}, \end{equation}
and observing that the convolution kernel associated to the operator $\frac{1}{(-\Delta)}$ is $\frac{C}{|x|}$, following the same ideas we have for the last term of (\ref{FormulePression_I}) the inequality
\begin{equation}\label{p14}
\left\|\frac{(\partial_i \partial_j \eta) C_{i,j}}{(-\Delta)}\right\|_{L_{t,x}^{\frac 32}(Q_\sigma)}\leq C\sigma^{\frac{1}{3}}\|\vu\|_{L^{\infty}_tL^{2}_x (Q_1)} \|\vn \otimes \vu\|_{L_{t,x}^{2} (Q_1)}.
\end{equation}

Therefore, combining the estimates \eqref{p11}, \eqref{p12},  \eqref{p13} and  \eqref{p14} and getting back to (\ref{FormulePression_I}) we finally have:
\begin{eqnarray}
(p_1)=\left\|\frac{1}{(-\Delta)}\Big(\eta \,  \sum^3_{i,j= 1}\partial_i \partial_j (u_i u_j)  \Big)\right\|_{L_{t,x}^{\frac 32}(Q_\sigma)}\leq C\left(\sigma^{\frac{1}{3}}\|\vu\|_{L^{\infty}_tL^{2}_x (Q_1)} \|\vn \otimes \vu\|_{L_{t,x}^{2} (Q_1)}\right).\label{FormulePression_I1}
\end{eqnarray}
\item For the term $(p_2)$ of (\ref{FormuleEtaPression}) we write
$$(p_2)=\left\|\frac{1}{(-\Delta)}\Big(\eta\partial_{x_3}\theta\Big)\right\|_{L_{t,x}^{\frac 32}(Q_\sigma)}\leq \underbrace{\left\|\frac{1}{(-\Delta)}\partial_{x_3}\big(\eta\theta\big)\right\|_{L_{t,x}^{\frac 32}(Q_\sigma)}}_{(p_{2a})}+\underbrace{\left\|\frac{1}{(-\Delta)}\Big(\big(\partial_{x_3}\eta\big)\theta\Big)\right\|_{L_{t,x}^{\frac 32}(Q_\sigma)}}_{(p_{2b})}.$$
For the quantity $(p_{2a})$ we start by considering the space variable and by the H\"older inequality with $\frac{2}{3}=\frac{1}{2}+\frac{1}{6}$ we write
\begin{eqnarray*}
\left\|\mathds{1}_{B_\sigma}\frac{1}{(-\Delta)}\partial_{x_3}\big(\eta\theta\big)(t,\cdot)\right\|_{L^{\frac 32}(B_\sigma)}&\leq &C\|\mathds{1}_{B_\sigma}\|_{L^2(B_\sigma)}\left\|\frac{1}{(-\Delta)}\partial_{x_3}\big(\eta\theta\big)(t,\cdot)\right\|_{L^{6}(B_\sigma)}\\
&\leq &C\sigma^{\frac{3}{2}}\left\|\frac{\partial_{x_3}}{(-\Delta)}\big(\eta\theta\big)(t,\cdot)\right\|_{L^{6}(\R)},
\end{eqnarray*}
and since the Riesz transforms are bounded in Lebesgue spaces we can write
\begin{eqnarray*}
\left\|\frac{\partial_{x_3}}{(-\Delta)}\big(\eta\theta\big)(t,\cdot)\right\|_{L^{6}(\R)}&\leq &C\sigma^{\frac{3}{2}}\left\|\frac{1}{\sqrt{(-\Delta)}}\big(\eta\theta\big)(t,\cdot)\right\|_{L_{t,x}^{6}(\R)}\\
&\leq& C\sigma^{\frac{3}{2}}\left\|(\eta\theta)(t,\cdot)\right\|_{L^{2}(\R)}\leq C\sigma^{\frac{3}{2}}\|\eta(t,\cdot)\|_{L^\infty}\|\theta(t, \cdot)\|_{L^2(B_1)},
\end{eqnarray*}
where we used the Hardy-Littlewood-Sobolev inequalities and the localizing properties of the function $\eta$. Now, we integrate with respect to the time variable, and with the previous controls we obtain
\begin{equation}\label{Estimate_P2A}
(p_{2a})=\left\|\frac{1}{(-\Delta)}\partial_{x_3}\big(\eta\theta\big)(t,\cdot)\right\|_{L_{t,x}^{\frac 32}(Q_\sigma)}\leq C\sigma^{\frac{3}{2}}\sigma^{\frac{4}{3}}\|\theta\|_{L^\infty_tL^2_x(Q_1)}= C\sigma^{\frac{17}{6}}\|\theta\|_{L^\infty_tL^2_x(Q_1)}.
\end{equation}
For the term $(p_{2b})$ we proceed as in (\ref{KernelEstimate1}) and due to the properties of the localizing function $\eta$ we write 
\begin{eqnarray*}
\left\|\frac{1}{(-\Delta)}\Big(\big(\partial_{x_3}\eta\big)\theta\Big)(t,\cdot)\right\|_{L^{\frac 32}(B_\sigma)}&\leq &C\sigma^2\left\|\frac{1}{(-\Delta)}\Big(\big(\partial_{x_3}\eta\big)\theta\Big)(t,\cdot)\right\|_{L^{\infty}(B_\sigma)}\\
&\leq &C\sigma^2\left\|\int_{\{\frac35<|y|<\frac45\}} \frac{1}{|x-y|} \Big(\big(\partial_{x_3}\eta\big)\theta\Big)(t,\cdot)\right\|_{L^{\infty}(B_\sigma)},
\end{eqnarray*}
since the kernel of convolution above is bounded (as $x\in B_\sigma$ with $0<\sigma<\frac12$ and $y\in \{\frac35<|y|<\frac45\}$), we obtain
$$\left\|\frac{1}{(-\Delta)}\Big(\big(\partial_{x_3}\eta\big)\theta\Big)(t,\cdot)\right\|_{L^{\frac 32}(B_\sigma)}\leq C\sigma^2\|\theta(t,\cdot)\|_{L^1(B_1)}\leq C\sigma^2\|\theta(t,\cdot)\|_{L^2(B_1)}.$$
Thus, taking the $L^{\frac32}$-norm in the time variable we have
\begin{equation}\label{Estimate_P2B}
(p_{2b})=\left\|\frac{1}{(-\Delta)}\Big(\big(\partial_{x_3}\eta\big)\theta\Big)\right\|_{L_{t,x}^{\frac 32}(Q_\sigma)}\leq C\sigma^{\frac{10}{3}}\|\theta\|_{L^\infty_tL^2_x(Q_1)}\leq C\sigma^{\frac{17}{6}}\|\theta\|_{L^\infty_tL^2_x(Q_1)},
\end{equation}
(since $0<\sigma<1$ and $\sigma^{\frac{10}{3}}\leq \sigma^{\frac{17}{6}}$). Now with the estimate (\ref{Estimate_P2A}) for $(p_{2a})$ and (\ref{Estimate_P2B}) for $(p_{2b})$ we can write
\begin{equation}\label{FormulePression_II}
(p_2)\leq  C\sigma^{\frac{17}{6}}\|\theta\|_{L^\infty_yL^2_x(Q_1)}.
\end{equation}
\item We continue our study of expression (\ref{FormuleEtaPression}) and for the  term $(p_3)$ we first treat the space variable. Recalling the support properties of the auxiliary function $\eta$ and properties of the convolution kernel associated to the operator $\frac{1}{(-\Delta)}$, we can write as before (see (\ref{KernelEstimate1})):
$$\left\|\frac{(\Delta \eta) p(t,\cdot)}{(-\Delta)}\right\|_{L^{\frac{3}{2}}(B_\sigma)}\leq C\sigma^{2} \| p(t,\cdot)\|_{L^1(B_1)}\leq C\sigma^{2} \| p(t,\cdot)\|_{L^{\frac{3}{2}}(B_1)},$$
and thus, taking the $L^{\frac{3}{2}}$-norm in the time variable we obtain:
\begin{equation}\label{FormulePression_III}
(p_3)=\left\|\frac{(\Delta \eta) p}{(-\Delta)}\right\|_{L^{\frac{3}{2}}_{t,x}(Q_\sigma)}\leq C\sigma^{2} \|p\|_{L^{\frac{3}{2}}_{t,x}(Q_1)}.
\end{equation}
\item For the last term of expression (\ref{FormuleEtaPression}), following the same ideas displayed in (\ref{KernelEstimate1}) we can write
$$\left\|\frac{\partial_i}{(-\Delta)}  (\partial_i \eta) p(t,\cdot) \right\|_{L^{\frac{3}{2}}(B_\sigma)}\leq C\sigma^{2} \| p(t,\cdot)\|_{L^1(B_1)}\leq C\sigma^{2} \|p(t,\cdot)\|_{L^{\frac{3}{2}}(B_1)},$$
and we obtain 
\begin{equation}\label{FormulePression_IIII}
(p_4)=\left\|\frac{\partial_i ( (\partial_i \eta) p )}{(-\Delta)}\right\|_{L_{t,x}^{\frac{3}{2}}(Q_\sigma)}\leq C\sigma^{2} \|p\|_{L^{\frac{3}{2}}_{t,x}(Q_1)}.
\end{equation}
\end{itemize}
Now, gathering the estimates (\ref{FormulePression_I1}), (\ref{FormulePression_II}),  (\ref{FormulePression_III}) and (\ref{FormulePression_IIII}) we obtain the inequality 
$$\|p\|_{L_{t,x}^{\frac{3}{2}}  (Q_\sigma)} \leq C \left( \sigma^{\frac{1}{3}}\|\vu\|_{L_t^{\infty}L_x^{2} (Q_1)} \|\vn \otimes \vu\|_{L_{t,x}^{2} (Q_1)} +\sigma^{\frac{17}{6}}\|\theta\|_{L_t^{\infty}L_x^{2} (Q_1)}+\sigma^{2}\|p\|_{L_{t,x}^{\frac{3}{2}} (Q_1)}\right),$$
which is (\ref{Estimation_Pression_Intermediaire}).\\

To deduce inequality \eqref{Estimate_Pressure}, if we fix $\sigma = \frac{r}{\rho} \leq \frac12$ and if we introduce the functions $p_\rho(t,x)=p(\rho^2t, \rho x)$, $\theta_\rho(t,x)=\theta(\rho^2t, \rho x)$ and $\vu_\rho(t,x)=\vu(\rho^2t, \rho x)$, then from the control (\ref{Estimation_Pression_Intermediaire}) we have
$$\|p_\rho\|_{L_{t,x}^{\frac{3}{2}}  (Q_{\frac{r}{\rho}})} \leq C \left( \left(\frac{r}{\rho}\right)^{\frac{1}{3}}\|\vu_\rho\|_{L_t^{\infty}L_x^{2} (Q_1)} \|\vn \otimes \vu_\rho\|_{L_{t,x}^{2} (Q_1)} + \left(\frac{r}{\rho}\right)^{\frac{17}{6}}\|\theta_\rho\|_{L_t^{\infty}L_x^{2} (Q_1)}+\left(\frac{r}{\rho}\right)^{2}\|p_\rho\|_{L_{t,x}^{\frac{3}{2}} (Q_1)}\right),$$
and by a change of variable we can write
\begin{eqnarray*}
\|p\|_{L_{t,x}^{\frac{3}{2}}(Q_{r})}\rho^{-\frac{10}{3}} &\leq &C \bigg( \left(\frac{r}{\rho}\right)^{\frac{1}{3}}\rho^{-\frac32}\|\vu\|_{L_t^{\infty}L_x^{2}(Q_\rho)} \rho^{-\frac32}\|\vn \otimes \vu\|_{L_{t,x}^{2} (Q_\rho)}+\left(\frac{r}{\rho}\right)^{\frac{17}{6}}\rho^{-\frac{3}{2}}\|\theta\|_{L_t^{\infty}L_x^{2} (Q_\rho)}\\
& & +\left(\frac{r}{\rho}\right)^{2}\rho^{-\frac{10}{3}}\|p\|_{L_{t,x}^{\frac{3}{2}} (Q_\rho)}\bigg).
\end{eqnarray*}
Now, recalling that by (\ref{Def_Invariants}) and (\ref{Formules_Definitions}) we have the identities 
$$r^{\frac{4}{3}}\mathcal{P}_r^{\frac{2}{3}}=\|p\|_{L_{t,x}^{\frac{3}{2}}  (Q_r)}, \quad \rho^{\frac{1}{2}}\mathbb{A}_\rho^{\frac{1}{2}}=\|\vu\|_{L_t^{\infty}L_x^{2} (Q_\rho)} \quad \mbox{and}\quad \rho^{\frac{1}{2}}\alpha_\rho^{\frac{1}{2}}=\|\vn \otimes \vu\|_{L_{t,x}^{2} (Q_\rho)},$$
we finally obtain
$$\mathcal{P}_r^{\frac{2}{3}}\leq C\left(\left(\frac{\rho}{r}\right)(\mathbb{A}_\rho\alpha_\rho)^{\frac{1}{2}}+\left(\frac{r}{\rho}\right)^{\frac{3}{2}}\rho^{\frac{5}{6}}\|\theta\|_{L_t^{\infty}L_x^{2} (Q_\rho)}+\left(\frac{r}{\rho}\right)^{\frac{2}{3}}\mathcal{P}_\rho^{\frac{2}{3}}\right),$$
and this finishes the proof of Proposition \ref{Proposition_EstimatePressure}.\hfill $\blacksquare$
\subsection{Study of the velocity field}
We continue our study with the treatment of the velocity field $\vu$ and we start with a relationship between some of the quantities defined in (\ref{Def_Invariants}) above:
\begin{Lemma}
For $0<r<1$ we have the relationship between $\mathbb{B}_r$, $\mathbb{A}_r$ and $\alpha_r$:
\begin{equation}\label{Estimation_normeLambda3}
\mathbb{B}_r^{\frac{1}{3}}\le C(\mathbb{A}_r+\alpha_r)^{\frac{1}{2}}.
\end{equation}
\end{Lemma}
{\bf Proof.} Indeed, using the definition of $ \mathbb{B}_r$ given in (\ref{Def_Invariants}) above and by H\"older inequality we have
$$\mathbb{B}_r^{\frac{1}{3}}=\frac{1}{r^{\frac{2}{3}}}\|\vu\|_{L_{t,x}^{3}(Q_r)}\le \frac{C}{r^\frac{1}{2}}\|\vu\|_{L_{t,x}^{\frac{10}{3}}(Q_r)}.$$
By interpolation we have  $\|\vu\|_{L_{t,x}^{\frac{10}{3}} (Q_r)} \leq \|\vu\|_{L_t^{\infty}L_x^{2} (Q_r)}^{\frac25}\|\vu\|_{L_t^{2}L_x^{6} (Q_r)}^{\frac{3}{5}}$. Now, for the $L_t^2L_x^6$ norm of $\vu$, we use the classical Gagliardo-Nirenberg inequality (see \cite{Brezis}) to obtain $\|\vu\|_{L_t^{2}L_x^{6} (Q_r)} \leq C\big(\|\vn \otimes\vu\|_{L_t^{2}L_x^{2} (Q_r)} +\|\vu\|_{L_t^{\infty}L_x^{2} (Q_r)}\big)$ and using Young's inequalities we have
\begin{eqnarray*}
\|\vu\|_{L_{t,x}^{\frac{10}{3}} (Q_r)} &\leq& C \|\vu\|_{L_t^{\infty}L_x^{2} (Q_r)}^{\frac 25}\big(\|\vn \otimes \vu\|_{L_t^{2}L_x^{2} (Q_r)}^{\frac 35}+\|\vu\|_{L_t^{\infty}L_x^{2} (Q_r)}^{\frac35} \big) \leq C\big(\|\vu\|_{L_t^{\infty}L_x^{2} (Q_r)}+\|\vn\otimes \vu\|_{L_t^{2}L_x^{2} (Q_r)}\big).
\end{eqnarray*}
Recalling that $\|\vu\|_{L_t^{\infty}L_x^{2} (Q_r)}=r^{\frac12}\mathbb{A}_r^{\frac12}$ and $\|\vn\otimes \vu\|_{L^2_tL_x^{2} (Q_r)}=r^{\frac12}\alpha_r^{\frac12}$, we finally obtain (\ref{Estimation_normeLambda3}).\hfill $\blacksquare$\\

Now, we will establish an estimate that relies in the energy estimate (\ref{Energy_Inequalities_U}):
\begin{Proposition}\label{Proposition_EstimateVelocity}
Under the hypotheses of Theorem \ref{Teo_HolderRegularity} and with the notations given in (\ref{Def_Invariants}) we have for any radius $0<r<\frac{\rho}{2}<1$ the inequality
\begin{equation}\label{Estimate_Velocity}
\mathbb{A}_r+\alpha_r \le C\frac{r^2}{\rho ^2} \mathbb{A}_\rho
+\frac{\rho ^2}{r^2}\alpha_\rho^{\frac{1}{2}}\mathbb{A}_\rho+C\frac{\rho ^2}{r^2}
\mathcal{P}_\rho^{\frac{2}{3}}(\mathbb{A}_\rho+\alpha_\rho)^{\frac{1}{2}}+C\frac{\rho^{\frac{7}{2}}}{r} 
\|\theta\|_{L^\infty_tL_{x}^2(Q_\rho)}\alpha_{\rho}^{\frac{1}{2}}.
\end{equation}
\end{Proposition}
\textbf{Proof.} Following the idea of Scheffer \cite{Scheffer}, we will consider the following test function:
\begin{Lemma}\label{Lemma_testfonc}
Let $0<r<\frac{\rho}{2}<R<1$. Let $\phi \in \mathcal{C}_0^{\infty}
(\mathbb{R}\times \R)$  be a function such that
\begin{equation*}
\phi(s,y)=r^2\omega\left(\frac{s-t}{\rho^2},\frac{y-x}{\rho}\right)
\theta\left(\frac{s-t}{r^2}\right)\mathfrak{g}_{(4r^2+t-s)}(x-y),
\end{equation*}
where $\omega \in \mathcal{C}_0^{\infty}(\mathbb{R}\times \R)$ is positive function whose support is in $Q_1(0,0)$ and equal to 1 in $Q_{\frac{1}{2}}(0,0)$. In addition $\theta$ is a smooth function non negative such that $\theta=1$ over $]-\infty,1[$ and $\theta=0$ over $]2,+\infty[$ and $\mathfrak{g}_t(\cdot)$ is the usual heat kernel. Then, we have the following points.
\begin{itemize}
\item[1)]the function $\phi$ is a bounded non-negative function, and its support is contained in the parabolic ball $Q_\rho$, and for all $(s,y)\in Q_r(t,x)$ we have the lower bound $\phi(s,y)\ge \frac{C}{r}$,
\item[2)] for all $(s,y)\in Q_\rho (t,x)$ with $0<s<t+r^2$ we have $\phi(s,y)\le \frac{C}{r}$,
\item[3)] for all $(s,y)\in Q_\rho(t,x) $ with $0<s<t+r^2$ we have $\vn \phi(s,y)\le \frac{C}{r^2}$,
\item [4)] moreover, for all $(s,y)\in Q_\rho(t,x) $ with $0<s<t+r^2$ we have $|(\partial_s+\Delta)\phi(s,y)|\le C\frac{r^2}{\rho^5}$.
\end{itemize}
\end{Lemma}
A detailed proof of this lemma can be found in \cite[Chapter 13]{PGLR1}.\\ 

With the all the properties of this function $\phi$, exploiting the fact that $(\vu, \theta)$ is a partial suitable solution (it satisfies a local energy inequality) and using the notations (\ref{Def_Invariants}) above, we can write:
\begin{eqnarray}
\mathbb{A}_r+\alpha_r &\leq &\underbrace{\int_{\R}(\partial_t\phi +\Delta \phi)|\vu|^2dxds}_{(1)}+2\underbrace{\int_{\mathbb{R}}\int_{\R}  p (\vu\cdot \vn \phi )dxds}_{(2)}+\underbrace{\int_{\mathbb{R}}\int _{\R} |\vu|^2(\vu \cdot \vn)\phi dx ds}_{(3)}\notag\\
&&+2\underbrace{\int_{\mathbb{R}}\int_{\R}\theta e_3\cdot (\phi\vu)dxds}_{(4)}. \label{Ineq_Energie2}
\end{eqnarray}
The terms of the right hand side above will be studied separately. Indeed, 
\begin{itemize}
\item[$\bullet$] For the quantity (1) in (\ref{Ineq_Energie2}), using the properties of the function $\phi$ given in Lemma \ref{Lemma_testfonc} and by the definition of the quantity $\mathbb{A}_\rho$ given in (\ref{Def_Invariants}) we have
$$\int_{\R}(\partial_t\phi +\Delta \phi)|\vu|^2dxds\leq C \frac{r^2}{\rho^5}\int_{Q_\rho}|\vu|^2dxds=C\frac{r^2}{\rho^5}\int_{t-\rho^2}^{t+\rho^2}\int_{B_\rho}|\vu|^2dxds\le C\frac{r^2}{\rho ^2} \mathbb{A}_\rho.$$
\item[$\bullet$] For the term (2) in (\ref{Ineq_Energie2}), by the properties of the function $\phi$ given in Lemma \ref{Lemma_testfonc} and by the H\"older inequality, we obtain
$$\int_{\mathbb{R}}\int_{\R}  p (\vu\cdot \vn \phi )dxds\leq \frac{C}{r^2}\int_{t-\rho^2}^{t+\rho^2}\int_{B_\rho} | p | |\vu|dxds\le \frac{C}{r^2}\|p\|_{L_{t,x}^{\frac{3}{2}}(Q_\rho)}\|\vu\|_{L_{t,x}^{3}(Q_\rho)},$$
noting that by (\ref{Def_Invariants}) we have $\|p\|_{L_{t,x}^{\frac{3}{2}}(Q_\rho)}=\rho^{\frac{4}{3}} \mathcal{P}_\rho^{\frac{2}{3}}$ and $\|\vu\|_{L_{t,x}^{3}(Q_\rho)}=\rho^{\frac{2}{3}}\lambda_\rho^{\frac{1}{3}}$, we can thus write
$$\int_{\mathbb{R}}\int_{\R}  p (\vu\cdot \vn \phi )dxds\leq  \frac{C}{r^2}\left(\rho^{\frac{4}{3}} \mathcal{P}_\rho^{\frac{2}{3}}\right)\left(\rho^{\frac{2}{3}}\lambda_\rho^{\frac{1}{3}}\right)\le C\frac{\rho ^2}{r^2} \mathcal{P}_\rho^{\frac{2}{3}}( \mathbb{A}_\rho+\alpha_\rho)^{\frac{1}{2}},$$
where in the last estimate we used the control (\ref{Estimation_normeLambda3}).
\item For the term (3) in (\ref{Ineq_Energie2}), let us first define the average $\displaystyle{(|\vu|^2)_\rho=\frac{1}{|B(x,\rho)|}\int_{B(x,\rho)}|\vu(t,y)|^2dy}$ and since $\vu$  is divergence free we have $\displaystyle{\int _{B_\rho} (|\vu|^2)_\rho(\vu \cdot \vn)\phi dx =0}$. Then, we can write by the properties of the function $\phi$ given in Lemma \ref{Lemma_testfonc} and by the H\"older inequality:
\begin{eqnarray*}
\int_{\mathbb{R}}\int_{\R} |\vu|^2(\vu \cdot \vn)\phi dxds&=&\int_{Q_\rho} [|\vu|^2-(|\vu|^2)_\rho](\vu\cdot 
\vn)\phi dx ds\leq \frac{C}{r^2}\int_{t-\rho^2}^{t+\rho^2}\int_{B_\rho} \left||\vu|^2
-(|\vu|^2)_\rho|\right|\vu|dx ds\\ 
&\le &\frac{C}{r^2}\int_{t-\rho^2}^{t+\rho^2}\||\vu|^2-(|\vu|^2)_\rho\|
_{L^{\frac{3}{2}}(B_\rho)}\|\vu(s,\cdot)\|_{L^{3}(B_\rho)}ds.
\end{eqnarray*}
Now, Poincare's inequality implies
\begin{align*}
&\le \frac{C}{r^2}\int_{t-\rho^2}^{t+\rho^2}\|\vn (|\vu(s,\cdot)|^2)\|_{L^{1}(B_\rho)}\|\vu(s,\cdot)\|_{L^{3}(B_\rho)}ds\\
&\le \frac{C}{r^2}\int_{t-\rho^2}^{t+\rho^2} \|\vu(s,\cdot)\|_{L^2(B_\rho)}\|\vn\otimes\vu(s,\cdot)\|_{L^2(B_\rho)}\|\vu(s,\cdot)\|_{L^{3}(B_\rho)}ds  \\
&\le \frac{C}{r^2} \|\vu\|_{L_t^6L_x^2(Q_\rho)}\|\vn \otimes \vu\|_{L_{t,x}
^2(Q_\rho)} \|\vu\|_{L_{t,x}^3(Q_\rho)},
\end{align*}
where in the last inequality we used the H\"older inequality in the time variable. We observe now that by the notations given in (\ref{Def_Invariants}) we can write 
$$\|\vu\|_{L_t^6L_x^2(Q_\rho)}\le C \rho ^{\frac{1}{3}}\|\vu\|_{L_t^{\infty}L_x^2(Q_\rho)}\le C\rho ^{\frac{5}{6}}\mathbb{A}_\rho^{\frac{1}{2}}, \quad \|\vn \otimes \vu\|_{L_{t,x}
^2(Q_\rho)}=\rho^{\frac{1}{2}}\alpha_\rho^{\frac{1}{2}}, \quad  \| \vu\|_{L_{t,x}
^3(Q_\rho)}=\rho^{\frac{2}{3}}\lambda_\rho^{\frac{1}{3}},$$
and we obtain, by (\ref{Estimation_normeLambda3}): 
$$\int_{\mathbb{R}}\int_{\R} |\vu|^2(\vu \cdot \vn)\phi dxds\leq C \frac{\rho^2}{r^2}\mathbb{A}_\rho^{\frac{1}{2}}\alpha_\rho^{\frac{1}{2}}\lambda_\rho^{\frac{1}{3}}\leq  C \frac{\rho^2}{r^2}\mathbb{A}_\rho^{\frac{1}{2}}\alpha_\rho^{\frac{1}{2}}(\mathbb{A}_\rho+\alpha_\rho)^{\frac{1}{2}}\leq C\frac{\rho^2}{r^2} \alpha_\rho^{\frac{1}{2}}
(\mathbb{A}_\rho+\alpha_\rho).$$
\item Finally, for the term (4) in (\ref{Ineq_Energie2}), by the H\"older inequality and by the properties of the function $\phi$ given in Lemma \ref{Lemma_testfonc}, we write
\begin{eqnarray*}
\int_{\mathbb{R}}\int_{\R}(\theta e_3)\cdot(\phi \vu)dxds&\le &\int_{t-\rho^2}^{t+\rho^2}\|\phi(s,\cdot)\|_{L^3_x(B_\rho)}\|\theta(s,\cdot)\|_{L_{x}^2(B_\rho)} \|\vu(s,\cdot)\|_{L_{x}^6(B_\rho)}ds\\
&\leq &C\frac{\rho}{r}\int_{t-\rho^2}^{t+\rho^2}\|\theta(s,\cdot)\|_{L_{x}^2(B_\rho)}\|\vu(s,\cdot)\|_{\dot{H}^1_x(B_\rho)}ds\\
&\leq &C\frac{\rho}{r}\|\theta\|_{L_{t,x}^2(Q_\rho)}\|\vu\|_{L_{t}^2 \dot{H}^1_x(Q_\rho)}\leq C\frac{\rho^3}{r}\|\theta\|_{L_{t}^\infty L_{x}^2(Q_\rho)}\|\vu\|_{L_{t}^2 \dot{H}^1_x(Q_\rho)},
\end{eqnarray*}
where we applied the Sobolev inequalities (see Corollary 9.14 of \cite{LibroBrezis}) and the Cauchy-Schwartz inequality in the time variable.
Since by (\ref{Def_Invariants}) we have $\|\vu\|_{L_{t}^2 \dot{H}^1_x(Q_\rho)}=\rho^{\frac{1}{2}} \alpha_\rho^\frac{1}{2}$, we conclude 
$$\int_{\mathbb{R}}\int_{\R}(\theta e_3)\cdot(\phi \vu)dxds\le C\frac{\rho^{\frac{7}{2}}}{r} 
\|\theta\|_{L_{t}^\infty L_{x}^2(Q_\rho)}\alpha_{\rho}^{\frac{1}{2}}.$$
\end{itemize}
Gathering all these estimates we obtain the inequality (\ref{Estimate_Velocity}) and this ends the proof of Proposition \ref{Proposition_EstimateVelocity}. \hfill$\blacksquare$
\subsection{Iterative process}
With the estimates  (\ref{Estimate_Pressure}) and (\ref{Estimate_Velocity}) and given in Propositions \ref{Proposition_EstimatePressure} and \ref{Proposition_EstimateVelocity}, respectively, we will set up a general inequality that will help us to deduce a gain of integrability. For this, we introduce the notations
\begin{equation}\label{Def_QuantiteIteration}
\mathscr{A}_r=\frac{1}{r^{2(1-\frac{5}{\tau_0})}}\left(\mathbb{A}_r+\alpha_r\right),\quad \mathscr{P}_r=\frac{1}{r^{\frac{3}{2}(1-\frac{5}{\tau_0})}}\mathcal{P}_r\quad \mbox{and}\quad 
\mathscr{O}_r=\mathscr{A}_r+\left(\left(\frac{r}{\rho}\right)^{\frac{15}{\tau_0}-\frac{15}{4}}\mathscr{P}_r\right)^{\frac{4}{3}},
\end{equation}
for a fixed $\tau_0$ such that $\frac{5}{1-\alpha}<\tau_0<6$, which is possible since $0<\alpha<\frac{1}{10}$. We have the following result:
\begin{Lemma}\label{Lema_iteration}
Under the hypotheses of Theorem \ref{Teo_HolderRegularity}, for $0<r<\frac{\rho}{2}<R<1$ there exists a constant $\epsilon>0$ such that
\begin{equation}\label{thetaite1}
\mathscr{O}_{r}(t_0,x_0)\le \frac{1}{2}\mathscr{O}_{\rho}(t_0,x_0)+\epsilon,
\end{equation}
where the point $(t_0,x_0)\in Q_R$ is given by the hypothesis (\ref{Condition_PetitEpsilon}).
\end{Lemma}
\textbf{Proof.} First, by the estimate (\ref{Estimate_Velocity}) we can write
\begin{eqnarray}
\mathscr{A}_r&=&\frac{1}{r^{2(1-\frac{5}{\tau_0})}}\left(\mathbb{A}_r+\alpha_r\right)\notag\\
&\leq& \frac{C}{r^{2(1-\frac{5}{\tau_0})}}\Big(\frac{r^2}{\rho ^2} \mathbb{A}_\rho+\frac{\rho^2}{r^2} \alpha_\rho^{\frac{1}{2}}\mathbb{A}_\rho+\frac{\rho ^2}{r^2}\mathcal{P}_\rho^{\frac{2}{3}}(\mathbb{A}_\rho+\alpha_\rho)^{\frac{1}{2}}+\frac{\rho^{\frac{7}{2}}}{r}\|\theta\|_{L^\infty_tL_{x}^2(Q_\rho)}\alpha_{\rho}^{\frac{1}{2}}\Big),\label{EstimationPourIteration}
\end{eqnarray}
and we will treat each one of the previous terms separately. Indeed, 
\begin{itemize}
\item For the first term of \eqref{EstimationPourIteration} we have
$$\frac{1}{r^{2(1-\frac{5}{\tau_0})}}\left(\frac{r^2}{\rho ^2}\mathbb{A}_\rho\right) \le \frac{1}{r^{2(1-\frac{5}{\tau_0})}}\frac{r^2}{\rho ^2}\rho^{2(1-\frac{5}{\tau_0})}\mathscr{A}_\rho =\left(\frac{r}{\rho}\right)^{\frac{10}{\tau_0}}\mathscr{A}_\rho.$$
\item For the second term of (\ref{EstimationPourIteration}), using the definition of $\mathbb{A}_\rho$ given in (\ref{Def_QuantiteIteration}), we obtain
$$\frac{1}{r^{2(1-\frac{5}{\tau_0})}}\left(\frac{\rho^2}{r^2} \alpha_\rho^{\frac{1}{2}}\mathbb{A}_\rho\right)\leq \frac{1}{r^{2(1-\frac{5}{\tau_0})}}\left(\frac{\rho^2}{r^2} \alpha_\rho^{\frac{1}{2}}\rho^{2(1-\frac{5}{\tau_0})}\mathscr{A}_\rho \right)=\left(\frac{\rho}{r}\right)^{4-\frac{10}{\tau_0}}\mathscr{A}_\rho \alpha_\rho^{\frac{1}{2}}.$$
\item The third term of (\ref{EstimationPourIteration}) follows essentially the same arguments as above and by the definition of the quantities  $\mathbb{A}_\rho$  and  $\mathbb{P}_\rho$ given in (\ref{Def_QuantiteIteration}) we can write
$$\frac{1}{r^{2(1-\frac{5}{\tau_0})}}\left(\frac{\rho ^2}{r^2}\mathcal{P}_\rho^{\frac{2}{3}}( \mathbb{A}_\rho+\alpha_\rho)^{\frac{1}{2}}\right)\leq \left(\frac{\rho}{r}\right)^{4-\frac{10}{\tau_0}}\mathscr{P}_\rho^{\frac{2}{3}}\mathscr{A}_\rho^{\frac{1}{2}}.$$
\item Finally, for the last term of (\ref{EstimationPourIteration}), we have
$$\frac{1}{r^{2(1-\frac{5}{\tau_0})}}\left(\frac{\rho^{\frac{7}{2}}}{r} 
\|\theta\|_{L_{t,x}^2(Q_\rho)}\alpha_{\rho}^{\frac{1}{2}}\right)\leq \left(\frac{\rho}{r}\right)^{3-\frac{10}{\tau_0}} \rho^{\frac12+\frac{10}{\tau_0}}\|\theta\|_{L^\infty_{t}L^2_{x}(Q_\rho)}\alpha_{\rho}^{\frac{1}{2}}.$$
\end{itemize}
Thus, gathering all these estimates, we have 
\begin{equation}\label{estimateA}
\mathscr{A}_r\leq C\Bigg(\left(\frac{r}{\rho}\right)^{\frac{10}{\tau_0}}\mathscr{A}_\rho+\left(\frac{\rho}{r}\right)^{4-\frac{10}{\tau_0}}\mathscr{A}_\rho \alpha_\rho^{\frac{1}{2}}+\left(\frac{\rho}{r}\right)^{4-\frac{10}{\tau_0}}\mathscr{P}_\rho^{\frac{2}{3}}\mathscr{A}_\rho^{\frac{1}{2}}+ \left(\frac{\rho}{r}\right)^{3-\frac{10}{\tau_0}} \rho^{\frac12+\frac{10}{\tau_0}}\|\theta\|_{L^\infty_{t}L^2_{x}(Q_\rho)}\alpha_{\rho}^{\frac{1}{2}}\Bigg).
\end{equation}
Now, for the pressure, from the inequality \eqref{Estimate_Pressure} we can write
\begin{equation}\label{Estimation_RecurrencePression1}
\mathscr{P}_r=\frac{1}{r^{\frac{3}{2}(1-\frac{5}{\tau_0})}}\mathcal{P}_r\le \frac{C}{r^{\frac{3}{2}(1-\frac{5}{\tau_0})}} \biggl(\left(\frac{\rho}{r}\right)^{\frac32}\mathbb{A}_\rho^{\frac{3}{4}}\alpha_\rho^{\frac{3}{4}}+\left(\frac{r}{\rho}\right)^{\frac{3}{2}}\rho^{\frac{5}{6}}\|\theta\|^{\frac{3}{2}}_{L_t^{\infty}L_x^{2} (Q_\rho)}+\left(\frac{r}{\rho}\right)\mathcal{P}_\rho \biggr).
\end{equation}
\begin{itemize}
\item Using the definition of $\mathscr{A}_\rho$ given in (\ref{Def_QuantiteIteration}) we obtain for the first term of the right-hand side above:
$$\frac{1}{r^{\frac{3}{2}(1-\frac{5}{\tau_0})}}\left(\frac{\rho}{r}\right)^{\frac{3}{2}}\mathbb{A}_\rho^{\frac{3}{4}} \alpha_\rho^{\frac{3}{4}}\leq\frac{1}{r^{\frac{3}{2}(1-\frac{5}{\tau_0})}}\left(\frac{\rho}{r}\right)^{\frac{3}{2}}\rho^{\frac{3}{2}(1-\frac{5}{\tau_0})}(\mathscr{A}_\rho\alpha_\rho)^{\frac{3}{4}}=\left(\frac{\rho}{r}\right)^{3-\frac{15}{2\tau_0}}(\mathscr{A}_\rho\alpha_\rho)^{\frac{3}{4}}.$$
\item Now, for the second term of (\ref{Estimation_RecurrencePression1}) we write
$$\frac{C}{r^{\frac{3}{2}(1-\frac{5}{\tau_0})}}\left(\frac{r}{\rho}\right)^{\frac{3}{2}}\rho^{\frac{5}{6}}\|\theta\|^{\frac{3}{2}}_{L_t^{\infty}L_x^{2} (Q_\rho)}\leq C\left(\frac{r}{\rho}\right)^{\frac{15}{2\tau_0}}\rho^{\frac{15}{2\tau_0}-\frac{2}{3}}\|\theta\|^{\frac{3}{2}}_{L^\infty_tL^2_x(Q_\rho)}.$$
\item Finally we use the fact that 
$$\frac{1}{r^{\frac{3}{2}(1-\frac{5}{\tau_0})}} \left(\frac{r}{\rho}\right)\mathcal{P}_\rho= \left(\frac{\rho}{r}\right)^{\frac{1}{2}-\frac{15}{2\tau_0}}\mathscr{P}_\rho$$ (by the definition of $\mathscr{P}_\rho$ given in (\ref{Def_QuantiteIteration})). 
\end{itemize}
We thus conclude that
\begin{align}\label{estimateP}
\mathscr{P}_r\le C\biggl(\left(\frac{\rho}{r}\right)^{3-\frac{15}{2\tau_0}}
(\mathscr{A}_\rho \alpha_\rho)^{\frac{3}{4}}+C\left(\frac{r}{\rho}\right)^{\frac{15}{2\tau_0}}\rho^{\frac{15}{2\tau_0}-\frac{2}{3}}\|\theta\|^{\frac{3}{2}}_{L^\infty_tL^2_x(Q_\rho)}+\left(\frac{\rho}{r}\right)^{\frac{1}{2}-\frac{15}{2\tau_0}}\mathscr{P}_\rho\biggr).
\end{align}
With the estimates (\ref{estimateA}) and (\ref{estimateP}) at hand, we will now introduce a relationship between the parameters $r$ and $\rho$: indeed, let us fix $0<\kappa\ll \frac{1}{2}$ a real number and consider $r=\kappa \rho$, then, by the definition of the quantity $\mathscr{O}_r$ given in (\ref{Def_QuantiteIteration}) we obtain:
\begin{eqnarray}
\mathscr{O}_r&=&\mathscr{A}_r+\left(\kappa^{\frac{15}{\tau_0}-\frac{15}{4}}\mathscr{P}_r\right)^{\frac{4}{3}}\leq C\Bigg(\underbrace{\kappa^{\frac{10}{\tau_0}}\mathscr{A}_\rho+\kappa^{\frac{10}{\tau_0}-4}\mathscr{A}_\rho \alpha_\rho^{\frac{1}{2}}}_{(1)}+\underbrace{\kappa^{\frac{10}{\tau_0}-4}\mathscr{P}_\rho^{\frac{2}{3}}\mathscr{A}_\rho^{\frac{1}{2}}}_{(2)}+\underbrace{\kappa^{\frac{10}{\tau_0}-3} \rho^{\frac{1}{2}+\frac{10}{\tau_0}}\|\theta\|_{L^\infty_{t}L^2_{x}(Q_\rho)}\alpha_{\rho}^{\frac{1}{2}}}_{(3)}\Bigg)\notag\\
&&+C\underbrace{\biggl(\kappa^{\frac{45}{2\tau_0}-\frac{27}{4}}(\mathscr{A}_\rho \alpha_\rho)^{\frac{3}{4}}+\kappa^{\frac{45}{2\tau_0}-\frac{15}{4}}\rho^{\frac{15}{2\tau_0}-\frac{2}{3}}\|\theta\|^{\frac{3}{2}}_{L^\infty_tL^2_x(Q_\rho)}+\kappa^{\frac{45}{2\tau_0}-\frac{17}{4}}\mathscr{P}_\rho\biggr)^{\frac43}}_{(4)}.\label{EstimationTheta}
\end{eqnarray}
We will rewrite now each one of the previous terms:
\begin{itemize}
\item Since by (\ref{Def_QuantiteIteration}) we have $\mathscr{A}_\rho\leq \mathscr{O}_\rho$, it is then easy to see that the term (1) above can be controlled in the following manner:
$$\kappa^{\frac{10}{\tau_0}}\mathscr{A}_\rho+\kappa^{\frac{10}{\tau_0}-4}\mathscr{A}_\rho \alpha_\rho^{\frac{1}{2}}\leq \big(\kappa^{\frac{10}{\tau_0}}+\kappa^{\frac{10}{\tau_0}-4} \alpha_\rho^{\frac{1}{2}}\big)\mathscr{O}_\rho.$$
\item For the quantity (2) in (\ref{EstimationTheta}), using Young's inequality and the relationships given in (\ref{Def_QuantiteIteration}), we observe that
\begin{eqnarray*}
\kappa^{\frac{10}{\tau_{0}}-4}  \mathscr{P}_\rho^{\frac{2}{3}} \mathscr{A}_\rho^{\frac{1}{2}}&=&\kappa^{\frac{10}{\tau_{0}}-4} \left(\kappa^{5(\frac{1}{\tau_{0}} -\frac{1}{2})} \mathscr{P}_\rho^{\frac{2}{3}} \times\kappa^{5(\frac{1}{2} - \frac{1}{\tau_{0}} )} \mathscr{A}_\rho^{\frac{1}{2}}\right)\leq\kappa^{\frac{10}{\tau_{0}}-4} \left( \kappa^{10(\frac{1}{2} - \frac{1}{\tau_{0}} )} \mathscr{A}_\rho + \kappa^{10(\frac{1}{\tau_{0}} - \frac{1}{2} )} \mathscr{P}_\rho^{\frac{4}{3}} \right)\notag \\
&\leq & \kappa \left(\mathscr{A}_\rho+ \left(\kappa^{\frac{15}{\tau_{0}} -\frac{15}{2}}\mathscr{P}_\rho\right)^{\frac{4}{3}} \right) \leq \kappa \, \mathscr{O}_\rho.
\end{eqnarray*}
\item For the term $(3)$ of (\ref{EstimationTheta}),  we just remark that the power of $\kappa$ is $\frac{10}{\tau_0}-3$ which is a negative number since $\frac{5}{1-\alpha}<\tau_0<6$.
\item For the last term of (\ref{EstimationTheta}), since $\left(\kappa^{\frac{15}{\tau_0}-\frac{15}{4}}\mathscr{P}_\rho\right)^{\frac{4}{3}}\leq \mathscr{O}_\rho$ and $\mathscr{A}_\rho\leq \mathscr{O}_\rho$, we have
\begin{eqnarray*}
\biggl(\kappa^{\frac{45}{2\tau_0}-\frac{27}{4}}(\mathscr{A}_\rho \alpha_\rho)^{\frac{3}{4}}+\kappa^{\frac{45}{2\tau_0}-\frac{15}{4}}\rho^{\frac{15}{2\tau_0}-\frac{2}{3}}\|\theta\|^{\frac{3}{2}}_{L^\infty_tL^2_x(Q_\rho)}+\kappa^{\frac{45}{2\tau_0}-\frac{17}{4}}\mathscr{P}_\rho\biggr)^{\frac43}\\
\leq C\biggl(\kappa^{\frac{30}{\tau_0}-9}\mathscr{A}_\rho \alpha_\rho +\kappa^{\frac{30}{\tau_0}-5}\rho^{\frac{10}{\tau_0}-\frac{8}{9}}\|\theta\|^{2}_{L^\infty_tL^2_x(Q_\rho)}+(\kappa^{\frac{45}{2\tau_0}-\frac{17}{4}}\mathscr{P}_\rho)^{\frac43}\biggr)\\
\leq C\bigg(\kappa^{\frac{30}{\tau_0}-9} \alpha_\rho+\kappa^{\frac{10}{\tau_0}-\frac{2}{3}}\bigg)\mathscr{O}_\rho+C\kappa^{\frac{30}{\tau_0}-5}\rho^{\frac{10}{\tau_0}-\frac{8}{9}}\|\theta\|^{2}_{L^\infty_tL^2_x(Q_\rho)}.
\end{eqnarray*}
\end{itemize}
Gathering these estimates we finally obtain
\begin{eqnarray}\label{thetaite2}
\mathscr{O}_r &\le& C\Biggl(\kappa^{\frac{10}{\tau_0}}+\kappa^{\frac{10}{\tau_0}-4}\alpha_\rho^{\frac{1}{2}}+\kappa+\kappa^{\frac{30}{\tau_0}-9}\alpha_\rho+\kappa^{\frac{10}{\tau_0}-\frac{2}{3}}\Biggr)\mathscr{O}_\rho\label{thetaite2}\\
&&+C\kappa^{\frac{10}{\tau_0}-3} \rho^{\frac{1}{2}+\frac{10}{\tau_0}}\|\theta\|_{L^\infty_tL_{x}^2(Q_\rho)}\alpha_{\rho}^{\frac{1}{2}}+C\kappa^{\frac{30}{\tau_0}-5}\rho^{\frac{10}{\tau_0}-\frac{8}{9}}\|\theta\|^{2}_{L^\infty_tL^2_x(Q_\rho)}.\label{thetaite3}
\end{eqnarray}
Futhermore, we claim that we have
\begin{equation}\label{LimsupHypo1}
C\Biggl(\kappa^{\frac{10}{\tau_0}}+\kappa^{\frac{10}{\tau_0}-4}\alpha_\rho^{\frac{1}{2}}+\kappa+\kappa^{\frac{30}{\tau_0}-9}\alpha_\rho+\kappa^{\frac{10}{\tau_0}-\frac{2}{3}} \Biggr)\le\frac{1}{2}.
\end{equation}
Indeed, since $\kappa=\frac{r}{\rho}\ll\frac{1}{2}$ is a fixed small parameter and since $\frac{10}{\tau_0}-\frac{2}{3}>0$ (recall again that $\frac{5}{1-\alpha}<\tau_0<6$), then the quantities $\kappa^{\frac{10}{\tau_0}}$, $\kappa$ and $\kappa^{\frac{10}{\tau_0}-\frac{2}{3}}$ in the previous formula are small. Now, using the fact that we have the control $\alpha_\rho \le \epsilon^*$ which is given in the hypothesis (\ref{Condition_PetitEpsilon}) where $\epsilon^*>0$ is small enough, then the terms $\kappa^{\frac{10}{\tau_0}-4}\alpha_\rho^{\frac{1}{2}}$ and $\kappa^{\frac{30}{\tau_0}-9}\alpha_\rho$   can be made small enough and thus we obtain the estimate \eqref{LimsupHypo1}.\\

To continue, we need to treat the two remaining terms given in (\ref{thetaite3}). For the first one we note that the quantity $\|\theta\|_{L^\infty_{t}L^2_{x}(Q_\rho)}$ is bounded since $\theta\in L_t^{\infty}L_x^2(Q_R)\cap L_t^2\dot H_x^1(Q_R)$, we can apply the same ideas used previously (\emph{i.e.} the fact that $\alpha_\rho \le \epsilon^*\ll1$) to obtain
$$C\kappa^{\frac{10}{\tau_0}-3}\rho^{\frac{1}{2}+\frac{10}{\tau_0}}\|\theta\|_{L^\infty_{t}L^2_{x}(Q_\rho)} \alpha_{\rho}^{\frac{1}{2}}<\frac{\epsilon}{2}.$$
For the last term of (\ref{thetaite3}), recalling that $\frac{5}{1-\alpha}<\tau<6$, we have $\frac{30}{\tau_0}-5>0$ and $\frac{10}{\tau_0}-\frac{8}{9}>0$. Thus, if $0<\rho\ll1$ and since $\kappa \ll\frac{1}{2}$ the quantity $\kappa^{\frac{30}{\tau_0}-5}\rho^{\frac{10}{\tau_0}-\frac{8}{9}}$ can be made small enough to absorb the term $\|\theta\|^{2}_{L^\infty_tL^2_x(Q_\rho)}$ and we obtain
$$\kappa^{\frac{30}{\tau_0}-5}\rho^{\frac{10}{\tau_0}-\frac{8}{9}}\|\theta\|^{2}_{L^\infty_tL^2_x(Q_\rho)}<\frac{\epsilon}{2}.$$

Then, with these estimates at hand and coming back to (\ref{thetaite2}) we conclude that $\mathscr{O}_r\le \frac{1}{2}\mathscr{O}_\rho+\epsilon$ and Lemma \ref{Lema_iteration} is proven. \hfill $\blacksquare$\\
\begin{Proposition}\label{Propo_FirstMorreySpace}
Under the hypotheses of Theorem \ref{Teo_HolderRegularity} consider $(\vu,p,\theta)$ a partial suitable solution for the Boussinesq equations (\ref{Equation_Boussinesq_intro}) over the set $Q_{R}$ given in (\ref{Def_BolaParabolica}). Then there exists a radius $0<R_1<\frac{R}{2}$ and an index $\tau_0$ such that $\frac{5}{1-\alpha}<\tau_0<6$ such that we have the following local Morrey information:
\begin{equation}\label{Conclusion_FirstMorreySpace}
\mathds{1}_{Q_{R_1}(t_0,x_0)} \vu\in \mathcal{M}_{t,x}^{3,\tau_0}(\mathbb{R}\times \R),
\end{equation}
where the point $(t_0,x_0)\in Q_R$ is given by the hypothesis (\ref{Condition_PetitEpsilon}).
\end{Proposition}
{\bf Proof of the Proposition \ref{Propo_FirstMorreySpace}.} 
Lemma \ref{Lema_iteration} paved the way to obtain the wished Morrey information for the velocity $\vu$. Indeed, from the definition of Morrey spaces given in (\ref{DefMorreyparabolico}) we only need to prove that for all radius $r>0$ such that $r<R_1\le \frac{R}{2}$ and  $(t,x)\in Q_{R_1}(t_0,x_0)$, we have
\begin{equation}\label{morrey1}
\int_{Q_r(t,x)} |\vu|^3 dy ds \le C r^{5(1-\frac{3}{\tau_0})},
\end{equation}
and this will imply that $\mathds{1}_{Q_{R_1}}\vu\in \M^{3,\tau_0}(\mathbb{R}\times \R)$. In order to obtain the control (\ref{morrey1}), by the definitions given in (\ref{Def_Invariants}) and by the estimate (\ref{Estimation_normeLambda3}), we observe that
\begin{equation*}
\int_{Q_r(t,x)} |\vu|^3dy ds=r^2\mathbb{B}_r(t,x) \le  r^2(\mathbb{A}_r(t,x)+\alpha_r(t,x))^{\frac{3}{2}}.
\end{equation*}
Hence, it is then enough to prove for all $0<r< {R_1}<\frac{R}{2}<R<1$ and  $(t,x)\in Q_{R_1}$ that one has the control
\begin{equation*}
\mathbb{A}_r(t,x)+\alpha_r(t,x) \le C r^{2(1-\frac{5}{\tau_0})}.
\end{equation*}
Recalling the definition of the quantity $\mathscr{A}_r$ given in (\ref{Def_QuantiteIteration}), we easily  see that the condition (\ref{morrey1}) above is equivalent to prove that there exists some $R_1$ and  $0<\kappa\ll \frac12$ such that for all $n\in \mathbb{N}$ and $(t,x)\in Q_{R_1}(t_0, x_0)$, we have estimates:
\begin{equation}\label{iterative}
\mathscr{A}_{\kappa^n R_1}(t,x)\le C.
\end{equation}
Note that, for any radius $r$ such that $0<r<R_1<\min\{\frac{R}{2},dist(\partial Q_R, (t_0,x_0))\}$ (and since we have $Q_{R_1}(t_0,x_0)\subset Q_{R}$) by the hypotheses of the Theorem \ref{Teo_HolderRegularity}, we have the bounds 
$$\|\vu\|_{L_t^\infty L_x^2(Q_r(t_0,x_0))}\le \|\vu\|_{L_t^\infty L_x^2(Q_R)}<+\infty,\quad \|\vn \otimes\vu\|_{L_{t,x}^2(Q_r(t_0,x_0))}\le \|\vn\otimes\vu\|_{L_{t,x}^2(Q_R)}<+\infty,$$
and $ \|p\|_{L_{t,x}^{\frac{3}{2}}(Q_r(t_0,x_0))}\le \|p\|_{L_{t,x}^{\frac{3}{2}}(Q_R)}<+\infty$. Then, by the notations introduced in (\ref{Def_Invariants}), we have the uniform bounds $\underset{0<r<R}{\sup}\biggl\{ r \mathbb{A}_r,r \alpha_r, r^2\mathcal{P}_r\biggr\}<+\infty$ from which we can deduce by the definition of the quantities $\mathscr{A}_\rho(t_0,x_0)$ and $\mathscr{P}_\rho (t_0,x_0)$ given in (\ref{Def_QuantiteIteration}), the uniform bounds
\begin{equation*}
\sup_{0<r<R}r^{3-\frac{10}{\tau_0}}\mathscr{A}_r(t_0,x_0)<+\infty,\quad\text{and} \quad \sup_{0<r<R}r^{5-\frac{3}{2}(1+\frac{5}{\tau_0})}\mathscr{P}_r (t_0,x_0)<+\infty.
\end{equation*}
Thus, there exists a radius $0<r_0< R$ small such that, by the estimates above, the quantities $\mathscr{A}_{r_0}$ and $\mathscr{P}_{r_0}$ are bounded: indeed, recall that we have $\tau_0>\frac{5}{1-\alpha}>5$ (where $0<\alpha< \frac{1}{10}$) and this implies that all the powers of $r$ in the expression above are positive. As a consequence of this fact, by (\ref{Def_QuantiteIteration}) the quantity $\mathscr{O}_{r_0}$ is itself bounded. Remark also that, if $r_0$ is small enough, then the inequality (\ref{thetaite1}) holds true and we can write $\mathscr{O}_{\kappa r_0}(t_0,x_0)\le \frac{1}{2}\mathscr{O}_{r_0}(t_0,x_0)+\epsilon$. We can iterate this process and we obtain for all $n>1$,
$$\mathscr{O}_{\kappa^nr_0}(t_0,x_0)\le \frac{1}{2^n}\mathscr{O}_{r_0}(t_0,x_0)+\epsilon\sum_{j=0}^{n-1}2^{-j},$$
and therefore there exists $N\ge 1$ such that for all $n\ge N$ we have
$\mathscr{O}_{\kappa^nr_0}(t_0,x_0)\le 4\epsilon$ from which we obtain (using the definition of $\mathscr{O}_r$ given in (\ref{Def_QuantiteIteration})) that
$$\mathscr{A}_{\kappa^Nr_0}(t_0,x_0)\le \frac{1}{8}C \quad \text{and}\quad
\mathscr{P}_{\kappa^Nr_0}(t_0,x_0)\le \frac{1}{32}C.$$
This information is centered at the point $(t_0,x_0)$, in order to treat
the uncentered bound, we can let $\frac{1}{2}\kappa^Nr_0$ to be the radius $R_1$
we want to find, thus for all points $(t,x)\in Q_{R_1}(t_0,x_0)	$ we have that
$Q_{R_1}\subset Q_{2R_1}(t_0,x_0)$, which implies 
\begin{equation*}
\mathscr{A}_{R_1}(t,x)\le 2^{3-\frac{10}{\tau_0}} \mathscr{A}_{2R_1}(t_0,x_0)\le 8 \mathscr{A}_{2R_1}(t_0,x_0)\le 8 \mathscr{A}_{\kappa^N \rho }(t_0,x_0)<C,
\end{equation*}
and $\mathscr{P}_{R_1}(t,x)\le 2^{5-\frac{3}{2}(1+\frac{5}{\tau_0})}\mathscr{P}_{2R_1}(t_0,x_0)\le 32\mathscr{P}_{2R_1}(t_0,x_0)\le 8 \mathscr{P}_{\kappa^N r }(t_0,x_0)<C$. Having obtained these bounds, by the definition of $\mathscr{O}_{R_1}$, we thus get $\mathscr{O}_{R_1}(t,x)\le C$. Applying the Lemma \ref{Lema_iteration} and iterating once more, we find that the same will be true for $\kappa R_1$ and then, for all $\kappa^n R_1$, $n\in \mathbb{N}$. Since by definition we have $\mathscr{A}_{\kappa^n R_1}(t,x)\le  \mathscr{O}_{\kappa^n R_1}(t,x)$ we have finally obtained the estimate $\mathscr{A}_{\kappa^n R_1}(t,x)\le C$ and the inequality (\ref{iterative}) is proven which implies the Proposition \ref{Propo_FirstMorreySpace}.\hfill$\blacksquare$
\begin{Corollary}\label{corolarioMorrey}
Under the hypotheses of Proposition \ref{Propo_FirstMorreySpace}, we also have the following local control:
\begin{equation}\label{ConlusioncorolarioMorrey}
\mathds{1}_{Q_{R_1}(t_0,x_0)}\vn\otimes \vu\in \mathcal{M}^{2,\tau_1}_{t,x}(\mathbb{R}\times \R),\quad \mbox{with}\quad  \frac{1}{\tau_1}=\frac{1}{\tau_0}+\frac{1}{5}.
\end{equation}
\end{Corollary}
\textbf{Proof.} In the previous results we have proved the estimate (\ref{iterative}). Let us recall now that, by the definition of the quantity $\mathscr{A}_r$ given in (\ref{Def_QuantiteIteration}), we can easily deduce for all $0<r\le R_1$ and $(t,x)\in Q_{R_1}$ the control $\alpha_r\leq Cr^{2(1-\frac{5}{\tau_0})}$ which can we rewritten as
$$\frac{1}{r}\bigg(\int _{Q_r(t,x)}|\vn \otimes \vu|^2dyds\bigg)\leq Cr^{2(1-\frac{5}{\tau_0})}.$$
Thus, since $\frac{1}{\tau_1}=\frac{1}{\tau_0}+\frac{1}{5}$, for all $0<r\le R_1$ and $(t,x)\in Q_{R_1}(t_0,x_0)$, we have the estimate
$$\int _{Q_r}|\vn \otimes \vu|^2dyds\leq C r^{3-\frac{10}{\tau_0}}= C r^{5(1-\frac{2}{\tau_1})},$$
and by the definition of Morrey spaces given in (\ref{DefMorreyparabolico}), we obtain that $\mathds{1}_{Q_{R_1}(t_0,x_0)}\vn \otimes \vu\in \mathcal{M}^{2,\tau_1}_{t,x}(\mathbb{R}\times \R)$.\hfill$\blacksquare$
\subsection{A partial gain of information for the variable $\vu$}\label{Sec_Propriedades_Vu2}
\begin{Proposition}\label{Propo_GainenMorreySigma}
Under the hypotheses of Theorem \ref{Teo_HolderRegularity} and within the framework of Proposition \ref{Propo_FirstMorreySpace}, there exists a radius $R_2$ with $0<R_2<R_1<R<1$ such that
\begin{equation*}
\mathds{1}_{Q_{R_2}(t_0,x_0)}\vu \in \M_{t,x}^{3,\sigma}(\mathbb{R}\times \R),
\end{equation*}
for some $\sigma$ close to $\frac{5}{1-\alpha}<\tau_0<6$ such that $\tau_0<\sigma$.
\end{Proposition}
\textbf{Proof of the Proposition \ref{Propo_GainenMorreySigma}.} In order to obtain this small additional gain of integrability we will first localize the variable $\vu$ in a suitable manner and then we will study its evolution: the wished result will then be deduced from the Duhamel formula and from all the available information over $\vu$. Let us start fixing the parameters $\mathfrak{R}_c, \mathfrak{R}_b, \mathfrak{R}_a$ such that
$$0<R_2<\mathfrak{R}_c<\mathfrak{R}_b<\mathfrak{R}_a<R_1,$$ 
with the associated parabolic balls $Q_{R_2}\subset Q_{\mathfrak{R}_c}\subset Q_{\mathfrak{R}_b}\subset Q_{\mathfrak{R}_a}\subset Q_{R_1}$ (all centered in the point $(t_0,x_0)$). Consider now $ \phi, \psi: \mathbb{R}\times\mathbb{R}^3\longrightarrow \mathbb{R}$ two  non-negative functions such that $\phi, \psi\in \mathcal{C}_0^{\infty}
(\mathbb{R}\times \mathbb{R}^3)$ and such that
\begin{equation}\label{ProprieteLocalisation_IterationVU}
 \phi \equiv 1\;\; \text{over}\; \; Q_{\mathfrak{R}_c},\; \; supp( \phi)\subset Q_{\mathfrak{R}_b} \quad \mbox{and}\quad \psi \equiv 1\;\; \text{over}\; \; Q_{\mathfrak{R}_a},\; \; supp(\psi)\subset Q_{R_1}.
\end{equation}
Using these auxiliar functions we will study the evolution of the variable $\vv= \phi\,\vu$ given by the system
\begin{equation}\label{equationV}
\begin{cases}
\partial_t\vv=\Delta\vv+\vV,\\[3mm]
\vv(0,x)=0,
\end{cases}
\end{equation}
where we have  
\begin{equation}\label{Formule_equationV}
\vV=(\partial_t\phi- \Delta \phi)\vu-2\sum_{i=1}^3 (\partial_i \phi)(\partial_i \vu)-\phi(\vu \cdot \vn)\vu-2\phi\vn p+\phi(\theta e_3).
\end{equation}
We will now rewrite the term $\phi\vn p$ above in order to avoid a direct derivative over the pressure. Indeed, as we have the identity $p= \psi p$ over $Q_{\mathfrak{R}_a}$, then over the smaller ball $Q_{R_2}$  (recalling that $ \psi=1$ over $Q_{R_2}$ by (\ref{ProprieteLocalisation_IterationVU}) since $Q_{R_2}\subset Q_{\mathfrak{R}_a}$), we can write $-\Delta( \psi p)=- \psi \Delta p+(\Delta  \psi)p-2\displaystyle{\sum_{i=1}^3}\partial_i((\partial_i \psi)p)$
from which we deduce the identity 
\begin{equation}\label{ExpressionPression6}
 \phi \vn p= \phi\frac{\vn (- \psi \Delta p)}{(-\Delta)}+ \phi\frac{\vn ((\Delta  \psi)p)}{(-\Delta)}-2\sum_{i=1}^3 \phi\frac{\vn (\partial_i((\partial_i \psi)p))}{(-\Delta)}.
\end{equation}
At this point we recall that we have by (\ref{Equation_Pressure_intro}) the following equation for the pressure 
$$\Delta p = -\displaystyle{\sum_{i,j=1}^3}\partial_i\partial_j\left(u_i u_j\right)+\partial_{x_3}\theta,$$ 
and thus, the first term of the right-hand side of the previous formula can be written in the following manner:
\begin{eqnarray*}
 \phi\frac{\vn (- \psi \Delta p)}{(-\Delta)}&= & \phi\frac{\vn}{(-\Delta)}\left( \sum_{i,j=1}^3 \psi\big(\partial_i\partial_j u_i u_j\big)\right)- \phi\frac{\vn}{(-\Delta)}\left(\psi\partial_{x_3}\theta \right)\notag\\
&=& \sum_{i,j=1}^3  \phi\frac{\vn}{(-\Delta)}\bigg(\partial_i\partial_j( \psi u_i u_j)\bigg)-\sum_{i,j=1}^3  \phi\frac{\vn}{(-\Delta)}\bigg(\partial_i((\partial_j  \psi)u_i u_j)+\partial_j((\partial_i  \psi)u_i u_j)-(\partial_i \partial_j  \psi)(u_i u_j)\bigg)\\
&&-\phi\frac{\vn}{(-\Delta)}\partial_{x_3}\left(\psi\theta \right)+\phi\frac{\vn}{(-\Delta)}\left((\partial_{x_3}\psi)\theta \right),
\end{eqnarray*}
Recalling that by construction of the auxiliar functions $ \phi$, $ \psi$ given in (\ref{ProprieteLocalisation_IterationVU}) we have the identity $ \phi  \psi= \phi$, we can write for the first term above:
$$ \phi\frac{\vn}{(-\Delta)}\partial_i\partial_j( \psi u_i u_j)=\left[ \phi, \frac{\vn\partial_i\partial_j}{(-\Delta)}\right]( \psi u_i u_j)+ \frac{\vn\partial_i\partial_j}{(-\Delta)}( \phi u_i u_j),$$
and we finally obtain the following expression for (\ref{ExpressionPression6}):
\begin{eqnarray*}
 \phi \vn p&=& \sum_{i,j=1}^3\left[ \phi, \frac{\vn\partial_i\partial_j}{(-\Delta)}\right]( \psi u_i u_j)+ \sum_{i,j=1}^3\frac{\vn\partial_i\partial_j}{(-\Delta)}( \phi u_i u_j)\\
&&-\sum_{i,j=1}^3  \phi\frac{\vn}{(-\Delta)}\bigg(\partial_i((\partial_j  \psi)u_i u_j)+\partial_j((\partial_i  \psi)u_i u_j)-(\partial_i \partial_j  \psi)(u_i u_j)\bigg)\\
&&-\phi\frac{\vn}{(-\Delta)}\partial_{x_3}\left(\psi\theta \right)+\phi\frac{\vn}{(-\Delta)}\left((\partial_{x_3}\psi)\theta \right)+ \phi\frac{\vn ((\Delta  \psi)p)}{(-\Delta)}-2\sum_{i=1}^3 \phi\frac{\vn (\partial_i((\partial_i \psi)p))}{(-\Delta)}.
\end{eqnarray*}
With this expression for the term that contains the pressure $p$, we obtain the (lengthy) formula for (\ref{Formule_equationV}):
\begin{align}
&\vV= \underbrace{(\partial_t  \phi - \Delta  \phi)\vu}_{(1)}-2\sum_{i=1}^3 \underbrace{(\partial_i  \phi)(\partial_i \vu)}_{(2)}-\underbrace{ \phi(\vu\cdot \vn)\vu}_{3}-\sum_{i,j=1}^3\underbrace{\left[ \phi,\frac{\vn\partial_i\partial_j}{(-\Delta)}\right]( \psi u_i u_j)}_{(4)}+\sum_{i,j=1}^3 \underbrace{\frac{\vn\partial_i\partial_j}{(-\Delta)}( \phi u_i u_j)}_{(5)}\notag\\
&- \sum_{i,j=1}^3 \frac{ \phi\vn}{(-\Delta)}\big[\underbrace{\partial_i((\partial_j  \psi)u_i u_j)}_{(6)}+\underbrace{\partial_j((\partial_i  \psi)u_i u_j)}_{(7)}-\underbrace{(\partial_i \partial_j  \psi)(u_i u_j)}_{(8)}\big]-\underbrace{\phi\frac{\vn}{(-\Delta)}\partial_{x_3}\left(\psi\theta \right)}_{(9)}\label{FormulePourN}\\
&+\underbrace{\phi\frac{\vn}{(-\Delta)}\left((\partial_{x_3}\psi)\theta \right)}_{(10)}+2\underbrace{ \phi\frac{\vn ((\Delta  \psi)p)}{(-\Delta)}}_{(11)}-2\sum_{i=1}^3\underbrace{ \phi\frac{\vn (\partial_i((\partial_i \psi)p))}{(-\Delta)}}_{(12)}+\underbrace{\phi(\theta e_3)}_{(13)}.\notag
\end{align}
Thus, by the Duhamel formula, the solution $\vv$ of the equation  (\ref{equationV}) is given by
$$\vv=\int_0^{t}e^{(t-s)\Delta}\vV(s,\cdot)ds=\sum_{k=1}^{13}\int_0^{t}e^{(t-s)\Delta}\vV_k(s,\cdot)ds=\sum_{k=1}^{13}\vec{\mathbb{V}}_k.$$
Since $\vv=\phi\vu$, and due to the support properties of $\phi$ (see (\ref{ProprieteLocalisation_IterationVU})), we have $\mathds{1}_{Q_{R_2}}\vv=\mathds{1}_{Q_{R_2}}\vu$ and to conclude that $\mathds{1}_{Q_{R_2}}\vu \in \M_{t,x}^{3,\sigma}(\mathbb{R}\times \R) $ we will study $\mathds{1}_{Q_{R_2}}\vec{\mathbb{V}}_k$ for all $1\leq k\leq 13$. 
\begin{itemize}
\item For $\vec{\mathbb{V}}_1$, by the term (1) in (\ref{FormulePourN}) we have
\begin{equation}\label{Estimation_Ponctuelle_Riesz0}
|\mathds{1}_{Q_{R_2}}\vec{\mathbb{V}}_1(t,x)|=\left|\mathds{1}_{Q_{R_2}}\int_{0}^{t} e^{(t-s)\Delta}[(\partial_t \phi - \Delta \phi)\vu](s,x)ds\right|,
\end{equation}
since the convolution kernel of the semi-group $e^{(t-s)\Delta}$ is the  usual 3D heat kernel $\mathfrak{g}_{t}$, we can write by the decay properties of the heat kernel as well as the properties of the test function $\phi$ (see (\ref{ProprieteLocalisation_IterationVU})), the estimate
$$|\mathds{1}_{Q_{R_2}}\vec{\mathbb{V}}_1(t,x)|\leq C\mathds{1}_{Q_{R_2}}\int_{\mathbb{R}} \int_{\mathbb{R}^3} \frac{1}{(|t-s|^{\frac{1}{2}}+|x-y|)^3}\left| \mathds{1}_{Q_{\mathfrak{R}_b}}\vu(s,y)
\right| \,dy \,ds,$$
Now, recalling the definition of the parabolic Riesz potential given in (\ref{Def_ParabolicRieszPotential}) and since $Q_{R_2}\subset Q_{\mathfrak{R}_b}$ we obtain the pointwise estimate 
\begin{equation}\label{Estimation_Ponctuelle_Riesz1}
|\mathds{1}_{Q_{R_2}}\vec{\mathbb{V}}_1(t,x)|\leq C\mathds{1}_{Q_{\mathfrak{R}_b}}\mathcal{L}_{2}(|\mathds{1}_{Q_{\mathfrak{R}_b}}\vu|)(t,x),
\end{equation}
and taking the Morrey $\mathcal{M}_{t,x}^{3, \sigma}$-norm we obtain
$$\|\mathds{1}_{Q_{R_2}}\vV_1(t,x)\|_{\mathcal{M}_{t,x}^{3, \sigma}}\leq C\|\mathds{1}_{Q_{\mathfrak{R}_b}}\mathcal{L}_{2}(  |\mathds{1}_{Q_{\mathfrak{R}_b}} \vu|)\|_{\mathcal{M}_{t,x}^{3, \sigma}}.$$
Now, for some $2<q<\frac{5}{2}$ we set $\lambda=1-\frac{2q}{5}$. Then, we have $ 3\le\frac{3}{\lambda}$ and $\sigma\leq \frac{q}{\lambda}$. Thus, by Lemma \ref{lemma_locindi} and by Lemma \ref{Lemma_Hed} we can write:
\begin{eqnarray*}
\|\mathds{1}_{Q_{\mathfrak{R}_b}}\mathcal{L}_{2}(  |\mathds{1}_{Q_{\mathfrak{R}_b}}  \vu |)\|_{\mathcal{M}_{t,x}^{3, \sigma}}&\leq &C\|\mathcal{L}_{2}(  |\mathds{1}_{Q_{\mathfrak{R}_b}} \vu|)\|_{\mathcal{M}_{t,x}^{\frac{3}{\lambda},\frac{q}{\lambda}}}\notag\\
&\leq &C\|\mathds{1}_{Q_{\mathfrak{R}_b}}  \vu\|_{\mathcal{M}_{t,x}^{3, q}}\leq C\|\mathds{1}_{Q_{R_1}} \vu\|_{\mathcal{M}_{t,x}^{3, \tau_0}}<+\infty,
\end{eqnarray*}
where in the last estimate we applied again Lemma \ref{lemma_locindi}  (noting that $q<\tau_0<6$) and we used the estimates over $\vu$ available in (\ref{Conclusion_FirstMorreySpace}).
\item For $\vec{\mathbb{V}}_2$, using the expression (2) in (\ref{FormulePourN}) we write $(\partial_{i}\phi) (\partial_i\vu)=\partial_i((\partial_i\phi)\vu)-(\partial_i^2\phi)\vu$ and we have
\begin{equation}\label{EstimationPonctuelleVV2}
|\mathds{1}_{Q_{R_2}}\vec{\mathbb{V}}_2(t,x)|\leq \sum_{i=1}^{3}\left|\mathds{1}_{Q_{R_2}}\int_{0}^{t} e^{ (t-s)\Delta} \partial_i\big((\partial_{i}\phi) \vu\big)ds\right|+\left|\mathds{1}_{Q_{R_2}}\int_{0}^{t} e^{ (t-s)\Delta} [(\partial_{i}^2\phi)\vu] ds\right|.
\end{equation}
Remark that the second term of the right-hand side of  (\ref{EstimationPonctuelleVV2}) can be treated in the same manner as the term $\vec{\mathbb{V}}_1$ so we will only study the first term: by the properties of the heat kernel and by the definition of the Riesz potential $\mathcal{L}_{1}$ (see (\ref{Def_ParabolicRieszPotential})), we obtain
\begin{eqnarray}
A_2:=\left|\mathds{1}_{Q_{R_2}}\int_{0}^{t} e^{ (t-s)\Delta} \partial_i\big((\partial_{i}\phi) \vu \big)ds\right|=\left|\mathds{1}_{Q_{R_2}}\int_{0}^{t} \int_{\mathbb{R}^3}\partial_i\mathfrak{g}_{t-s}(x-y)(\partial_{i}\phi) \vu(s,y)dyds\right|\notag\\
\leq C\mathds{1}_{Q_{R_2}}\int_{\mathbb{R}}\int_{\mathbb{R}^3} \frac{|\mathds{1}_{Q_{\mathfrak{R}_b}} \vu(s,y)|}{(|t-s|^{\frac{1}{2}}+|x-y|)^4}dyds
\leq C\mathds{1}_{Q_{R_2}}(\mathcal{L}_1(|\mathds{1}_{Q_{\mathfrak{R}_b}} \vu |))(t,x).\label{Estimation_Ponctuelle_Riesz2}
\end{eqnarray}
Taking the Morrey $\mathcal{M}_{t,x}^{3, \sigma}$ norm we obtain $\|A_2\|_{\mathcal{M}_{t,x}^{3, \sigma}}\leq C\|\mathds{1}_{Q_{R_2}}(\mathcal{L}_1(|\mathds{1}_{Q_{\mathfrak{R}_b}} \vu|))\|_{\mathcal{M}_{t,x}^{3, \sigma}}$.
Now, for some $4\leq q <5<\frac{5}{1-\alpha}<\tau_0<6$ we define $\lambda=1-\frac{q}{5}$, noting that $3\leq \frac{3}{\lambda}$ and $\sigma \leq\frac{q}{\lambda}$, by Lemma \ref{Lemma_Hed}, we can write
\begin{eqnarray*}
\|\mathds{1}_{Q_{R_2}}(\mathcal{L}_1(|\mathds{1}_{Q_{\mathfrak{R}_b}} \vu|))\|_{\mathcal{M}_{t,x}^{3, \sigma}}&\leq &C\|\mathcal{L}_1(|\mathds{1}_{Q_{\mathfrak{R}_b}} \vu|)\|_{\mathcal{M}_{t,x}^{\frac{3}{\lambda}, \frac{q}{\lambda}}}\leq C\|\mathds{1}_{Q_{\mathfrak{R}_b}} \vu\|_{\mathcal{M}_{t,x}^{3, q}}\\
&\leq & C\|\mathds{1}_{Q_{R_1}} \vu\|_{\mathcal{M}_{t,x}^{3, \tau_0}}<+\infty,
\end{eqnarray*}
from which we deduce that $\|\mathds{1}_{Q_{R_2}}\vec{\mathbb{V}}_2\|_{\mathcal{M}_{t,x}^{3, \sigma}}<+\infty$.

\item For the term $\vec{\mathbb{V}}_3$, by the same arguments given to obtain the pointwise estimate (\ref{Estimation_Ponctuelle_Riesz1}), we have
\begin{eqnarray*}
|\mathds{1}_{Q_{R_2}}\vec{\mathbb{V}}_3(t,x)|&=&\left|\mathds{1}_{Q_{R_2}}\int_{0}^{t} \int_{\mathbb{R}^3}  \mathfrak{g}_{t-s}(x-y)\left[\phi \left( (\vu\cdot\vn)\vu\right)\right](s,y)dyds\right|\\
&\leq & C\mathds{1}_{Q_{R_2}}\mathcal{L}_2\left(\left|\mathds{1}_{Q_{\mathfrak{R}_b}}\left( (\vu\cdot\vn)\vu\right)\right|\right)(t,x),
\end{eqnarray*}
(recall (\ref{ProprieteLocalisation_IterationVU})) from which we deduce 
\begin{equation}\label{Decomposition2termesVV3}
\|\mathds{1}_{Q_{R_2}}\vec{\mathbb{V}}_3\|_{\mathcal{M}_{t,x}^{3, \sigma}}\leq C\left\|\mathds{1}_{Q_{R_2}}\mathcal{L}_2\left(|\mathds{1}_{Q_{\mathfrak{R}_b}} (\vu\cdot\vn)\vu|\right)\right\|_{\mathcal{M}_{t,x}^{3, \sigma}}.
\end{equation}
We set now $\frac{5}{3-\alpha}<q<\frac{5}{2}$ and $\lambda=1-\frac{2q}{5}$. Since $3\leq \frac{6}{5\lambda}$ and $\tau_0<6<\sigma\leq \frac{q}{\lambda}$, applying Lemma \ref{lemma_locindi} and Lemma \ref{Lemma_Hed} we have
$$\left\|\mathds{1}_{Q_{R_2}}\mathcal{L}_2\left(|\mathds{1}_{Q_{\mathfrak{R}_b}} (\vu\cdot\vn)\vu|\right)\right\|_{\mathcal{M}_{t,x}^{3, \sigma}}\leq C\left\|\mathds{1}_{Q_{R_2}}\mathcal{L}_2\left(|\mathds{1}_{Q_{\mathfrak{R}_b}} (\vu\cdot\vn)\vu|\right)\right\|_{\mathcal{M}_{t,x}^{\frac{6}{5\lambda},\frac{q}{\lambda}}}\leq C\left\|\mathds{1}_{Q_{\mathfrak{R}_b}} (\vu\cdot\vn)\vu\right\|_{\mathcal{M}_{t,x}^{\frac{6}{5}, q}}.$$
Recall that we have $\tau_0<6<\sigma$ and by the H\"older inequality in Morrey spaces (see Lemma \ref{lemma_Product}) we obtain
$$\left\|\mathds{1}_{Q_{\mathfrak{R}_b}} (\vu\cdot\vn)\vu\right\|_{\mathcal{M}_{t,x}^{\frac{6}{5}, q}}\leq \left\|\mathds{1}_{Q_{R_1}}\vu\right\|_{\mathcal{M}_{t,x}^{3, \tau_0}}\left\|\mathds{1}_{Q_{R_1}}\vn \otimes \vu\right\|_{\mathcal{M}_{t,x}^{2, \tau_1}}<+\infty,$$
where $\frac{1}{q}=\frac{1}{\tau_0}+\frac{1}{\tau_1}=\frac{2}{\tau_0}+\frac{1}{5}$. These two last quantities are bounded by (\ref{Conclusion_FirstMorreySpace}) and (\ref{ConlusioncorolarioMorrey}). Note that the condition $\tau_0<6<\sigma$ and the relationship $\frac{1}{q}=\frac{2}{\tau_0}+\frac{1}{5}$ are compatible with the fact that $\frac{5}{3-\alpha}<q < \frac{5}{2}$ (recall that $0<\alpha\ll\frac{1}{10}$). 
\item The term $\vec{\mathbb{V}}_4$ is the most technical one. Indeed, by the expression of $\vV_4$ given in (\ref{FormulePourN}), we write
\begin{eqnarray*}
|\mathds{1}_{Q_{R_2}}\vec{\mathbb{V}}_4|\leq \sum^3_{i,j= 1} \mathds{1}_{Q_{R_2}}\int_{\mathbb{R}} \int_{\mathbb{R}^3}\frac{\left|\left[\phi, \, \frac{ \vn \partial_i\partial_j}{(-\Delta)}\right] (\psi u_i u_j)(s,y)\right|}{(|t-s|^{\frac{1}{2}}+|x-y|)^3}dyds\leq \sum^3_{i,j= 1}\mathds{1}_{Q_{R_2}}\mathcal{L}_{2}\left(\left|\left[\phi, \, \frac{ \vn \partial_i\partial_j}{(-\Delta)}\right] (\psi u_i u_j)\right|\right),
\end{eqnarray*}
and taking the $\mathcal{M}_{t,x}^{3, \sigma}$-norm we have
$\|\mathds{1}_{Q_{R_2}}\vec{\mathbb{V}}_4\|_{\mathcal{M}_{t,x}^{3, \sigma}}\leq \sum^3_{i,j= 1}\left\|\mathds{1}_{Q_{R_2}}\mathcal{L}_{2}\left(\left|\left[\phi, \, \frac{ \vn \partial_i\partial_j}{(-\Delta)}\right] (\psi u_i u_j)\right|\right)\right\|_{\mathcal{M}_{t,x}^{3, \sigma}}$. If we set $\frac{1}{q}=  \frac{2}{\tau_0} + \frac{1}{5}$ and $\lambda=1-\frac{2q}{5}$ then we have $3\leq \frac{3}{2\lambda}$ and for
\begin{equation}\label{UpperBound_sigma}
\sigma \leq \frac{q}{\lambda} = \frac{5 \tau_0}{10- \tau_0},
\end{equation}
by Lemmas \ref{lemma_locindi} and \ref{Lemma_Hed} we obtain:
\begin{eqnarray*}
\left\|\mathds{1}_{Q_{R_2}}\mathcal{L}_{2}\left(\left|\left[\phi, \, \frac{ \vn \partial_i\partial_j}{(-\Delta)}\right] (\psi u_i u_j)\right|\right)\right\|_{\mathcal{M}_{t,x}^{3, \sigma}}&\leq &C\left\|\mathds{1}_{Q_{R_2}}\mathcal{L}_{2}\left(\left|\left[\phi, \, \frac{ \vn \partial_i\partial_j}{(-\Delta)}\right] (\psi u_i u_j)\right|\right)\right\|_{\mathcal{M}_{t,x}^{\frac{3}{2\lambda}, \frac{q}{\lambda}}}\\
&\leq &C\left\|\left[\phi, \, \frac{ \vn \partial_i\partial_j}{(-\Delta)}\right] (\psi u_i u_j)\right\|_{\mathcal{M}_{t,x}^{\frac{3}{2}, q}},
\end{eqnarray*}
We will study this norm and by the definition of Morrey spaces (\ref{DefMorreyparabolico}), if we introduce a threshold $\mathfrak{r}=\frac{\mathfrak{R}_b-R_2}{2}$, we have
\begin{equation}\label{CommutatorEstimatevVV4}
\begin{split}
\left\|\left[\phi, \, \frac{ \vn \partial_i\partial_j}{(-\Delta)}\right] (\psi u_i u_j)\right\|_{\mathcal{M}_{t,x}^{\frac{3}{2}, q}}^{\frac{3}{2}}&\leq\underset{\underset{0<r<\mathfrak{r}}{(\mathfrak{t},\bar{x})}}{\sup}\;\frac{1}{r^{5(1-\frac{3}{2q})}}\int_{Q_r(\mathfrak{t},\bar x)}\left|\left[\phi, \, \frac{ \vn \partial_i\partial_j}{(-\Delta)}\right] (\psi u_i u_j)\right|^{\frac{3}{2}}dxdt\qquad\\
&+\underset{\underset{\mathfrak{r}< r}{(\mathfrak{t},\bar{x})}}{\sup}\;\frac{1}{r^{5(1-\frac{3}{2q})}}\int_{Q_r(\mathfrak{t},\bar{x})}\left|\left[\phi, \, \frac{ \vn \partial_i\partial_j}{(-\Delta)}\right] (\psi u_i u_j)\right|^{\frac{3}{2}}dxdt.\qquad
\end{split}
\end{equation}
Now, we study the second term of the right-hand side above, which is easy to handle as we have $\mathfrak{r}<r$ and we can write
$$\underset{\underset{\mathfrak{r}< r}{(\mathfrak{t},\bar{x}) \in \mathbb{R}\times \R}}{\sup}\;\frac{1}{r^{5(1-\frac{3}{2q})}}\int_{Q_r(\mathfrak{t},\bar{x})}\left|\left[\phi, \, \frac{ \vn \partial_i\partial_j}{(-\Delta)}\right] (\psi u_i u_j)\right|^{\frac{3}{2}}dxdt\leq C_{\mathfrak{r}}\left\|\left[\phi, \, \frac{ \vn \partial_i\partial_j}{(-\Delta)}\right] (\psi u_i u_j)\right\|_{L^{\frac{3}{2}}_{t,x}}^{\frac{3}{2}},$$
and since $\bar\phi$ is a regular function and $\frac{ \vn \partial_i\partial_j}{(-\Delta)}$ is a Calder\'on-Zydmund operator, by the Calder\'on commutator theorem (see the book \cite{PGLR0}), we have that the operator $\left[\phi, \, \frac{ \vn \partial_i\partial_j}{(-\Delta)}\right]$ is bounded in the space $L^{\frac{3}{2}}_{t,x}$ and we can write (using the support properties of $\psi$ given in (\ref{ProprieteLocalisation_IterationVU}) and the information given in (\ref{Conclusion_FirstMorreySpace})):
\begin{eqnarray*}
\left\|\left[\bar\phi, \, \frac{ \vn \partial_i\partial_j}{(-\Delta)}\right] (\psi u_i u_j)\right\|_{L^{\frac{3}{2}}_{t,x}}&\leq &C\left\|\psi u_i u_j\right\|_{L^{\frac{3}{2}}_{t,x}}\leq C\|\mathds{1}_{Q_{R_1}} u_i u_j\|_{\mathcal{M}^{\frac{3}{2}, \frac{3}{2}}_{t,x}}\\
&\leq & C\|\mathds{1}_{Q_{R_1}} \vu\|_{\mathcal{M}^{3,3}_{t,x}}\|\mathds{1}_{Q_{R_1}} \vu\|_{\mathcal{M}^{3, 3}_{t,x}}\leq C\|\mathds{1}_{Q_{R_1}} \vu\|_{\mathcal{M}^{3,\tau_0}_{t,x}}\|\mathds{1}_{Q_{R_1}} \vu\|_{\mathcal{M}^{3, \tau_0}_{t,x}}<+\infty,
\end{eqnarray*}
where in the last line we used H\"older inequalities in Morrey spaces and we applied Lemma \ref{lemma_locindi}.\\

The first term of the right-hand side of  (\ref{CommutatorEstimatevVV4}) requires some extra computations: indeed, as we are interested to obtain information over the parabolic ball $Q_{r}(\mathfrak{t}, \bar{x})$ we can write
for some  $0<r<\mathfrak{r}$:
\begin{equation}\label{CommutatorEstimatevVV401}
\mathds{1}_{Q_{r}}\left[\phi, \, \frac{ \vn \partial_i\partial_j}{(-\Delta)}\right] (\psi u_i u_j))=\mathds{1}_{Q_{r}}\left[\phi, \, \frac{ \vn \partial_i\partial_j}{(-\Delta)}\right] (\mathds{1}_{Q_{2r}}\psi u_i u_j)+\mathds{1}_{Q_{r}}\left[\phi, \, \frac{ \vn \partial_i\partial_j}{(-\Delta)}\right] ((\mathbb{I}-\mathds{1}_{Q_{2r}})\psi u_i u_j),
\end{equation}
and as before we will study the $L^{\frac{3}{2}}_{t,x}$ norm of these two terms. For the first quantity in the right-hand side of (\ref{CommutatorEstimatevVV401}), by the Calder\'on commutator theorem, by the definition of Morrey spaces and by the H\"older inequalities we have 
\begin{eqnarray*}
\left\|\mathds{1}_{Q_{r}}\left[\phi, \, \frac{ \vn \partial_i\partial_j}{(-\Delta)}\right] (\mathds{1}_{Q_{2r}}\psi u_i u_j)\right\|_{L^{\frac{3}{2}}_{t,x}}^{\frac{3}{2}}&\leq& C\|\mathds{1}_{Q_{2r}}\psi u_i u_j\|_{L^{\frac{3}{2}}_{t,x}}^{\frac{3}{2}}\leq Cr^{5 (1-\frac{3}{\tau_0})} \|\mathds{1}_{Q_{R_1}}u_i u_j\|_{\mathcal{M}^{\frac{3}{2}, \frac{\tau_0}{2}}_{t,x}}^{\frac{3}{2}}\\
&\leq & Cr^{5 (1-\frac{3}{\tau_0})} \|\mathds{1}_{Q_{R_1}}\vu\|_{\mathcal{M}^{3, \tau_0}_{t,x}}^{\frac{3}{2}}\|\mathds{1}_{Q_{R_1}}\vu\|_{\mathcal{M}^{3, \tau_0}_{t,x}}^{\frac{3}{2}},
\end{eqnarray*}
for all $0<r<\mathfrak{r}$, from which we deduce that 
$$\underset{\underset{0<r<\mathfrak{r}}{(\mathfrak{t},\bar{x}) }}{\sup}\;\frac{1}{r^{5(1-\frac{3}{2q})}}\int_{Q_r(\mathfrak{t},\bar{x})}\left|\mathds{1}_{Q_{r}}\left[\phi, \, \frac{ \vn \partial_i\partial_j}{(-\Delta)}\right] (\mathds{1}_{Q_{2r}}\psi u_i u_j)\right|^{\frac{3}{2}}dxdt\leq C \|\mathds{1}_{Q_{R_1}}\vu\|_{\mathcal{M}^{3, \tau_0}_{t,x}}^{\frac{3}{2}}\|\mathds{1}_{Q_{R_1}}\vu\|_{\mathcal{M}^{3, \tau_0}_{t,x}}^{\frac{3}{2}}<+\infty.$$
We study now the second term of the right-hand side of (\ref{CommutatorEstimatevVV401}) and for this we consider the following operator:
$$T: f \mapsto  \left(\mathds{1}_{Q_{r}}  \left[\phi, \, \frac{ \vn \partial_i \partial_j}{- \Delta }\right] (\mathbb{I} - \mathds{1}_{Q_{2r}})  \psi
\right) f,$$
and by the properties of the convolution kernel of the operator $\frac{1}{(-\Delta)}$ we obtain
$$|T(f)(x)|\leq C\mathds{1}_{Q_{r}}(x)\int_{\mathbb{R}^3}\frac{(\mathbb{I} - \mathds{1}_{Q_{2r}})(y) \mathds{1}_{Q_{R_1}}(y) |f(y)| | \phi(x)- \phi(y)|}{|x-y|^4} dy.$$
Recalling that $0<r<\mathfrak{r}=\frac{\mathfrak{R}_b-R_2}{2}$, by the support properties of the test function $\phi$ (see (\ref{ProprieteLocalisation_IterationVU})), the integral above is meaningful if $|x-y|>r$ and thus we can write
\begin{eqnarray*}
\left\|\mathds{1}_{Q_{r}}\left[\phi, \, \frac{ \vn \partial_i\partial_j}{(-\Delta)}\right] ((\mathbb{I}-\mathds{1}_{Q_{2r}})\psi u_i u_j)\right\|_{L^\frac{3}{2}_{t,x}}^\frac{3}{2}\leq C\left\|\mathds{1}_{Q_{r}} \int_{\mathbb{R}^3} \frac{\mathds{1}_{|x-y| > r}}{|x-y|^4}(\mathbb{I} - \mathds{1}_{Q_{2r}})(y) \mathds{1}_{Q_{R_1}}(y)|u_i u_j |dy\right\|_{L^\frac{3}{2}_{t,x}}^\frac{3}{2}\\
\leq C\left(\int_{|y|>r}\frac{1}{|y|^4}\| \mathds{1}_{Q_{R_1}}
|u_i u_j |(\cdot-y)\|_{L^\frac{3}{2}_{t,x}(Q_r)}dy\right)^\frac{3}{2}\leq Cr^{-\frac{3}{2}}\| \mathds{1}_{Q_{R_1}}u_i u_j \|_{L^\frac{3}{2}_{t,x}(Q_{r})}^\frac{3}{2},
\end{eqnarray*}
with this estimate at hand and using the definition of Morrey spaces, we can write
\begin{eqnarray*}
\int_{Q_r(\mathfrak{t},\bar{x})}\left|\mathds{1}_{Q_{r}}\left[\phi, \, \frac{ \vn \partial_i\partial_j}{(-\Delta)}\right] ((\mathbb{I}-\mathds{1}_{Q_{2r}})\psi u_i u_j)\right|^{\frac{3}{2}}dxdt&\leq &Cr^{-\frac{3}{2}}r^{5 (1-\frac{3}{\tau_0})}\| \mathds{1}_{Q_{R_1}}u_i u_j \|_{\mathcal{M}^{\frac{3}{2}, \frac{\tau_0}{2}}_{t,x}}^\frac{3}{2}\\
&\leq &Cr^{5(1-\frac{3}{2q})}\| \mathds{1}_{Q_{R_1}}u_i u_j \|_{\mathcal{M}^{\frac{3}{2}, \frac{\tau_0}{2}}_{t,x}}^\frac{3}{2},
\end{eqnarray*}
where in the last inequality we used the fact that $\frac{1}{q}=  \frac{2}{\tau_0} + \frac{1}{5}$, which implies  $r^{-\frac{3}{2}}r^{5 (1-\frac{3}{\tau_0})}= r^{5(1-\frac{3}{2q})}$. Thus we finally obtain
$$\underset{\underset{0<r<\mathfrak{r}}{(\mathfrak{t},\bar{x}) }}{\sup}\;\frac{1}{r^{5(1-\frac{3}{2q})}}\int_{Q_r(\mathfrak{t},\bar{x})}\left|\mathds{1}_{Q_{r}}\left[\phi, \, \frac{ \vn \partial_i\partial_j}{(-\Delta)}\right] ((\mathbb{I}-\mathds{1}_{Q_{2r}})\psi u_i u_j)\right|^{\frac{3}{2}}dxdt\leq C \|\mathds{1}_{Q_{R_1}}\vu\|_{\mathcal{M}^{3, \tau_0}_{t,x}}^{\frac{3}{2}}\|\mathds{1}_{Q_{R_1}}\vu\|_{\mathcal{M}^{3, \tau_0}_{t,x}}^{\frac{3}{2}}<+\infty.$$
We have proven that all the term in (\ref{CommutatorEstimatevVV4}) are bounded and we can conclude that $\|\mathds{1}_{Q_{R_2}}\vec{\mathbb{V}}_4\|_{\mathcal{M}_{t,x}^{3, \sigma}}<+\infty$.
\begin{Remark}
The condition \eqref{UpperBound_sigma} implies an upper bound for $\sigma$ depending on the current Morrey information of $\vu$, which a priori is close to $\tau_0$ with $\frac{5}{1-\alpha}<\tau_0<6$. Nevertheless it is clear that whether we obtain a better Morrey information on integrability for $\vu$, the value of $\sigma$ can increase.
\end{Remark}
\item For the quantity $\vec{\mathbb{V}}_5$, based in the expression (\ref{FormulePourN}) we write
\begin{eqnarray*}
|\mathds{1}_{Q_{R_2}}\vec{\mathbb{V}}_5(t,x)|&\leq &C\sum^3_{i,j= 1} \mathds{1}_{Q_{R_2}}\int_{\mathbb{R}}\int_{\mathbb{R}^3} \frac{|\mathcal{R}_i\mathcal{R}_j( \phi u_i u_j )(s,y)|}{(|t-s|^{\frac{1}{2}}+|x-y|)^4}dyds\leq C\sum^3_{i,j= 1} \mathds{1}_{Q_{R_2}} \mathcal{L}_{1}\left(|\mathcal{R}_i\mathcal{R}_j( \phi u_i u_j )|\right)(t,x),
\end{eqnarray*}
where we used the decaying properties of the heat kernel (recall that $\mathcal{R}_i=\frac{\partial_i}{\sqrt{- \Delta}}$ are the Riesz transforms). 
Now taking the Morrey $\mathcal{M}^{3, \sigma}_{t,x}$ norm and by Lemma \ref{lemma_locindi} (with $\nu=\frac{4\tau_0+5}{5\tau_0}$, $p=3$, $q=\tau_0$ such that $\frac{p}{\nu}>3$ and $\frac{q}{\nu}>\sigma$ which is compatible with the condition $\tau_0<\sigma$) we have
\begin{eqnarray*}
\|\mathds{1}_{Q_{R_2}}\vV_5\|_{\mathcal{M}^{3, \sigma}_{t,x}}&\leq &C \sum^3_{i,j= 1}\|  \mathds{1}_{Q_{R_2}}\mathcal{L}_{1}\left(|\mathcal{R}_i\mathcal{R}_j( \phi u_i u_j )|\right)\|_{\mathcal{M}^{\frac{p}{\nu}, \frac{q}{\nu}}_{t,x}}
\end{eqnarray*}
Then by Lemma \ref{Lemma_Hed} with $\lambda= 1 - \tfrac{\tau_0 /2}{5}$ (recall $\tau_0<6 <10$ so that $\nu > 2 \lambda$) and by the boundedness of Riesz transforms in Morrey spaces we obtain:
\begin{eqnarray*}
\|\mathds{1}_{Q_{R_2}} \mathcal{L}_{1}\left(|\mathcal{R}_i\mathcal{R}_j( \phi u_i u_j )|\right)\|_{\mathcal{M}^{\frac{p}{\nu}, \frac{q}{\nu}}_{t,x}}
&\leq&C\|\mathcal{L}_{1}\left(|\mathcal{R}_i\mathcal{R}_j( \phi u_i u_j )|\right)\|_{\mathcal{M}^{\frac{p}{2\lambda}, \frac{q}{2\lambda}}_{t,x}}
\leq C \|\mathcal{R}_i\mathcal{R}_j( \phi u_i u_j )\|_{\mathcal{M}^{\frac{3}{2},\frac{\tau_0}{2}}_{t,x}}\\
&\leq &\| \mathds{1}_{Q_{R_1}} u_i u_j \|_{\mathcal{M}^{\frac{3}{2},\frac{\tau_0}{2}}_{t,x}}
\leq  C\|\mathds{1}_{Q_{R_1}} \vu\|_{\mathcal{M}^{3,\tau_0}_{t,x}}\|  \mathds{1}_{Q_{R_1}} \vu\|_{\mathcal{M}^{3,\tau_0}_{t,x}}<+\infty,
\end{eqnarray*}
and we obtain $\|\mathds{1}_{Q_{R_2}}\vec{\mathbb{V}}_5\|_{\mathcal{M}_{t,x}^{3, \sigma}}<+\infty$.
\item For the term $\vec{\mathbb{V}}_6$ and following the same ideas we have
$$|\mathds{1}_{Q_{R_2}}\vec{\mathbb{V}}_6|\leq C\sum^3_{i,j= 1}  \mathds{1}_{Q_{R_2}}\int_{\mathbb{R}}  \int_{\mathbb{R}^3}\frac{\left|\frac{\phi \vn\partial_i}{(- \Delta )} (\partial_j \psi) u_i u_j(s,y)\right|}{(|t-s|^{\frac{1}{2}}+|x-y|)^3} dyds=C\sum^3_{i,j= 1}  \mathds{1}_{Q_{R_2}}\mathcal{L}_{2}\left(\left|\frac{\phi \vn\partial_i}{(- \Delta )} (\partial_j \psi) u_i u_j\right|\right).$$
For $2<q<\frac{5}{2}$, define $\lambda=1-\frac{2q}{5}$, we thus have $3\leq \frac{3}{2\lambda}$ and $\sigma \leq \frac{q}{\lambda}$. Then, by Lemma \ref{lemma_locindi} and Lemma \ref{Lemma_Hed} we can write
\begin{eqnarray*}
\left\|\mathds{1}_{Q_{R_2}}\mathcal{L}_{2}\left|\frac{\phi \vn\partial_i}{(- \Delta )} (\partial_j \psi) u_i u_j\right|\right\|_{\mathcal{M}^{3, \sigma}_{t,x}}\leq C\left\|\mathds{1}_{Q_{R_2}}\mathcal{L}_{2}\left|\frac{\phi \vn\partial_i}{(- \Delta )} (\partial_j \psi) u_i u_j\right|\right\|_{\mathcal{M}^{\frac{3}{2\lambda}, \frac{q}{\lambda}}_{t,x}}\leq  C\left\|\frac{\phi \vn\partial_i}{(- \Delta )} (\partial_j \psi) u_i u_j\right\|_{\mathcal{M}^{\frac{3}{2}, q}_{t,x}},
\end{eqnarray*}
but since the operator $\frac{\phi \vn\partial_i}{(- \Delta )}$ is bounded in Morrey spaces and since $2<q<\frac{5}{2}< \tfrac{\tau_0}{2}<3$ (since $\tau_0<6$), one has by Lemma \ref{lemma_locindi} and by the H\"older inequalities
\begin{eqnarray*}
\left\|\frac{\phi \vn\partial_i}{(- \Delta )} (\partial_j \psi) u_i u_j\right\|_{\mathcal{M}^{\frac{3}{2}, q}_{t,x}}\leq C\left\| \mathds{1}_{Q_{R_1}}u_i u_j\right\|_{\mathcal{M}^{\frac{3}{2}, q}_{t,x}}\leq C\| \mathds{1}_{Q_{R_1}}u_i u_j\|_{\mathcal{M}^{\frac{3}{2}, \frac{\tau_0}{2}}_{t,x}}\leq C\|\mathds{1}_{Q_{R_1}} \vu\|_{\mathcal{M}^{3,\tau_0}_{t,x}}\|  \mathds{1}_{Q_{R_1}} \vu\|_{\mathcal{M}^{3,\tau_0}_{t,x}},
\end{eqnarray*}
from which we deduce $\|\mathds{1}_{Q_{R_2}}\vec{\mathbb{V}}_6\|_{\mathcal{M}^{3, \sigma}_{t,x}}<+\infty$. Note that the same computations can be performed to obtain that $\|\mathds{1}_{Q_{R_2}}\vec{\mathbb{V}}_7\|_{\mathcal{M}^{3, \sigma}_{t,x}}<+\infty$.
\item The quantity $\vec{\mathbb{V}}_8$ based in the term (8) of (\ref{FormulePourN}) is treated in the following manner: we first write
$$\|\mathds{1}_{Q_{R_2}}\vec{\mathbb{V}}_8\|_{\mathcal{M}^{3, \sigma}_{t,x}}\leq C \sum^3_{i,j= 1}\left\| \mathds{1}_{Q_{R_2}}\left(\mathcal{L}_2\left|\phi \frac{\vn}{(-\Delta )}(\partial_i \partial_j \psi) (u_i u_j)\right|\right)\right\|_{\mathcal{M}^{3,\sigma}_{t,x}}.$$
We set $1<\nu<\frac{3}{2}$, $2\nu <q<\frac{5\nu}{2}$ and $\lambda=1-\frac{2q}{5\nu}$, thus we have $3\leq \frac{\nu}{\lambda}$ and $\sigma \leq \frac{q}{\lambda}$, then, by Lemma \ref{lemma_locindi} and by Lemma \ref{Lemma_Hed} we can write
\begin{eqnarray}
\left\|\mathds{1}_{Q_{R_2}}\left(\mathcal{L}_2\left|\phi \frac{\vn}{(-\Delta )}(\partial_i \partial_j \psi) (u_i u_j)\right|\right)\right\|_{\mathcal{M}^{3,\sigma}_{t,x}}\leq C\left\| \mathds{1}_{Q_{R_2}}\left(\mathcal{L}_2\left|\phi \frac{\vn}{(-\Delta )}(\partial_i \partial_j \psi) (u_i u_j)\right|\right)\right\|_{\mathcal{M}^{\frac{\nu}{\lambda},\frac{q}{\lambda}}_{t,x}}\notag\\
\leq C\left\| \phi \frac{\vn}{(-\Delta )}(\partial_i \partial_j \psi) (u_i u_j)\right\|_{\mathcal{M}^{\nu,q}_{t,x}}\leq C\left\| \phi \frac{\vn}{(-\Delta )}(\partial_i \partial_j \psi) (u_i u_j)\right\|_{\mathcal{M}^{\nu,\frac{5\nu}{2}}_{t,x}}\notag\\
\leq C\left\| \phi \frac{\vn}{(-\Delta )}(\partial_i \partial_j \psi) (u_i u_j)\right\|_{L^{\nu}_tL^{\infty}_x}\label{Formula_intermediairevVV80}
\end{eqnarray}
where in the last estimate we used the space inclusion $L^{\nu}_tL^{\infty}_x\subset \mathcal{M}^{\nu,\frac{5\nu}{2}}_{t,x}$. 
\begin{Remark}\label{Remarque_Iteration}
Note that if the parameter $q$ above is close to the value $\frac{5\nu}{2}$, then $\lambda=1-\frac{2q}{5\nu}$ is close to $0$ and thus the value $\frac{q}{\lambda}$ can be made very big: in the estimates (\ref{Formula_intermediairevVV80}) we can consider a Morrey space $\mathcal{M}^{3,\sigma}_{t,x}$ with $\sigma\gg 1$.
\end{Remark}
Let us focus now in the $L^\infty$ norm above (\emph{i.e.} without considering the time variable). Remark that due to the support properties of the auxiliary function $\psi$ given in (\ref{ProprieteLocalisation_IterationVU}) we have $supp(\partial_i \partial_j \psi) \subset Q_{R_1}\setminus Q_{\mathfrak{R}_a}$ and recall by (\ref{ProprieteLocalisation_IterationVU}) we have $supp\;  \phi = Q_{\mathfrak{R}_b}$ where $\mathfrak{R}_b<\mathfrak{R}_a<R_1$, thus by the properties of the kernel of the operator $\frac{\vn}{(-\Delta)}$ we can write
\begin{eqnarray}
\left| \phi \frac{\vn}{(-\Delta )}(\partial_i \partial_j \psi) (u_i u_j)\right|\leq C\left|\int_{\mathbb{R}^3} \frac{1}{|x-y|^2}\mathds{1}_{Q_{\mathfrak{R}_b}}(x)\mathds{1}_{Q_{R_1}\setminus Q_{\mathfrak{R}_a}}(y)(\partial_i \partial_j \psi) (u_i u_j)(\cdot,y)dy\right|\notag\\
\leq  C\left|\int_{\mathbb{R}^3} \frac{\mathds{1}_{|x-y|>(\mathfrak{R}_a-\mathfrak{R}_b)}}{|x-y|^2}\mathds{1}_{Q_{\mathfrak{R}_b}}(x)\mathds{1}_{Q_{R_1}\setminus Q_{\mathfrak{R}_a}}(y)(\partial_i \partial_j \psi) (u_i u_j)(\cdot,y)dy\right|,\label{Formula_intermediairevVV801}
\end{eqnarray}
and the previous expression is nothing but the convolution between the function $(\partial_i \partial_j \psi) (u_i u_j)$ and a $L^\infty$-function, thus we have 
\begin{equation}\label{Formula_intermediairevVV81}
\left\| \phi \frac{\vn}{(-\Delta )}(\partial_i \partial_j \psi) (u_i u_j) (t,\cdot)\right\|_{L^\infty}\leq C\|(\partial_i \partial_j \psi) (u_i u_j)(t,\cdot)\|_{L^1}\leq C\|\mathds{1}_{Q_{R_1}}(u_i u_j)(t,\cdot)\|_{L^{\nu}},
\end{equation}
and taking the $L^\nu$-norm in the time variable we obtain
\begin{eqnarray*}
\left\| \phi \frac{\vn}{(-\Delta )}(\partial_i \partial_j \psi) (u_i u_j)\right\|_{L^{\nu}_t L^{\infty}_{x}}
&\leq &C\|\mathds{1}_{Q_{R_1}}u_i u_j\|_{L^{\nu}_{t,x}}\leq C\|\mathds{1}_{Q_{R_1}}\vu\|_{\mathcal{M}^{3,\tau_0}_{t,x}}\|\mathds{1}_{Q_{R_1}}\vu\|_{\mathcal{M}^{3,\tau_0}_{t,x}}<+\infty,
\end{eqnarray*}
where we used the fact that $1<\nu<\frac{3}{2}<\frac{\tau_0}{2}$ and we applied Hölder's inequality. Gathering together all these estimates we obtain 
$\|\mathds{1}_{Q_{R_2}}\vec{\mathbb{V}}_8\|_{\mathcal{M}^{3, \sigma}_{t,x}}<+\infty$.\\
\item For the quantity $\vec{\mathbb{V}}_{9}$ based in the term (9) of (\ref{FormulePourN}) we have
$$\left\|\mathds{1}_{Q_{R_2}}\vec{\mathbb{V}}_9\right\|_{\mathcal{M}^{3, \sigma}_{t,x}}=\left\|\mathds{1}_{Q_{R_2}}\int_{0}^{t} e^{(t-s)\Delta}\phi\frac{\vn}{(-\Delta)}\partial_{x_3}\left(\psi\theta \right)(s,\cdot)ds\right\|_{\mathcal{M}^{3, \sigma}_{t,x}},$$
and following the same ideas as before we can write
\begin{eqnarray*}
\left\|\mathds{1}_{Q_{R_2}}\vec{\mathbb{V}}_9\right\|_{\mathcal{M}^{3, \sigma}_{t,x}}&\leq& \left\|\mathds{1}_{Q_{R_2}}\mathcal{L}_2\left(\frac{\phi\vn \partial_{x_3}}{(-\Delta)}\left(\psi\theta \right)\right)\right\|_{\mathcal{M}^{3, \sigma}_{t,x}}\leq \left\|\mathcal{L}_2\left(\frac{\phi\vn \partial_{x_3}}{(-\Delta)}\left(\psi\theta \right)\right)\right\|_{\mathcal{M}^{310, 310}_{t,x}},
\end{eqnarray*}
if we set $q=\frac{62}{25}$, we have $2<\frac{q}{5}=\frac{125}{62}$, $\lambda=1-\frac{2q}{5}=\frac{1}{125}$ and $\frac{q}{\lambda}=310$, so we can write
\begin{eqnarray}
\left\|\mathcal{L}_2\left(\frac{\phi\vn \partial_{x_3}}{(-\Delta)}\left(\psi\theta \right)\right)\right\|_{\mathcal{M}^{310, 310}_{t,x}}&=&\left\|\mathcal{L}_2\left(\frac{\phi\vn \partial_{x_3}}{(-\Delta)}\left(\psi\theta \right)\right)\right\|_{\mathcal{M}^{\frac{q}{\lambda}, \frac{q}{\lambda}}_{t,x}}\notag\\
&\leq &C\left\|\frac{\phi\vn \partial_{x_3}}{(-\Delta)}\left(\psi\theta \right)\right\|_{\mathcal{M}^{q, q}_{t,x}}\leq C\left\|\psi\theta \right\|_{\mathcal{M}^{q, q}_{t,x}}\label{Valeur_Max_Iteration}
\end{eqnarray}
where we used the Lemma \ref{Lemma_Hed} as well as the fact that Riesz transforms are bounded in Lebesgue spaces. Since $q=\frac{62}{25}<\frac{10}{3}$ (and due to the support properties of the function $\psi$) we obtain (recall (\ref{Information_Interpolation_Theta})):
$$\left\|\psi\theta \right\|_{\mathcal{M}^{q, q}_{t,x}}\leq C\left\|\psi\theta \right\|_{\mathcal{M}^{\frac{10}{3}, \frac{10}{3}}_{t,x}}=C\left\|\psi\theta \right\|_{L^{\frac{10}{3}}_{t,x}}\leq C\|\theta\|_{L^{\frac{10}{3}}_{t,x}}<+\infty.$$

\item The quantity $\vec{\mathbb{V}}_{10}$ based in the term (10) of (\ref{FormulePourN}) and by the same arguments displayed to deduce (\ref{Formula_intermediairevVV80}), we can write (recall that $1<\nu<\frac{3}{2}$): 
$$\|\mathds{1}_{Q_{R_2}}\vec{\mathbb{V}}_{10}\|_{\mathcal{M}^{3, \sigma}_{t,x}}\leq C\left\|\phi\frac{\vn}{(-\Delta)}\left((\partial_{x_3}\psi)\theta \right) \right\|_{L^{\nu}_t L^{\infty}_{x}}.$$ 
If we study the $L^\infty$-norm in the space variable of this term, by the same ideas used in (\ref{Formula_intermediairevVV801})-(\ref{Formula_intermediairevVV81}) we obtain $\left\| \phi \frac{\vn}{(-\Delta )}\left((\partial_{x_3}\psi)\theta \right)(t,\cdot)\right\|_{L^{\infty}}\leq C\|(\partial_{x_3}\psi)\theta (t,\cdot)\|_{L^1}\leq C\|\mathds{1}_{Q_{R_1}}\theta (t,\cdot)\|_{L^{\nu}}$. Thus, taking the $L^{\nu}$-norm in the time variable we have
$$\|\mathds{1}_{Q_{R_2}}\vec{\mathbb{V}}_{10}\|_{\mathcal{M}^{3, \sigma}_{t,x}}
\leq C\left\| \phi \frac{\vn}{(-\Delta )}\left((\partial_{x_3}\psi)\theta \right)\right\|_{L^{\nu}_t L^{\infty}_{x}}\leq C \|\mathds{1}_{Q_{R_1}}\theta\|_{L^{\nu}_{t,x}}\leq C \|\mathds{1}_{Q_{R_1}}\theta\|_{L^{\infty}_{t}L^{2}_x}<+\infty.$$
\item The quantity $\vec{\mathbb{V}}_{11}$, givin in the term (11) of (\ref{FormulePourN}), can be treated in a similar manner. Indeed, by the same arguments displayed to deduce (\ref{Formula_intermediairevVV80}), we can write (recall that $1<\nu<\frac{3}{2}$): 
$$\|\mathds{1}_{Q_{R_2}}\vec{\mathbb{V}}_{9}\|_{\mathcal{M}^{3, \sigma}_{t,x}}\leq C\left\| \phi \frac{\vn}{(-\Delta )}( (\Delta \psi) p)\right\|_{L^{\nu}_t L^{\infty}_{x}},$$ 
and if we study the $L^\infty$-norm in the space variable of this term, by the same ideas used in (\ref{Formula_intermediairevVV801})-(\ref{Formula_intermediairevVV81}) we obtain $\left\| \phi \frac{\vn}{(-\Delta )}( (\Delta \psi) p) (t,\cdot)\right\|_{L^{\infty}}\leq C\|(\Delta \psi) p (t,\cdot)\|_{L^1}\leq C\|\mathds{1}_{Q_{R_1}}p (t,\cdot)\|_{L^{\nu}}$. Thus, taking the $L^{\nu}$-norm in the time variable we have
$$\|\mathds{1}_{Q_{R_2}}\vec{\mathbb{V}}_{9}\|_{\mathcal{M}^{3, \sigma}_{t,x}}
\leq C\left\| \phi \frac{\vn}{(-\Delta )}( (\Delta \psi) p)\right\|_{L^{\nu}_t L^{\infty}_{x}}\leq C \|\mathds{1}_{Q_{R_1}}p\|_{L^{\nu}_{t,x}}\leq C \|\mathds{1}_{Q_{R_1}}p\|_{L^{\frac32}_{t,x}}<+\infty.$$
\item The study of the quantity $\vec{\mathbb{V}}_{12}$ follows almost the same lines as the terms $\vec{\mathbb{V}}_{8}$ and $\vec{\mathbb{V}}_{11}$. However instead of (\ref{Formula_intermediairevVV801}) we have
$$\left|\phi \frac{\vn\partial_i }{(-\Delta)} ( (\partial_i \psi)p)\right|
\leq  C\left|\int_{\mathbb{R}^3} \frac{\mathds{1}_{|x-y|>(\mathfrak{R}_a-\mathfrak{R}_b)}}{|x-y|^3}\mathds{1}_{Q_{\mathfrak{R}_b}}(x)\mathds{1}_{Q_{R_1}\setminus Q_{\mathfrak{R}_a}}(y)(\partial_i  \psi) p(t,y)dy\right|,$$
and thus we can write:
$$\|\mathds{1}_{Q_{R_2}}\vec{\mathbb{V}}_{10}\|_{\mathcal{M}^{3, \sigma}_{t,x}}\leq  \left\|\phi \frac{\vn\partial_i }{(-\Delta)} ( (\partial_i \psi) p)\right\|_{L^{\nu}_t L^{\infty}_{x}}\leq C \|\mathds{1}_{Q_{R_1}}p\|_{L^{\nu}_{t,x}}\leq C \|\mathds{1}_{Q_{R_1}}p\|_{L^{\frac32}_{t,x}}<+\infty.$$
\item Finally, for the term $\vec{\mathbb{V}}_{13}$ based in the term (13) of (\ref{FormulePourN}) we write:
\begin{eqnarray*}
|\mathds{1}_{Q_{R_2}}\vec{\mathbb{V}}_{13}|=\left|\mathds{1}_{Q_{R_2}}\int_{0}^{t} e^{(t-s)\Delta}\phi(\theta e_3)(s,x)ds\right|.
\end{eqnarray*}
Using the same ideas as in  (\ref{Estimation_Ponctuelle_Riesz0})-(\ref{Estimation_Ponctuelle_Riesz1}) and applying again the Lemma \ref{lemma_locindi},  we obtain
\begin{eqnarray}
\|\mathds{1}_{Q_{R_2}}\vec{\mathbb{V}}_{13}\|_{\mathcal{M}^{3, \sigma}_{t,x}}&\leq &C\|\mathds{1}_{Q_{R_2}}(\mathcal{L}_2(|\mathds{1}_{Q_{\mathfrak{R}_b}} \theta |))\|_{\mathcal{M}^{3, \sigma}_{t,x}}\notag\\
&\leq &C\|\mathds{1}_{Q_{R_2}}(\mathcal{L}_2(|\mathds{1}_{Q_{\mathfrak{R}_b}} \theta |))\|_{\mathcal{M}^{310, 310}_{t,x}}\label{FormuleMaxSigma}\\
&\leq &C\|\mathcal{L}_2(|\mathds{1}_{Q_{\mathfrak{R}_b}} \theta |)\|_{\mathcal{M}^{\frac{q}{\lambda}, \frac{q}{\lambda}}_{t,x}},\notag
\end{eqnarray}
where we used the same parameters as in (\ref{Valeur_Max_Iteration}) and we can write (recall (\ref{Information_Interpolation_Theta})):
\begin{eqnarray*}
\|\mathcal{L}_2(|\mathds{1}_{Q_{\mathfrak{R}_b}} \theta |)\|_{\mathcal{M}^{\frac{q}{\lambda}, \frac{q}{\lambda}}_{t,x}}&\leq &C\|\mathds{1}_{Q_{\mathfrak{R}_b}} \theta \|_{\mathcal{M}^{q, q}_{t,x}}= C\|\mathds{1}_{Q_{\mathfrak{R}_b}} \theta \|_{\mathcal{M}^{\frac{10}{3}, \frac{10}{3}}_{t,x}}=\|\mathds{1}_{Q_{\mathfrak{R}_b}} \theta \|_{L^{\frac{10}{3}}_{t,x}}<+\infty.
\end{eqnarray*}
We can thus conclude that 
$$\|\mathds{1}_{Q_{R_2}}\vec{\mathbb{V}}_{13}\|_{\mathcal{M}^{3, \sigma}_{t,x}}<+\infty.$$
\end{itemize}
With all these estimates Proposition \ref{Propo_GainenMorreySigma} is now proven. \hfill$\blacksquare$
\begin{Remark}\label{Rem_LimitIteration_vu}
Note that the value of the index $\sigma$ of the Morrey space $\mathcal{M}_{t,x}^{3, \sigma}(\mathbb{R}\times \R)$ is potentially bounded by the information available over the variable $\theta$  (recall (\ref{Information_Interpolation_Theta})) and the maximal possible value for this parameter is close to $\sigma=310$ (see the expression (\ref{FormuleMaxSigma}) above).
\end{Remark}
This result gives a small gain of integrability as we pass from an information on the Morrey space $\mathcal{M}_{t,x}^{3, \tau_0}$ to a control over the space $\mathcal{M}_{t,x}^{3, \sigma}$ with $\tau_0<\sigma$ with $\sigma$ close to $\tau_0$. This is of course not enough and we need to repeat the arguments above in order to obtain a better control. In this sense we have the following proposition:
\begin{Proposition}\label{MorreyPropo2}
Under the hypotheses of Theorem \ref{Teo_HolderRegularity} and within the framework of Proposition \ref{Propo_FirstMorreySpace}, there exists a radius $\bar R_2$ with $0<\bar R_2<R_2$ such that
\begin{equation}\label{Conclusion_SecondMorreySpace}
\mathds{1}_{Q_{\bar R_2}(t_0,x_0)}\vu\in \M_{t,x}^{3,310}(\mathbb{R}\times \R),
\end{equation}
\end{Proposition}
\noindent \textbf{Proof.}
By the Proposition \ref{Propo_GainenMorreySigma} above it follows that $\mathds{1}_{Q_{R_2}} \vu \in \mathcal{M}_{t,x}^{3,\sigma}(\mathbb{R}\times \R)$ with $\sigma$ very close to $\tau_0$ (say $\sigma=\tau_0+\epsilon$). Hence, with the information $\mathds{1}_{Q_{R_2}} \vu \in \mathcal{M}_{t,x}^{3, \tau_0+\epsilon}(\mathbb{R}\times \R)$ at hand, we can reapply the Proposition \ref{Propo_GainenMorreySigma} to obtain for some smaller radius $\bar {R}_2< R_2$ that $\mathds{1}_{Q_{\bar{R}_2}} \vu \in \mathcal{M}_{t,x}^{3, \sigma_1}(\mathbb{R}\times \R)$ where $\sigma_1=\sigma+\epsilon=\tau_0+2\epsilon$. Iterating these arguments as long as necessary, we obtain the information $\mathds{1}_{Q_{R_2}} \vu \in \mathcal{M}_{t,x}^{3, 60}(\mathbb{R}\times \R)$ where the value $\sigma=60$ is fixed by the information available for the quantity $\vw$ which is the only term that is fixed: see the computation leading to the estimate (\ref{FormuleMaxSigma}) and Remark \ref{Rem_LimitIteration_vu}. Let us note that a slight abuse of language has been used for the radius $\bar R_2$: at each iteration this radius is smaller and smaller, but in order to maintain the notations we still denote the final radius by $\bar R_2$. \hfill $\blacksquare$

\section{A  gain of information for the variable $\theta$}\label{Secc_GainTheta}

Note that the Proposition \ref{MorreyPropo2} and the Corollary \ref{corolarioMorrey} give interesting control (on a small neighborhood of a point $(t_0, x_0)$) for the variable $\vu$. Indeed, we have obtained so far the information
\begin{equation}\label{Information_PremiereEtape}
\mathds{1}_{Q_{\bar R_2}(t_0,x_0)} \vu\in \mathcal{M}_{t,x}^{3,310}(\mathbb{R}\times \R),\qquad \mbox{and} \qquad
\mathds{1}_{Q_{R_1}(t_0,x_0)}\vn\otimes \vu\in \mathcal{M}^{2,\tau_1}_{t,x}(\mathbb{R}\times \R),
\end{equation}
where 
\begin{equation}\label{Info_Radios}
0<\bar R_2<R_1<R<1,
\end{equation}
with $\frac{5}{1-\alpha}<\tau_0<6$ and $\tau_1$ is given by the condition $\frac{1}{\tau_1}=\frac{1}{\tau_0}+\frac{1}{5}$ (see the Corollary \ref{corolarioMorrey}).\\

Now we will exploit all this information in order to derive some Morrey control for the variable $\theta$, indeed, we have:
\begin{Proposition}\label{Propo_Information_Morrey_Theta}
Under the general hypotheses of Theorem \ref{Teo_HolderRegularity}, if we have the controls (\ref{Information_PremiereEtape}) over $\vu$ then we have, for some radius $0<R_3<\bar R_2$, we have
$$\mathds{1}_{Q_{R_3}(t_0,x_0)}\theta\in\M_{t,x}^{\frac{12}{5},\frac{50}{9}}(\mathbb{R}\times\R).$$ 
\end{Proposition}
\textbf{Proof.} Consider  $\bar \phi: \mathbb{R}\times\mathbb{R}^3\longrightarrow \mathbb{R}$ a non-negative function such that $\bar \phi\in \mathcal{C}_0^{\infty}(\mathbb{R}\times \mathbb{R}^3)$ and such that
\begin{equation}\label{Def_FuncAuxiliarPourDivW}
\bar \phi \equiv 1\;\; \text{over}\; \; Q_{\rho_b}(t_0,x_0),\; \; supp(\bar \phi)\subset Q_{\rho_a}(t_0,x_0),
\end{equation}
where we have
\begin{equation}\label{Info_Radios1}
0<R_3<\rho_b<\rho_a<\bar R_2, 
\end{equation} 
where the radius $\bar R_2$ is fixed in (\ref{Info_Radios}).  With the help of this auxiliar function we define the variable $\Theta$ by 
$$\Theta=\bar \phi \theta,$$ 
note that, due to the support properties of the function $\bar \phi$ we have $\mathds{1}_{Q_{R_3}}\Theta=\mathds{1}_{Q_{R_3}}\theta$. \\

Thus, if we study the evolution of $\Theta$ we obtain (by (\ref{Equation_Boussinesq_intro}))
\begin{eqnarray}
\partial_t\Theta&=&(\partial_t \bar \phi )\theta+ \bar \phi \bigg(\Delta \theta -\vu\cdot \vn \theta\bigg)\notag\\
&=&\Delta \Theta + (\partial_t \bar \phi +\Delta \bar\phi)\theta-2\sum_{i=1}^3  \partial_i\big((\partial_i \bar \phi) \theta\big)-div(\bar\phi \vu\theta)+(\vn \bar\phi)\cdot(\vu\theta),\label{Equation_localisee_Theta}
\end{eqnarray}
where we used the identities $\displaystyle{\bar \phi \Delta \theta=\Delta(\bar \phi \theta)+\Delta \bar \phi \theta-2\sum_{i=1}^3  \partial_i\big((\partial_i \bar \phi) \theta\big)}$ and the fact that
$$\bar\phi(\vu\cdot\vn\theta)=\bar\phi div(\vu\theta)=div(\bar\phi\vu\theta)-\vn\bar\phi\cdot(\vu\theta),$$
since $div(\vu)=0$. As we have $\Theta(0, \cdot)=0$ (by the properties of the localizing function $\bar\phi$ given in (\ref{Def_FuncAuxiliarPourDivW})), applying the Duhamel formula we can write:
\begin{eqnarray}
\Theta(t,x)&=&\underbrace{\int_0^t e^{(t-s)\Delta}(\partial_t \bar \phi+\Delta \bar\phi)\theta ds}_{\Theta_1}-2\sum_{i=1}^3\underbrace{\int_0^t e^{(t-s)\Delta}\partial_i \big((\partial_i \bar\phi)\theta \big)ds}_{\Theta_2}\label{Ws}-\underbrace{\int_0^t e^{(t-s)\Delta} div\big(\bar\phi\vu \theta\big)ds}_{\Theta_3}\\
&&+\underbrace{\int_0^t e^{(t-s)\Delta} (\vn \bar\phi)\cdot(\vu\theta) ds}_{\Theta_4},\notag
\end{eqnarray}
and we will estimate each one of the terms above to prove that we have the wished Morrey control over the set $Q_{R_3}$.
\begin{itemize}
\item For the first term $\Theta_1$ we write, 
\begin{eqnarray}
|\mathds{1}_{Q_{R_3}}\Theta_1|&=&\left|\mathds{1}_{Q_{R_3}}\int_0^t e^{(t-s)\Delta}(\partial_t \bar \phi+\Delta \bar\phi)\theta ds\right|\label{Estimation_W1}
\end{eqnarray}
by the decay properties of the heat kernel, by the properties of the test function $\bar\phi$ (see (\ref{Def_FuncAuxiliarPourDivW})) and by the definition of the parabolic Riesz potential $\mathcal{L}_2$ given in (\ref{Def_ParabolicRieszPotential}), we can write the estimate
 \begin{eqnarray}
|\mathds{1}_{Q_{R_3}}\Theta_1|&\leq &C\mathds{1}_{Q_{R_3}}\int_{\mathbb{R}}\int_{\mathbb{R}^3} \frac{|\mathds{1}_{Q_{\rho_a}} \theta(s,y)|}{(|t-s|^{\frac{1}{2}}+|x-y|)^3}dyds\notag\\
&&\leq C\mathds{1}_{Q_{R_3}}(\mathcal{L}_2(|\mathds{1}_{Q_{\rho_a}} \theta |))(t,x).\label{Estimation_W11}
\end{eqnarray}
If we fix $p=\frac{12}{5}$, $q=\frac{62}{25}$ and $\lambda=\frac{1}{125}$, by applying Lemma \ref{lemma_locindi} and  Lemma \ref{Lemma_Hed} we obtain
\begin{eqnarray*}
\left\|\mathds{1}_{Q_{R_3}}(\mathcal{L}_2(|\mathds{1}_{Q_{\rho_a}} \theta |))\right\|_{\M_{t,x}^{\frac{12}{5},\frac{50}{9}}}&\leq &C\|\mathcal{L}_2(|\mathds{1}_{Q_{\rho_a}}\theta|)\|_{\M_{t,x}^{\frac{12}{\lambda 5},\frac{q}{\lambda }}}\\
&\leq &C\|\mathds{1}_{Q_{\rho_a}}\theta \|_{\M_{t,x}^{\frac{12}{5},\frac{62}{25}}}\leq C\|\mathds{1}_{Q_{\rho_a}}\theta \|_{\M_{t,x}^{\frac{10}{3},\frac{10}{3}}}\leq C \|\mathds{1}_{Q_{\rho_a}}\theta \|_{L_{t,x} ^{\frac{10}{3}}}<+\infty,
\end{eqnarray*}
where we used the information (\ref{Information_Interpolation_Theta}).\\
\item For the term $\Theta_2$ of (\ref{Ws}) we need to study, for all $1\le i\le 3$, the quantities
$$\mathbb{O}_i=\left|\mathds{1}_{Q_{R_3}}\int_0^t e^{(t-s)\Delta}\partial_i \big((\partial_i \bar\phi)\theta\big)ds\right|.$$
Thus, using the definition of the parabolic Riesz potential $\mathcal{L}_1$ given in (\ref{Def_ParabolicRieszPotential}) and applying the Lemma \ref{lemma_locindi}  we can write
$$\left\|\mathbb{O}_i\right\|_{\M_{t,x}^{\frac{12}{5},\frac{50}{9}}}\leq C\left\|\mathds{1}_{Q_{R_3}}\mathcal{L}_1 \big((\partial_i \bar\phi)\theta\big)\right\|_{\M_{t,x}^{\frac{12}{5},\frac{50}{9}}}\leq C\left\|\mathcal{L}_1 \big((\partial_i \bar\phi)\theta\big)\right\|_{\M_{t,x}^{\frac{50}{9},\frac{50}{9}}}.$$
If we set now $q=\frac{10}{3}$, $0<1<\frac{5}{q}$, $\lambda=1-\frac{q}{5}=\frac{3}{5}$, we have $\frac{q}{\lambda}=\frac{50}{9}$, we can use the Lemma \ref{Lemma_Hed} to obtain
\begin{eqnarray*}
\left\|\mathbb{O}_i\right\|_{\M_{t,x}^{\frac{12}{5},\frac{50}{9}}}&\leq &C\left\|\mathcal{L}_1 \big((\partial_i \bar\phi)\theta\big)\right\|_{\M_{t,x}^{\frac{q}{\lambda},\frac{q}{\lambda}}}\\
&\leq &C\left\|(\partial_i \bar\phi)\theta\right\|_{\M_{t,x}^{q,q}}= C\left\|(\partial_i \bar\phi)\theta\right\|_{\M_{t,x}^{\frac{10}{3},\frac{10}{3}}}\leq C\|\theta\|_{L_{t,x}^{\frac{10}{3}}}<+\infty.
\end{eqnarray*}
Note that the value $\frac{50}{9}$ is related to the information available over $\theta$ stated in (\ref{Information_Interpolation_Theta}). With these estimates at hand we deduce that 
$$\left\|\Theta_2\right\|_{\M_{t,x}^{\frac{12}{5},\frac{50}{9}}}<+\infty.$$
\item We study now the term $\Theta_3$ defined in (\ref{Ws}) and we write (using the same arguments as before)
\begin{eqnarray*}
\|\mathds{1}_{Q_{R_3}}\Theta_3\|_{\M_{t,x}^{\frac{12}{5}, \frac{50}{9}}}&= &\left\|\mathds{1}_{Q_{R_3}}\int_0^t e^{(t-s)\Delta} div\big(\bar\phi\vu \theta\big)ds\right\|_{\M_{t,x}^{\frac{12}{5}, \frac{50}{9}}}\leq  C\left\|\mathds{1}_{Q_{R_3}}\mathcal{L}_1\big(\bar\phi\vu \theta\big)\right\|_{\M_{t,x}^{\frac{12}{5}, \frac{50}{9}}}.
\end{eqnarray*}
We set $q=\frac{310}{94}$, $\lambda=1-\frac{q}{5}=\frac{16}{47}$ and $p=\frac{30}{19}$. Since $\frac{12}{5}<\frac{p}{\lambda}$ and $\frac{50}{9}<\frac{q}{\lambda}$ we obtain (by the Lemma \ref{Lemma_Hed})
\begin{eqnarray*}
\left\|\mathds{1}_{Q_{R_3}}\mathcal{L}_1\big(\bar\phi\vu \theta\big)\right\|_{\M_{t,x}^{\frac{12}{5}, \frac{50}{9}}}&\leq & C\left\|\mathcal{L}_1\big(\bar\phi\vu \theta\big)\right\|_{\M_{t,x}^{\frac{p}{\lambda}, \frac{q}{\lambda}}} \leq  C\left\| \bar\phi\vu \theta\right\|_{\M_{t,x}^{\frac{30}{19}, \frac{310}{94}}}\\
&\leq &C\|\mathds{1}_Q\vu\|_{\M_{t,x}^{3, 310}}\|\mathds{1}_Q\theta\|_{\M_{t,x}^{\frac{10}{3},\frac{10}{3}}}<+\infty,
\end{eqnarray*}
where in the last estimate we used the H\"older inequalities in Morrey spaces with $\frac{19}{30}=\frac{3}{10}+\frac{1}{3}$ and $\frac{94}{310}=\frac{3}{10}+\frac{1}{310}$. We conclude that 
$$\|\mathds{1}_{Q_{R_3}}\Theta_3\|_{\M_{t,x}^{\frac{12}{5}, \frac{50}{9}}}<+\infty.$$
\item For the term $\Theta_4$ given in (\ref{Ws}), by the same arguments as before we obtain the estimates
$$\|\mathds{1}_{Q_{R_3}}\Theta_4\|_{\M_{t,x}^{\frac{12}{5},\frac{50}{9}}}\leq C\left\|\mathds{1}_{Q_{R_3}}(\mathcal{L}_2(|\mathds{1}_{Q_{\rho_a}} \vu\theta |))\right\|_{\M_{t,x}^{\frac{12}{5},\frac{50}{9}}}.$$
Setting $p=\frac{30}{19}$, $q=\frac{62}{25}$ and $\lambda=\frac{1}{125}$, since $\frac{12}{5}<\frac{p}{\lambda}$ and $\frac{50}{9}<\frac{q}{\lambda}$, applying Lemma \ref{Lemma_Hed} and Lemma \ref{lemma_locindi} we have
\begin{eqnarray*}
\left\|\mathds{1}_{Q_{R_3}}(\mathcal{L}_2(|\mathds{1}_{Q_{\rho_a}} \vu\theta |))\right\|_{\M_{t,x}^{\frac{12}{5},\frac{50}{9}}}&\leq &\left\|\mathcal{L}_2(|\mathds{1}_{Q_{\rho_a}} \vu\theta |)\right\|_{\M_{t,x}^{\frac{p}{\lambda},\frac{q}{\lambda}}}\leq  C\left\|\mathds{1}_{Q_{\rho_a}} \vu\theta \right\|_{\M_{t,x}^{\frac{30}{19},\frac{62}{25}}}\\
&\leq & C\left\|\mathds{1}_{Q_{\rho_a}} \vu\theta \right\|_{\M_{t,x}^{\frac{30}{19},\frac{310}{94}}}\leq C\|\mathds{1}_Q\vu\|_{\M_{t,x}^{3, 310}}\|\mathds{1}_Q\theta\|_{\M_{t,x}^{\frac{10}{3},\frac{10}{3}}}<+\infty,
\end{eqnarray*}
where we used the H\"older inequalities with $\frac{19}{30}=\frac{3}{10}+\frac{1}{3}$ and $\frac{94}{310}=\frac{3}{10}+\frac{1}{310}$. We finally obtain
$$\|\mathds{1}_{Q_{R_3}}\Theta_4\|_{\M_{t,x}^{\frac{12}{5}, \frac{50}{9}}}<+\infty.$$
\end{itemize}
With all these controls, Proposition \ref{Propo_Information_Morrey_Theta} is proven. \hfill$\blacksquare$
\section{The end of the proof of Theorem \ref{Teo_HolderRegularity}}\label{Secc_Final}
The key result for obtaining a gain of regularity is the following lemma coming from the theory of parabolic equations (see \cite{Lady68,PGLR1}).
\begin{Lemma}\label{Lemma_parabolicHolder}
Let $\sigma$  be a smooth homogeneous function over $\R\setminus \{0\}$,  of exponent 1 with 
$\sigma(D)$ the Fourier multiplier associated.
Consider the functions $ \Phi \in \mathcal{M}_{t,x}^{\p_0,\q_0}(\mathbb{R}\times\R)$ and $h\in \mathcal{M}_{t,x}^{\p_0,\q_1}(\mathbb{R}\times\R)$ such that $1\le \p_0\le \q_0$, with  $\frac{1}{\q_0}=\frac{2-\alpha}{5}$,$\frac{1}{\q_1}=\frac{1-\alpha}{5}$, for $0<\alpha<1$. Then, the function $v$ equal to $0$ for $t\le 0$ and
\begin{equation*}
v(t,x)=\int_0^t e^{(t-s)\Delta }(\Phi (s,\cdot)+\sigma(D) h(s,\cdot))ds,
\end{equation*}
for $t>0$,  is H\"older continuous  of exponent $\alpha$ with respect to the parabolic distance.  
\end{Lemma} 
In order to apply this lemma to the proof of Theorem \ref{Teo_HolderRegularity}, we will first localize each one of the equations of the Boussinesq system (\ref{Equation_Boussinesq_intro}) and then we will show the terms of the corresponding Duhamel formula belongs either to the space $\mathcal{M}_{t,x}^{\p_0,\q_0}(\mathbb{R}\times\R)$ or to the space $ \mathcal{M}_{t,x}^{\p_0,\q_1}(\mathbb{R}\times\R)$.\\

We start by localizing the problem and for this we consider  $\phi: \mathbb{R}\times \R\longrightarrow \mathbb{R}$ a test function such that $supp(\phi)\subset ]-\frac{1}{4},\frac{1}{4}[\times B(0,\frac{1}{2})$ and $\phi \equiv 1\;\; \text{over}\; \; ]-\frac{1}{16},\frac{1}{16}[\times B(0,\frac{1}{4})$. We consider next a radius $\mathbf{R}>0$ such that 
\begin{equation}\label{Info_Radios2}
4\mathbf{R}<R_3<\bar R_2<R_1<R<1,
\end{equation}
where $R_3$ is the radius of Proposition \ref{Propo_Information_Morrey_Theta}, $\bar R_2$ is the radius of Proposition \ref{MorreyPropo2} and $R_1$ is the radius obtained in Proposition \ref{Propo_FirstMorreySpace}. We then define
\begin{equation}\label{Def_Eta}
\eta(t,x)=\phi\left(\frac{t-t_0}{\mathbf{R}^2},\frac{x-x_0}{\mathbf{R}}\right).
\end{equation}
We start with the equation over the velocity field and we consider the variable $\vUU$ defined by the formula
\begin{equation}\label{DefU}
\vUU=\eta\vu,
\end{equation}
then, by the properties of the auxiliar function $\eta$, we have the identity $\vUU=\vu$ over a small neighborhood of the point $(t_0,x_0)$, the support of the variable $\vUU$ is contained in the parabolic ball $Q_{\mathbf{R}}$ and moreover we also have $\vUU(0,x)=0$. Thus, if we study the evolution of this variable, following the system (\ref{Equation_Boussinesq_intro}), we have
\begin{eqnarray*}
\partial_t\vUU&=&(\partial_t\eta)\vu+\eta \Delta \vu-\eta ((\vu\cdot \vn)\vu)-\eta\vn p+\eta\theta e_3.
\end{eqnarray*}
We use now the identity $\displaystyle{\eta \Delta \vu=\Delta\vUU +(\Delta \eta)\vu-2\sum_{i=1}^3 \partial_i((\partial_i \eta) (\vu))}$ to obtain the equation
$$\partial_t\vUU=\Delta\vUU+(\partial_t\eta+\Delta \eta)(\vu+\vw)-2\sum_{i=1}^3\partial_i((\partial_i \eta) (\vu))-\eta ((\vu\cdot \vn)\vu)-\eta\vn p+\eta\theta e_3.$$
Noting that we have (since $div(\vu)=0$) the identity
$$\eta((\vu\cdot\vn)\vu)=\eta div(\vu\otimes \vu)=div(\eta \vu\otimes \vu)-\vu\otimes \vu\cdot\vn\eta,$$
we finally deduce the following equation for the variable $\vUU$:
\begin{equation}\label{Equation_Evolution_Final_vUU}
\partial_t\vUU=\Delta\vUU+(\partial_t\eta+\Delta \eta)\vu-2\sum_{i=1}^3\partial_i((\partial_i \eta) (\vu))-div(\eta \vu\otimes \vu)+\vu\otimes \vu\cdot\vn\eta-\eta\vn p+\eta\theta e_3.
\end{equation}
We consider now the equation for the variable $\theta$ and we consider the variable $\mathcal{O}$ defined by the formula
\begin{equation}\label{DefU}
\mathcal{O}=\eta\theta.
\end{equation}
Thus, following the same ideas used to deduce the equation (\ref{Equation_localisee_Theta}) we can write
\begin{equation}\label{Equation_Evolution_Final_Theta}
\partial_t\mathcal{O}=\Delta \mathcal{O} + (\partial_t \eta +\Delta \eta)\theta-2\sum_{i=1}^3  \partial_i\big((\partial_i \eta) \theta\big)-div(\eta \vu\theta)+(\vn \eta)\cdot(\vu\theta).
\end{equation}
At this point, we define the $(3+1)$D vector $\vV=\displaystyle{\left[\begin{matrix}
\vUU\\[1mm]
\mathcal{O}
\end{matrix}\right]}=\displaystyle{\left[\begin{matrix}
\eta\vu\\[1mm]
\eta\theta
\end{matrix}\right]}$, thus with the equations (\ref{Equation_Evolution_Final_vUU}) and (\ref{Equation_Evolution_Final_Theta}) we obtain the system
\begin{equation}\label{Equation_Boussinesq_Vectorial}
\begin{cases}
\partial_t\vV=  \Delta\vV+\vA -2\displaystyle{\sum_{i=1}^3}\partial_i\vB-div(\mathbb{C}) + \vD - \vE +\vF,\\[2mm] 
\vV(0,x)=0.
\end{cases}
\end{equation}
where the $(3+1)D$ vectors $\vA,...,\vH$ are defined by 
\begin{equation}\label{Definition_QuantitiesMorreyHolder}
\begin{split}
&\vA=(\partial_t \eta +\Delta \eta)\left[\begin{matrix}
\vu\\[1mm]
\theta
\end{matrix}\right], \quad
\vB=\left[\begin{matrix}
(\partial_i \eta) \vu\\[1mm]
(\partial_i \eta) \theta
\end{matrix}\right],
\quad
\mathbb C=\left[\begin{matrix}
\eta \vu\otimes \vu\\[1mm]
\eta \vu\theta
\end{matrix}\right],\\[2mm]
&\vD=\left[\begin{matrix}
\vu\otimes \vu\cdot\vn\eta\\[1mm]
(\vn \eta)\cdot(\vu\theta)
\end{matrix}\right],\quad
\vE=\left[\begin{matrix}
\eta\vn p\\[1mm]
0
\end{matrix}\right],
\quad
\vF=\left[\begin{matrix}
\eta \theta e_3\\
0
\end{matrix}\right],
\end{split}
\end{equation}
note that the term $\mathbb{C}$ is not exactly a $(3+1)D$ vector (the first bloc is a tensor) and the quantity $div(\mathbb{C})$ must be understood in the following sense:  $div(\mathbb {C})=\left[\begin{matrix}
div(\eta \vu\otimes \vu)\\[1mm]
div(\eta \vu\theta)
\end{matrix}\right]$. This slight abuse of notation can be easily understood if we work component by component.\\

Thus, by the Duhamel formula, the solution of the equation (\ref{Equation_Boussinesq_Vectorial}) can be written in the following manner:
\begin{equation}\label{Duhamel_FullBoussinesq}
\vV(t,x)=\int_0^t e^{(t-s)\Delta }\bigg(\vA-2\displaystyle{\sum_{i=1}^3}\vB-div(\mathbb{C})+\vD-\vE+\vF\bigg)ds,
\end{equation}
thus, in order to apply the Lemma \ref{Lemma_parabolicHolder} to this system and obtain a parabolic gain of regularity, we only need to prove that the quantities $\vA,...,\vF$, defined in (\ref{Definition_QuantitiesMorreyHolder}) satisfy:

\begin{equation}\label{Objetivo_MorreyLadyz}
\vA, \vD, \vE, \vF\in \M_{t,x}^{\p_0,\q_0}(\mathbb{R}\times\R)\quad \mbox{and}\quad \vB, \mathbb{C} \in \M_{t,x}^{\p_0,\q_1}(\mathbb{R}\times\R),
\end{equation}
where we will assume $1\le \p_0\leq \frac{6}{5}\le \q_0$, with  $\frac{1}{\q_0}=\frac{2-\alpha}{5}$,$\frac{1}{\q_1}=\frac{1-\alpha}{5}$, for some $0<\alpha\ll \frac{1}{10}$.\\

\noindent To prove (\ref{Objetivo_MorreyLadyz}) we recall that we have the following controls
\begin{equation}\label{Recap_Controls}
\begin{split}
\mathds{1}_{Q_{\bar R_2}} \vu\in \mathcal{M}_{t,x}^{3,310}(\mathbb{R}\times \R),\qquad
\mathds{1}_{Q_{R_1}}\vn\otimes \vu\in \mathcal{M}^{2,\tau_1}_{t,x}(\mathbb{R}\times \R)\\
\mbox{and}\qquad \mathds{1}_{Q_{R_3}}\theta\in\M_{t,x}^{\frac{12}{5},\frac{50}{9}}(\mathbb{R}\times\R),\qquad \mathds{1}_{Q_{R}}\theta\in\M_{t,x}^{\frac{10}{3},\frac{10}{3}}(\mathbb{R}\times\R),
\end{split}
\end{equation}
where we have $4\mathbf{R}<R_3<\bar R_2< R_1<R$ and $\frac{5}{1-\alpha}<\tau_0<6$ and $\tau_1$ is given by the condition $\frac{1}{\tau_1}=\frac{1}{\tau_0}+\frac{1}{5}$.
\begin{itemize}
\item Let us start with the quantity $\vA$. We write, by Lemma \ref{lemma_locindi} and since $1\leq \p_0\leq \frac{6}{5}$: 
$$\|\vA\|_{\M_{t,x}^{\p_0,\q_0}}=\left\|(\partial_t \eta +\Delta \eta)\left[\begin{matrix}
\vu\\[1mm]
\theta
\end{matrix}\right]\right\|_{\M_{t,x}^{\p_0,\q_0}}\leq C\left\|\mathds{1}_{Q_\mathbf{R}}\left[\begin{matrix}
\vu\\[1mm]
\theta
\end{matrix}\right]\right\|_{\M_{t,x}^{\frac{6}{5},\q_0}} \leq C\|\mathds{1}_{Q_\mathbf{R}}\vu\|_{\M_{t,x}^{\frac{6}{5},\q_0}}+C\|\mathds{1}_{Q_\mathbf{R}}\theta\|_{\M_{t,x}^{\frac{6}{5},\q_0}},$$
noting that $\q_0=\frac{5}{2-\alpha}<\frac{50}{9}$ as $0<\alpha\ll\frac{1}{10}$, we thus have 
$$\|\vA\|_{\M_{t,x}^{\p_0,\q_0}}\leq C\|\mathds{1}_{Q_{\bar R_2}}\vu\|_{\M_{t,x}^{3,310}}+C\|\mathds{1}_{Q_{R_3}}\theta\|_{\M_{t,x}^{\frac{12}{5},\frac{50}{9}}}<+\infty,$$
since we have the controls (\ref{Recap_Controls}).
\item For the term $\vB$ given in (\ref{Definition_QuantitiesMorreyHolder}), we have 
$$\|\vB\|_{\M_{t,x}^{\p_0,\q_1}}=\left\|\left[\begin{matrix}
(\partial_i \eta) \vu\\[1mm]
(\partial_i \eta) \theta
\end{matrix}\right]\right\|_{\M_{t,x}^{\p_0,\q_1}}\leq C\|\mathds{1}_{Q_\mathbf{R}}\vu\|_{\M_{t,x}^{\p_0,\q_1}}+C\|\mathds{1}_{Q_\mathbf{R}}\theta\|_{\M_{t,x}^{\p_0,\q_1}},$$
since $1\leq \p_0\leq \frac{6}{5}$ and $\q_1=\frac{5}{1-\alpha}<\frac{50}{9}$ (as $0<\alpha\ll\frac{1}{10})$, by Lemma \ref{lemma_locindi} we have
$$\|\vB\|_{\M_{t,x}^{\p_0,\q_1}}\leq C\|\mathds{1}_{Q_{\bar R_2}}\vu\|_{\M_{t,x}^{3,310}}+C\|\mathds{1}_{Q_{R_3}}\theta\|_{\M_{t,x}^{\frac{12}{5},\frac{50}{9}}}<+\infty,$$
since we have the controls (\ref{Recap_Controls}).
\begin{Remark}
Note that in this particular stage the information $\theta\in L^{\frac{10}{3}}_{t,x}$ is not enough and we need the sharper control $\mathds{1}_{Q_{R_3}}\theta\in \M_{t,x}^{\frac{12}{5},\frac{50}{9}}$.
\end{Remark}
\item For the quantity $\mathbb{C}$ defined in  (\ref{Definition_QuantitiesMorreyHolder}) we have
\begin{equation}\label{Estimate_C}
\|\mathbb{C}\|_{\M_{t,x}^{\p_0,\q_1}}=\left\|\left[\begin{matrix}
\eta \vu\otimes \vu\\[1mm]
\eta \vu\theta
\end{matrix}\right]\right\|_{\M_{t,x}^{\p_0,\q_1}}\leq C\|\mathds{1}_{Q_\mathbf{R}} \vu\otimes \vu\|_{\M_{t,x}^{\p_0,\q_1}}+C\|\mathds{1}_{Q_\mathbf{R}}\vu\theta\|_{\M_{t,x}^{\p_0,\q_1}}.
\end{equation}
For the first term in the right-hand side above we write, since $\p_0\leq \frac{6}{5}<\frac{3}{2}$ and $\q_1=\frac{5}{1-\alpha}<6$ as $0<\alpha\ll1$: 
$$\|\mathds{1}_{Q_\mathbf{R}} \vu\otimes \vu\|_{\M_{t,x}^{\p_0,\q_1}}\leq C\|\mathds{1}_{Q_\mathbf{R}} \vu\otimes \vu\|_{\M_{t,x}^{\frac{3}{2},155}}\leq C\|\mathds{1}_{\bar Q_2} \vu\otimes \vu\|_{\M_{t,x}^{3,310}}\|\mathds{1}_{\bar Q_2} \vu\otimes \vu\|_{\M_{t,x}^{3,310}}<+\infty,$$
where we used the H\"older inequalities in Morrey spaces and the information given in (\ref{Recap_Controls}).\\

For the second term in the right-hand side of (\ref{Estimate_C}) we have (since $\p_0\leq \frac{6}{5}<\frac{4}{3}$ and $\q_1=\frac{5}{1-\alpha}<\frac{775}{142}$, which is possible if $0<\alpha \ll1$)
$$\|\mathds{1}_{Q_\mathbf{R}}\vu\theta\|_{\M_{t,x}^{\p_0,\q_1}}\leq \|\mathds{1}_{Q_\mathbf{R}}\vu\theta\|_{\M_{t,x}^{\frac{4}{3},\frac{775}{142}}}\leq C\|\mathds{1}_{Q_{\bar R_2}}\vu\|_{\M_{t,x}^{3,310}}\|\mathds{1}_{Q_{R_3}}\theta\|_{\M_{t,x}^{\frac{12}{5},\frac{50}{9}}}<+\infty,$$
where we used the H\"older inequalities in the Morrey space setting with $\frac{3}{4}=\frac{1}{3}+\frac{5}{12}$ and $\frac{142}{775}=\frac{1}{310}+\frac{9}{50}$.\\

We thus obtain $\|\mathbb{C}\|_{\M_{t,x}^{\p_0,\q_1}}<+\infty$.
\item The term $\vD$ given in (\ref{Definition_QuantitiesMorreyHolder}) we write
\begin{equation}\label{Estimate_D}
\|\vD\|_{\M_{t,x}^{\p_0,\q_0}} =\left\|\left[\begin{matrix}
\vu\otimes \vu\cdot\vn\eta\\[1mm]
(\vn \eta)\cdot(\vu\theta)
\end{matrix}\right]\right\|_{\M_{t,x}^{\p_0,\q_0}}\leq C\|\mathds{1}_{Q_\mathbf{R}}
\vu\otimes \vu\|_{\M_{t,x}^{\p_0,\q_0}}+C\|\mathds{1}_{Q_\mathbf{R}}\vu\theta\|_{\M_{t,x}^{\p_0,\q_0}}.\end{equation}
For the first term of the right-hand side above we write, by Lemma \ref{lemma_locindi}, since $1\leq \p_0\leq \frac{6}{5}\leq \frac{3}{2}$ and since $\q_0=\frac{5}{2-\alpha}<155$:
\begin{eqnarray*}
\|\mathds{1}_{Q_\mathbf{R}}\vu\otimes \vu\|_{\M_{t,x}^{\p_0,\q_0}}&\leq &C\|\mathds{1}_{Q_\mathbf{R}}\vu\otimes \vu\|_{\M_{t,x}^{\frac{3}{2},\q_0}}\leq C \|\mathds{1}_{Q_\mathbf{R}}
\vu\|_{\M_{t,x}^{3,2\q_0}}\|\mathds{1}_{Q_\mathbf{R}}\vu\|_{\M_{t,x}^{3,2\q_0}}\\
&\leq & C\|\mathds{1}_{Q_{\bar R_2}}\vu\|_{\M_{t,x}^{3,310}}\|\mathds{1}_{Q_{\bar R_2}}\vu\|_{\M_{t,x}^{3,310}}<+\infty,
\end{eqnarray*}
as we have the controls (\ref{Recap_Controls}).
For the second term of (\ref{Estimate_D}), by the H\"older inequalities we have (since $\frac{5}{6}=\frac{8}{15}+\frac{3}{10}$ and
$\frac{1}{\q_0}=\frac{1-2\alpha}{10}+\frac{3}{10}$ as $\q_0=\frac{5}{2-\alpha}$)
$$\|\mathds{1}_{Q_\mathbf{R}}\vu\theta\|_{\M_{t,x}^{\p_0,\q_0}}\leq C\|\mathds{1}_{Q_\mathbf{R}}\vu\theta\|_{\M_{t,x}^{\frac{6}{5},\q_0}}\leq C\|\mathds{1}_{Q_\mathbf{R}}\vu\|_{\M_{t,x}^{\frac{15}{8},\frac{10}{1-2\alpha}}}\|\mathds{1}_{Q_\mathbf{R}}\theta\|_{\M_{t,x}^{\frac{10}{3},\frac{10}{3}}},$$
since $\frac{15}{8}<3$ and $\frac{10}{1-2\alpha}<310$ (recall that $0<\alpha\ll1$, see Remark \ref{Rem_ValueAlphaSmall} below), we have
$$\|\mathds{1}_{Q_\mathbf{R}}\vu\theta\|_{\M_{t,x}^{\p_0,\q_0}}\leq C\|\mathds{1}_{Q_{\bar R_2}}\vu\|_{\M_{t,x}^{3,310}}\|\mathds{1}_{Q_{R_3}}\theta\|_{\M_{t,x}^{\frac{10}{3},\frac{10}{3}}}<+\infty,$$
since we have the controls (\ref{Recap_Controls}).\\

We deduce that $\|\vD\|_{\M_{t,x}^{\p_0,\q_0}}<+\infty$.
\item The term $\vE$ defined in (\ref{Definition_QuantitiesMorreyHolder}) is treated as follows. Recall that by the equation (\ref{Equation_Pressure_intro}) over the pressure we have the expression 
$$p = \displaystyle{\sum_{i,j=1}^3}\frac{\partial_i\partial_j}{(-\Delta)}\left(u_i u_j\right)-\frac{\partial_{x_3}}{(-\Delta)}\theta.$$  
We consider now a positive test function $\varphi$ such that
$$\varphi \equiv 1\;\; \text{over}\; \; Q_{r_a}(t_0,x_0)\; \mbox{and }\; supp(\varphi)\subset Q_{R_3}(t_0,x_0),$$
where $4\mathbf{R}<r_a<R_3$. Note in particular that by definition of the auxiliary functions $\phi$ and $\eta$ (see (\ref{Def_Eta})) we have the identities $\eta=\eta \varphi$ and $\eta\vn \varphi=0$. Thus for the term $\eta \vn p$ we have
$$\eta \vn p= \eta\varphi \vn p= \eta \vn (\varphi p)=\eta  \vn\left(\sum^3_{i,j= 1}\varphi\frac{\partial_i  \partial_j  }{(- \Delta )} (u_i u_j )-\varphi\frac{\partial_{x_3}}{(-\Delta)}\theta\right),$$ 
and this expression can be rewritten in the following manner
\begin{equation}\label{Formule_Pression_Ladyz}
\begin{split}
\eta \vn p=&\underbrace{\sum^3_{i,j= 1} \eta \frac{\vn \partial_i  \partial_j }{(- \Delta )} (\varphi u_i u_j )}_{(a)}-\underbrace{\sum^3_{i,j= 1}  \frac{\eta \vn\partial_i}{(- \Delta )} (\partial_j \varphi) u_i u_j}_{(b)}-\underbrace{\sum^3_{i,j= 1}  \frac{\eta \vn\partial_j}{(- \Delta )} (\partial_i \varphi) u_i u_j}_{(c)}\\
&+ 2 \underbrace{\sum^3_{i,j= 1} \eta \frac{\vn}{(-\Delta )}(\partial_i \partial_j\varphi) (u_i u_j)}_{(d)}+ \underbrace{\eta \frac{\vn \big(  (\Delta \varphi) p \big) }{(-\Delta)}}_{(e)}
-2 \underbrace{\sum^3_{i= 1}\eta\frac{\vn\big(\partial_i ( (\partial_i \varphi) p )\big)}{(-\Delta)}}_{(f)}\\
& -\underbrace{\eta\frac{\vn}{(-\Delta)}\partial_{x_3}\left(\varphi\theta \right)}_{(g)}+\underbrace{\eta\frac{\vn}{(-\Delta)}\left((\partial_{x_3}\varphi)\theta \right)}_{(h)}.
\end{split}
\end{equation}
We study each term above separately. 

\begin{itemize}
\item[$*$] For the term $(a)$ in (\ref{Formule_Pression_Ladyz}), since the Riesz transforms are bounded in Morrey spaces, we obtain
$$\left\|\eta\frac{\vn  \partial_i \partial_j }{(- \Delta )}( \varphi u_i u_j )\right\|_{\mathcal{M}^{\p_0,\q_0}_{t,x}}\leq C\left\|\vn( \varphi u_i u_j )\right\|_{\mathcal{M}^{\p_0,\q_0}_{t,x}},$$
now, for $1\leq k\leq 3$, using all the information available over $\vu$ (see (\ref{Recap_Controls})), by Lemma \ref{lemma_locindi} (recall that $\p_0\leq \frac{6}{5}<\frac{3}{2}$ and $\q_0=\frac{5}{2-\alpha}<155$) and by the H\"older inequality in Morrey spaces, we have
$$\left\|(\partial_k\varphi)u_i u_j \right\|_{\mathcal{M}^{\p_0,\q_0}_{t,x}}\leq C\left\|\mathds{1}_{Q_{\bar R_2}}u_i u_j \right\|_{\mathcal{M}^{\frac{3}{2}, 155}_{t,x}}\leq C\|\mathds{1}_{Q_{\bar R_2}}u_i\|_{\mathcal{M}^{3, 310}_{t,x}}\|\mathds{1}_{Q_{\bar R_2}}u_j\|_{\mathcal{M}^{3, 310}_{t,x}}<+\infty.$$
By essentially the same arguments (recall the informations over $\vu$ given in (\ref{Recap_Controls})) we have (note that $\p_0\leq \frac{6}{5}$ and $\q_0=\frac{5}{2-\alpha}<\frac{120}{47}$ since $0<\alpha\ll1$):
$$\|\varphi(\partial_k u_i) u_j \|_{\mathcal{M}^{\p_0, \q_0}_{t,x}}\leq C\|\varphi(\partial_k u_i) u_j \|_{\mathcal{M}^{\frac{6}{5}, \frac{120}{47}}_{t,x}}\leq C\|\mathds{1}_{Q_{R_1}}\vn \otimes \vu\|_{\mathcal{M}^{2, \frac{3720}{1445}}_{t,x}}\|\mathds{1}_{Q_{R_3}} u_j\|_{\mathcal{M}^{3, 310}_{t,x}},
$$
since $\frac{5}{6}=\frac{1}{2}+\frac{1}{3}$ and $\frac{47}{120}=\frac{1445}{3720}+\frac{1}{310}$. Recall now that we have $\tau_1=\frac{5\tau_0}{\tau_0+5}$ (see (\ref{Recap_Controls})) and since $\frac{5}{1-\alpha}<\tau_0<6$, the parameter $\tau_0$ can be chosen so that $\frac{3720}{1445}<\tau_1$ and we obtain
$$\|\varphi(\partial_k u_i) u_j \|_{\mathcal{M}^{\p_0, \q_0}_{t,x}}\leq C\|\mathds{1}_{Q_{R_1}}\vn \otimes \vu\|_{\mathcal{M}^{2,\tau_1}_{t,x}}\|\mathds{1}_{Q_{R_3}} u_j\|_{\mathcal{M}^{3, 310}_{t,x}}<+\infty,$$
and a symmetric argument gives
\begin{eqnarray*}
\|\varphi u_i(\partial_k u_j) \|_{\mathcal{M}^{\p_0, \q_0}_{t,x}} &\leq & C\|\mathds{1}_{Q_{R_3}}u_i\|_{\mathcal{M}^{3, 310}_{t,x}}\|\mathds{1}_{Q_{R_1}} \vn \otimes \vu\|_{\mathcal{M}^{2, \tau_1}_{t,x}}<+\infty.
\end{eqnarray*}
Thus we can deduce that we have the estimate 
$$\left\|\eta\frac{\vn  \partial_i \partial_j }{(- \Delta )}( \varphi u_i u_j )\right\|_{\mathcal{M}^{\p_0,\q_0}_{t,x}}<+\infty.$$
\begin{Remark}\label{Rem_ValueAlphaSmall}
Note that the condition $\q_1=\frac{5}{2-\alpha}<\frac{120}{47}$ is the most restrictive constraint over the parameter $\alpha$ and it implies that $0<\alpha<\frac{1}{24}\sim 0.04166$.
\end{Remark}

\item[$*$] The terms $(b)$ and $(c)$ of (\ref{Formule_Pression_Ladyz}) can be treated in a similar manner and using the information available in (\ref{Recap_Controls}) we have:
\begin{eqnarray*}
\left\|\frac{\eta \vn\partial_i}{(- \Delta )} (\partial_j \varphi) u_i u_j\right\|_{\mathcal{M}^{\p_0, \q_0}_{t,x}}&\leq &\left\|\frac{\eta \vn\partial_i}{(- \Delta )} (\partial_j \varphi) u_i u_j\right\|_{\mathcal{M}^{\frac{6}{5}, \frac{120}{47}}_{t,x}}\leq  C \left\|(\partial_j \varphi) u_i u_j\right\|_{\mathcal{M}^{\frac{6}{5}, \frac{120}{47}}_{t,x}}\\
&\leq &C\|\mathds{1}_{Q_{\bar R_2}}u_i u_j\|_{\mathcal{M}^{\frac{3}{2}, 155}_{t,x}}\leq C\|\mathds{1}_{Q_{\bar R_2}}u_i\|_{\mathcal{M}^{3, 310}_{t,x}}\|\mathds{1}_{Q_{\bar R_2}}u_j\|_{\mathcal{M}^{3, 310}_{t,x}}<+\infty.
\end{eqnarray*}
\item[$*$] The term $(d)$ is treated as follows. By Lemma \ref{lemma_locindi}, since $\p_0\leq \frac{6}{5}<\frac{3}{2}$ and $\q_0<\frac{120}{47}< \frac{15}{4}$, we have
$$\left\|\eta \frac{\vn}{(-\Delta )}(\partial_i \partial_j \varphi) (u_i u_j)\right\|_{\mathcal{M}^{\p_0, \q_0}_{t,x}}\leq C\left\|\eta \frac{\vn}{(-\Delta )}(\partial_i \partial_j \varphi) (u_i u_j)\right\|_{\mathcal{M}^{\frac{3}{2}, \frac{15}{4}}_{t,x}}.$$
Now, by the space inclusion $L^{\frac{3}{2}}_t L^{\infty}_{x}\subset \mathcal{M}^{\frac{3}{2}, \frac{15}{4}}_{t,x}$ we obtain
$$\left\|\eta \frac{\vn}{(-\Delta )}(\partial_i \partial_j \varphi) (u_i u_j)\right\|_{\mathcal{M}^{\frac{3}{2}, \frac{15}{4}}_{t,x}}\leq C\left\|\eta \frac{\vn}{(-\Delta )}(\partial_i \partial_j  \varphi) (u_i u_j)\right\|_{L^{\frac{3}{2}}_t L^{\infty}_{x}}.$$
Following the same ideas displayed in formulas (\ref{Formula_intermediairevVV80})-(\ref{Formula_intermediairevVV81}), due to the support properties of the auxiliary functions we obtain
$$\left\|\eta\frac{\vn}{(-\Delta )}(\partial_i\partial_j \varphi) (u_i u_j)\right\|_{L^{\frac{3}{2}}_tL^{\infty}_{x}}\leq C\|\mathds{1}_{Q_{\bar R_2}}u_i u_j\|_{L^{\frac{3}{2}}_{t,x}}\leq C\|\mathds{1}_{Q_{\bar R_2}}\vu\|_{\mathcal{M}^{3,310}_{t,x}}\|\mathds{1}_{Q_{\bar R_2}}\vu\|_{\mathcal{M}^{3,310}_{t,x}}<+\infty.$$

\item[$*$] The term $(e)$ of (\ref{Formule_Pression_Ladyz}) follows the same ideas as before, indeed we have
$$\left\|\eta \frac{\vn \big(  (\Delta \varphi) p \big) }{(-\Delta)}\right\|_{\mathcal{M}^{\p_0, \q_0}_{t,x}}\leq C\left\|\eta \frac{\vn \big(  (\Delta \varphi) p \big) }{(-\Delta)}\right\|_{\mathcal{M}^{\frac{3}{2}, \frac{15}{4}}_{t,x}}\leq C\left\|\eta \frac{\vn \big(  (\Delta \varphi) p \big) }{(-\Delta)}\right\|_{L^{\frac{3}{2}}_t L^{\infty}_{x}}\leq C \|\mathds{1}_{Q_{R}}p\|_{L^{\frac{3}{2}}_{t,x}}<+\infty,$$
since we have by hypothesis that $\mathds{1}_{Q_{R}}p\in L^{\frac{3}{2}}_{t,x}(\mathbb{R}\times \R)$.

\item[$*$] The term (f) of (\ref{Formule_Pression_Ladyz}) is estimated in a very similar manner:
\begin{eqnarray*}
\left\|\eta\frac{\vn\big(\partial_i ( (\partial_i \varphi) p )\big)}{(-\Delta)}\right\|_{\mathcal{M}^{\p_0, \q_0}_{t,x}}&\leq& C\left\|\eta\frac{\vn\big(\partial_i ( (\partial_i \varphi) p )\big)}{(-\Delta)}\right\|_{\mathcal{M}^{\frac{6}{5}, \frac{120}{47}}_{t,x}}\leq C\left\|\eta\frac{\vn\big(\partial_i ( (\partial_i \varphi) p )\big)}{(-\Delta)}\right\|_{\mathcal{M}^{\frac{3}{2}, \frac{15}{4}}_{t,x}}\\
&\leq &C\left\|\eta\frac{\vn\big(\partial_i ( (\partial_i \varphi) p )\big)}{(-\Delta)}\right\|_{L^{\frac{3}{2}}_t L^{\infty}_{x}}\leq C \|\mathds{1}_{Q_{R}}p\|_{L^{\frac{3}{2}}_{t,x}}<+\infty.
\end{eqnarray*}

\item[$*$] For the quantity (g) of (\ref{Formule_Pression_Ladyz}), since $\p_0\leq \frac{6}{5}<\frac{10}{3}$ and $\q_0<\frac{120}{47}<\frac{10}{3}$, we have (since the Riesz transforms are bounded in Lebesgue spaces)
$$\left\|\eta\frac{\vn\partial_{x_3}}{(-\Delta)}\left(\varphi\theta \right)\right\|_{\mathcal{M}^{\p_0, \q_0}_{t,x}}\leq C\left\|\frac{\vn\partial_{x_3}}{(-\Delta)}\left(\varphi\theta \right)\right\|_{L^{\frac{10}{3}}_{t,x}}\leq C\left\|\varphi\theta \right\|_{L^{\frac{10}{3}}_{t,x}}<+\infty.
$$
\item[$*$] Finally, the term (h) of (\ref{Formule_Pression_Ladyz}), is treated as follows
\begin{eqnarray*}
\left\|\eta\frac{\vn}{(-\Delta)}\left((\partial_{x_3}\varphi)\theta \right)\right\|_{\mathcal{M}^{\p_0, \q_0}_{t,x}}&\leq &C\left\|\eta\frac{\vn}{(-\Delta)}\left((\partial_{x_3}\varphi)\theta \right)\right\|_{\mathcal{M}^{\frac{3}{2}, \frac{15}{4}}_{t,x}}\\
&\leq& C\left\|\eta\frac{\vn}{(-\Delta)}\left((\partial_{x_3}\varphi)\theta \right)\right\|_{L^{\frac{3}{2}}_t L^{\infty}_{x}}\leq C \|\mathds{1}_{Q_{R}}\theta\|_{L^{\frac{10}{3}}_{t,x}}<+\infty.
\end{eqnarray*}
Gathering all these estimates, we finally obtain that $\|\vE\|_{\M_{t,x}^{\p_0,\q_0}}<+\infty$.
\end{itemize}
\item The term $\vF$ given in (\ref{Definition_QuantitiesMorreyHolder}) is treated as follows:
$$\|\vF\|_{\M_{t,x}^{\p_0,\q_0}} =\left\|\left[\begin{matrix}
\eta \theta e_3\\
0
\end{matrix}\right]\right\|_{\M_{t,x}^{\p_0,\q_0}} \leq \|\eta \theta\|_{\M_{t,x}^{\p_0,\q_0}}\leq \|\mathds{1}_{Q_{R_3}}\theta\|_{\M_{t,x}^{\frac{12}{5},\frac{50}{9}}}<+\infty,$$
where we used the Lemma \ref{lemma_locindi} and $\p_0\leq \frac{6}{5}<\frac{12}{5}$, $\q_0=\frac{5}{2-\alpha}<\frac{50}{9}$ as well as the controls (\ref{Recap_Controls}). We finally obtain that $\|\vF\|_{\M_{t,x}^{\p_0,\q_0}}<+\infty$.\\
\end{itemize}
With all the previous computations we have proven all the information stated in (\ref{Objetivo_MorreyLadyz}), which applied in the integral representation formula (\ref{Duhamel_FullBoussinesq}) allows us, with Lemma \ref{Lemma_parabolicHolder}, to conclude that  $\vV\in \dot{\mathcal{C}}^\alpha(\mathbb{R}\times \R)$ with $0<\alpha\ll1$, and since we have $$\vV=\displaystyle{\left[\begin{matrix}
\eta\vu\\[1mm]
\eta\theta
\end{matrix}\right]},$$
we deduce that $\vu$ and $\theta$ are also H\"older regular over a small neighborhood of the point $(t_0, x_0)$ and this finishes the proof of Theorem \ref{Teo_HolderRegularity}. \hfill $\blacksquare$\\

\noindent{\bf Acknowledgement:} C. M.  thanks the support of the ERC-CZ Grant LL2105 CONTACT of the Faculty of Mathematics and Physics of Charles University. Also , C.M.  acknowledges for the support of the project GAUK No.  45612, of Charles University.



\end{document}